\documentclass[sn-mathphys,Numbered]{sn-jnl}

\usepackage{graphicx}%
\usepackage{multirow}%
\usepackage{amsmath,amssymb,amsfonts,stmaryrd}%
\usepackage{amsthm}%
\usepackage{mathrsfs}%
\usepackage[title]{appendix}%
\usepackage{xcolor}%
\usepackage{textcomp}%
\usepackage{manyfoot}%
\usepackage{booktabs}%
\usepackage{algorithm}%
\usepackage{algorithmicx}%
\usepackage{algpseudocode}%
\usepackage{listings}%
\usepackage{hyperref,bm}
\newcommand\norm[1]{\left\lVert#1\right\rVert}

\usepackage{subcaption}

\newtheorem{remark}{Remark}%

\begin{document}

\title[Fast parametric analysis of Kirchhoff-Love shells]{Fast parametric analysis of trimmed multi-patch isogeometric Kirchhoff-Love shells using a local reduced basis method}

\author*[1]{\fnm{Margarita} \sur{Chasapi}}\email{margarita.chasapi@epfl.ch}

\author[1]{\fnm{Pablo} \sur{Antolin}}

\author[1,2]{\fnm{Annalisa} \sur{Buffa}}

\affil*[1]{\orgdiv{Institute of Mathematics}, \orgname{\'Ecole Polytechnique F\'ed\'erale de Lausanne},
\city{Lausanne},  \country{Switzerland}}

\affil[2]{\orgdiv{Instituto di Matematica Applicata e Tecnologie Informatiche 'E. Magenes' (CNR)}, \city{Pavia}, \country{Italy}}


\abstract{This contribution presents a model order reduction framework for real-time efficient solution of trimmed, multi-patch isogeometric Kirchhoff-Love shells. In several scenarios, such as design and shape optimization, multiple simulations need to be performed for a given set of physical or geometrical parameters. This step can be computationally expensive in particular for real world, practical applications. We are interested in geometrical parameters and take advantage of the flexibility of splines in representing complex geometries. In this case, the operators are geometry-dependent and generally depend on the parameters in a non-affine way. Moreover, the solutions obtained from trimmed domains may vary highly with respect to different values of the parameters. Therefore, we employ a local reduced basis method based on clustering techniques and the Discrete Empirical Interpolation Method to construct affine approximations and efficient reduced order models. In addition, we discuss the application of the reduction strategy to parametric shape optimization. Finally, we demonstrate the performance of the proposed framework to parameterized Kirchhoff-Love shells through benchmark tests on trimmed, multi-patch meshes including a complex geometry. The proposed approach is accurate and achieves a significant reduction of the online computational cost in comparison to the standard reduced basis method.}


\maketitle

\section{Introduction}\label{sec1}
In the past decades, the integration of geometric design and finite element analysis has attracted a lot of attention within the computational engineering community. The introduction of Isogeometric Analysis (IGA) \cite{Hughes2005} paved the way for an integrated framework from Computer Aided Design (CAD) to numerical simulation followed by a wide range of successful applications in several fields. Adopting the isogeometric paradigm, the effort involved in meshing and geometry clean-up can be circumvented by employing the same representation for both geometric description and numerical analysis. The prevailing technology for geometric representation in CAD is based on splines, namely B-splines and in particular non-uniform rational B-splines (NURBS). The reader is further referred to  \cite{Hughes2005,Cottrell2007} for a detailed overview of the method.

In the context of shell analysis, IGA offers clear advantages. In particular, Kirchhoff-Love shell formulations are based on fourth order partial differential equations (PDEs) and require $C^1$-continuity that hinders the use of classical $C^0$-continuous finite elements. Due to the higher continuity of splines, isogeometric discretizations are well suited for the modeling of higher-order PDEs. An overview of isogeometric methods for Kirchhoff-Love shell analysis is given in \cite{Kiendl2009}.

Nevertheless, real-world applications of complex geometries require the treatment of further issues. Typically, complex shapes are represented in CAD using Boolean operations such as intersection and union. This results in trimmed meshes, where the underlying discretization is unfitted with respect to the physical object. First, tackling trimmed surfaces needs to be properly addressed in the numerical simulation. A detailed review on trimming and related challenges such as numerical integration, conditioning, and others, is given in \cite{Marussig2018}. We also refer to \cite{Breitenberger.2015, Guo2018} in the context of trimmed shell analysis. Besides trimming, complex shapes are commonly represented by multiple, non-conforming patches in CAD and require suitable coupling strategies. In particular for shells, $C^1$-continuity between adjacent patches is required for the analysis. There are several methods in the literature to achieve displacement and rotational continuity in a weak sense, such as the penalty, Nitsche's, and mortar methods. A comprehensive overview of these methods in the context of IGA is given in \cite{Apostolatos2014}. We also refer to previous works on mortar methods \cite{brivadis2015,Horger2019,Dittmann2019,Chasapi2020}, penalty approaches \cite{Herrema2018,Coradello2021,Prosperio2022}, and Nitsche's method \cite{Nguyen-Thanh2017,Guo2015,Guo2017,Benzaken2021} including their application to the analysis of trimmed shells.
For the case of conforming patches, $C^1$-continuity across interfaces can be even imposed strongly, as shown recently in~\cite{farahat2023isogeometric}.

On the other hand, there exist several applications where the analysis needs to be performed rapidly for many different parametric configurations. This is, for example, the case in design and shape optimization, where the geometry is updated on the fly within the optimization loop. The reader is further referred to \cite{Nagy2013,Kiendl2014,Hirschler2021}, for an overview on isogeometric shape optimization. To achieve a computational speedup, research efforts have been devoted to the development of efficient reduced order models (ROMs) for the solution of parameterized problems. In particular, we refer to previous works in the context of combining IGA with ROMs \cite{Manzoni2015,Salmoiraghi2016,Devaud2017}. Most of these works employ the reduced basis method to construct projection-based ROMs in combination with hyper-reduction techniques such as the Empirical Interpolation Method \cite{Barrault2004} to handle geometrical parameterizations for non-affine problems. A detailed overview on reduced basis methods is given in  \cite{Hesthaven2016,QMN_RBspringer}. We also refer to \cite{Maquart2020,Chasapi2022} in the context of isogeometric ROMs for multi-patch geometries. Nevertheless, the study of efficient ROMs in the context of isogeometric shell analysis is a subject that still lacks thorough investigation. The application of reduced basis methods to isogeometric Kirchhoff-Love shells has been explored in \cite{Rinaldi2015} for simple geometries. However, a general ROM framework suitable for trimmed, multi-patch geometries of industrial relevance is still missing.

This contribution aims to investigate the use of efficient ROM strategies for fast isogeometric Kirchhoff-Love shell analysis formulated on parameterized trimmed and multi-patch geometries. Based on our previous work \cite{Chasapi2023}, we employ a local reduced basis method to construct efficient ROMs for problems formulated on parameterized unfitted geometries obtained by trimming operations. Local ROMs have been a subject of previous research works within the model reduction community \cite{Haasdonk2011,Amsallem2012,Peherstorfer2014,Pagani2018}. The method we use in this work is based on extension and parameter-based clustering of snapshots to apply reduction techniques such as the Proper Orthogonal Decomposition and Discrete Empirical Interpolation Method \cite{Chaturantabut2010,Negri2015}.

We recast the Kirchhoff-Love shell formulation into parameterized trimmed geometries within a multi-patch setting. The proposed reduction framework enables an efficient offline/online decomposition and is agnostic to the underlying coupling method, therefore different strategies can be in principle considered to enforce displacement and rotational continuity in a weak sense. In addition, we discuss the application of the reduction framework to parametric shape optimization problems. We also refer to previous works on ROM-based parametric optimization \cite{Amsallem2015,Manzoni2012,Antil2012}.

The manuscript is structured as follows: Section 2 provides a brief overview of basic concepts related to B-splines and trimming formulated on parameterized domains. In Section 3 we present the parameterized formulation for the Kirchhoff-Love shell, while Section 4 provides the necessary definitions related to the coupling of multi-patches.  We discuss the local reduced basis method in Section 5 and its application to parametric shape optimization in Section 6. In Section 7 we study several numerical experiments including non-conforming discretizations and a complex geometry to assess the performance of the proposed framework. The main conclusions that can be drawn from this study are finally summarized in Section 8.

\section{Parameterized trimmed geometries in isogeometric analysis}\label{sec2}
In this section, we provide a brief review of some basic concepts related to B-splines and trimming in isogeometric analysis. In particular, we formulate the trimming operation in the context of parameterized geometries for single-patch domains.
For a detailed review of splines in the context of isogeometric analysis, the reader is referred to \cite{Hughes2005,Cottrell2007,Piegl1995}. Moreover, we refer to \cite{Marussig2018} for a comprehensive review of trimming in isogeometric analysis.

\subsection{B-spline basis functions}
To illustrate the basic concept of B-splines, we introduce a \emph{knot vector} in the parametric space $[0,1]$ as a non-decreasing sequence of real values denoted by $\Xi = \{\xi_1,\dots,\xi_{n+p+1}\}$. Here, the integer $n$ denotes the number of basis functions and $p$ the degree. We further introduce a univariate B-spline basis function $b^j_{i_j,p_j}$ where $p_j$ denotes the degree and $i_j$ the index of the function in the $j$-th parametric direction. The B-spline functions $\mathcal{B}_{\boldsymbol{i},\boldsymbol{p}}(\boldsymbol{\xi})$ can be easily defined in multiple dimensions exploiting the tensor product of univariate B-splines:
\begin{equation}\label{eq11}
	\mathcal{B}_{\bf{i},\bf{p}}(\boldsymbol{\xi}) = \prod_{j=1}^{{\hat{d}}} b_{i_j,p_j}^j(\xi^j),
\end{equation}
where $\hat{d}$ is the dimension of the parametric space. Moreover, the vector ${\bf{i}}=(i_1,\dots,i_{\hat{d}})$ is a multi-index denoting the position in the tensor-product structure and ${\bf{p}}=(p_1,\dots,p_{\hat{d}})$ the polynomial degree corresponding to the parametric coordinate $\boldsymbol{\xi}=(\xi^1,\dots,\xi^{\hat{d}})$. For simplicity, we will assume from now on that the vector $\bf{p}$ is identical in all parametric directions and can be replaced by a scalar value $p$. The B-spline basis is $C^{p-k}$-continuous at every knot, where $k$ is the multiplicity of the knot, and is $C^{\infty}$ elsewhere. 
The concept of B-splines can be easily extended to rational B-splines and the interested reader is further referred to \cite{Cottrell2007} for a detailed exposition. 

\subsection{Parameterized trimming} \label{sec:paramtrimming}
In the following, we briefly present the mathematical foundation of trimming in isogeometric analysis and recast its formulation into the context of parameterized geometries. For this purpose, we adopt the notation previously introduced in \cite{Antolin2019,Coradello2020}. 

Let us first define the parameterized physical domain $\Omega(\bm{\mu}) \subset \mathbb{R}^d$ described by the geometrical parameters $\bm{\mu} \in \mathcal{P} \subset \mathbb{R}^M$, where $d$ is the dimension of the physical space for our problem, $\mathcal{P}$ is the space of parameters, and $M$ is the number of parameters. We will show later how to obtain $\Omega(\bm{\mu})$ through trimming operations. Now we consider the non-trimmed, single-patch physical domain $\Omega_0 \subset \mathbb{R}^d$ and its counterpart in the parameter space $\hat{\Omega}_0 = [0,1]^{\hat{d}}$. In the following, we assume that the non-trimmed domain $\Omega_0$ is parameter-independent without loss of generality. In principle, the extension to the parameter dependent case is straightforward. The multi-patch setting will be discussed below in Section~\ref{sec4}.

\begin{figure}[!h]
	\centering
	\includegraphics[width=1.0\textwidth]{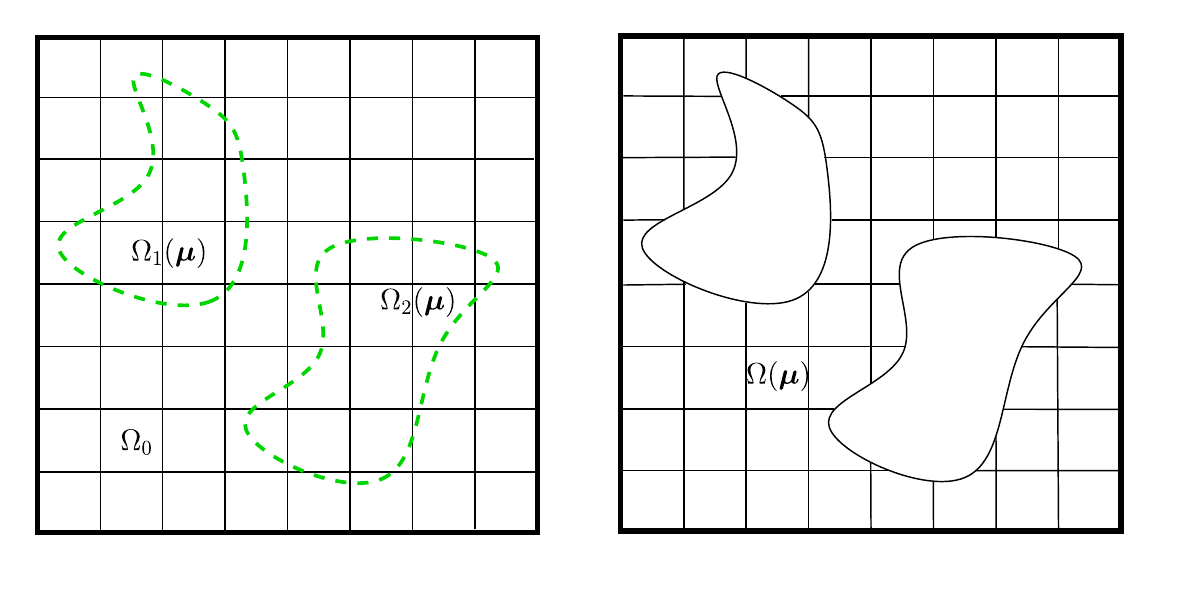} \\
	\caption{Exemplary trimmed domain. The final rectangular domain $\Omega({\bm{\mu}})$ is a result of trimming away the green regions $\Omega_1({\bm{\mu}})$ and $\Omega_2({\bm{\mu}})$ from the non-trimmed domain $\Omega_0$.}\label{fig:trimming_example}
\end{figure}

Given a control point mesh $\bf{P_i}$, a spline geometric map ${\bf{F}} : \hat{\Omega}_0 \to \Omega_0$ can be defined as
\begin{equation}\label{eq12}
	{\bf{F}}(\boldsymbol{\xi}) = \sum_{\bf{i}} \mathcal{B}_{{\bf{i}},p}(\boldsymbol{\xi})\bf{P_i}.
\end{equation}
As a result, by an abuse of notation, the non-trimmed physical domain is obtained by $\Omega_0 = {\bf{F}}(\hat{\Omega}_0)$. 
In what follows we are interested in parameterizing the trimmed regions. For this purpose, we introduce $K$ parameter-dependent trimmed regions $\Omega_1(\bm{\mu}),\dots,\Omega_K(\bm{\mu}) \subset \mathbb{R}^d$ that are cut away by the trimming operation. The obtained physical domain can then be expressed as:
\begin{equation}\label{eq13}
	\Omega({\bm{\mu}}) = \Omega_0 \backslash \bigcup_{i=1}^K {\overline{\Omega_i}}(\bm{\mu}).
\end{equation}
The concept is illustrated in Fig.~\ref{fig:trimming_example}. Following this definition, the boundary of the domain consists of a trimmed part $\partial{\Omega({\bm{\mu}})}\backslash\partial{\Omega_0}$ and a part that corresponds to the original domain $\partial{\Omega({\bm{\mu}})}\cap\partial{\Omega_0}$. For the sake of simplicity, we will assume from now on that Dirichlet boundary conditions are not applied on the trimming boundary. The interested reader may refer to \cite{Elfverson2019,Buffa2020} regarding the weak imposition of Dirichlet boundary conditions on the trimming boundary in combination with stabilization techniques. 

Furthermore, it should be noted that the elements and basis functions are defined on the non-trimmed domain $\Omega_0$ and the original domain remains unchanged by the trimming operation. Thus, we define the B-spline space $S^p_{0}$, which is constructed upon $\Omega_0$ and is independent of the parameters ${\bm{\mu}}$, as
\begin{equation}\label{eq14a}
	S^p_{0} = \text{span}\{\mathcal{B}_{{\bf{i}},p} \circ {\bf{F}}^{-1}\}.
\end{equation}
In fact, let us now rewrite the trimming operation in the parametric domain $\hat{\Omega}_0$. The counterpart of the parameter-dependent trimmed regions $\Omega_i(\bm{\mu}), i=1,\dots,K$ (see also~\eqref{eq13}) in the parametric domain are defined as $\hat{\Omega}_i \subset \hat{\Omega}_0$, while it holds that $\Omega_i(\bm{\mu}) = {\bf{F}}(\hat{\Omega}_i), i=1,\dots,K$. Thus, Equation~\eqref{eq13} can be reinterpreted as
\begin{equation}\label{eq13b}
	\hat{\Omega}({\bm{\mu}}) = \hat{\Omega}_0 \backslash \bigcup_{i=1}^K {\overline{\hat{\Omega}_i}}(\bm{\mu}).
\end{equation}
Let us now introduce the B-spline space  of degree $p$ restricted to the trimmed domain ${\Omega}({\bm{\mu}})$:
\begin{equation}\label{eq14}
	S^p({{{{\Omega}}({\bm{\mu}})}}) = \text{span}\{\mathcal{B}_{{\bf{i}},p} \circ {\bf{F}}^{-1}, \text{supp}({\mathcal{B}_{{\bf{i}},p}}) \cap  \hat{\Omega}({\bm{\mu}})  \ne \varnothing \},
\end{equation}
where $\text{supp}({\mathcal{B}_{{\bf{i}},p}})$ is the support of the non-trimmed basis functions. 
We remark that for numerical integration, a re-parameterization of cut elements is performed in the trimmed parametric domain $\hat{\Omega}({\bm{\mu}})$. In this work we follow the high-order reparameterization procedure discussed in \cite{Wei2021}. The trimming operation can be understood as restricting the map ${\bf{F}}$ to an \emph{active} region of the original domain, which is the visible part after the trimmed regions are cut away. Indeed, only the basis functions whose support intersects $\hat{\Omega}({\bm{\mu}}) $ are active. The remaining basis functions are inactive and do not contribute to the solution discretization of the problem. Moreover, the active basis functions and the dimension of the space may change for different values of the parameters ${\bm{\mu}}$. Therefore, to apply reduce order modeling techniques we will rely on the definition of the non-trimmed domain $\Omega_0$ and associated B-spline space $S^p_{0}$ following \cite{Chasapi2023}. This aspect will be further discussed  in Section~\ref{sec5}.

\section{Parameterized Kirchhoff-Love shell formulation}\label{sec3}
In the following we introduce the weak formulation of the Kirchhoff-Love shell problem following closely the notation set in \cite{Benzaken2021,Coradello2021}. Let us consider a parameterized single-patch computational domain $\Omega({\bm{\mu}}) \subset \mathbb{R}^3$ representing a two-dimensional manifold with smooth boundary $\Gamma({\bm{\mu}}) = \partial \Omega({\bm{\mu}})$. The boundary is partitioned such that both displacements and rotations, respectively, as well as their energetically conjugate shears forces and bending moments can be prescribed on the boundary. Thus, we partition  $\Gamma({\bm{\mu}})$ into a Dirichlet boundary $\Gamma_D({\bm{\mu}}) = \Gamma_{D,u}({\bm{\mu}}) \cup \Gamma_{D,\theta}({\bm{\mu}})$ associated with prescribed transverse displacements and normal rotations as well as a Neumann boundary $\Gamma_N({\bm{\mu}}) = \Gamma_{N,s}({\bm{\mu}}) \cup \Gamma_{N,b}({\bm{\mu}})$ associated with applied transverse shear forces and bending moments, respectively. Note that it holds $\Gamma({\bm{\mu}}) = \overline{\Gamma_D({\bm{\mu}}) \cup \Gamma_N({\bm{\mu}})}$, $ \Gamma_{D,u}({\bm{\mu}}) \cap \Gamma_{N,s}({\bm{\mu}})= \varnothing$ and $\Gamma_{D,\theta}({\bm{\mu}}) \cap \Gamma_{N,b}({\bm{\mu}}) = \varnothing$. Now let us further introduce a set of corners $\chi({\bm{\mu}})\subset \Gamma({\bm{\mu}})$ that can be  decomposed into a Dirichlet part $\chi_D ({\bm{\mu}})= \chi ({\bm{\mu}})\cap \Gamma_D({\bm{\mu}})$ and a Neumann part $\chi_N({\bm{\mu}}) = \chi ({\bm{\mu}})\cap \Gamma_N({\bm{\mu}})$. Note that $\chi({\bm{\mu}})= \chi_D({\bm{\mu}}) \cup  \chi_N({\bm{\mu}})$ and $\chi_D ({\bm{\mu}})\cap  \chi_N ({\bm{\mu}})= \varnothing$.
Let us also assume an applied body load $\tilde{\bf{f}} \in [L^2(\Omega({\bm{\mu}}))]^3$,  a prescribed bending moment $\tilde{{B}}_{nn} \in L^2(\Gamma_{N,b}({\bm{\mu}}))$, an applied twisting moment $\tilde{S}\big\vert_C \in \mathbb{R}$ at corner $C \in \chi_N$, and a prescribed transverse shear or ersatz traction $\tilde{\bf{T}} \in L^2(\Gamma_{N,s}({\bm{\mu}}))^3$ as defined in \cite{Benzaken2021}.

Hereinafter, for the sake of simplicity of exposition, and without loss of generality, we consider that no normal rotations are prescribed on the boundary, i.e., $\Gamma_{D,\theta}({\bm{\mu}})=\emptyset$, and that only homogeneous Dirichlet boundary conditions are applied on $\Gamma_{D}({\bm{\mu}})=\Gamma_{D,u}({\bm{\mu}})$.

We further recall the spline geometric map ${\bf{F}}$ in~\eqref{eq12} and construct a covariant basis, where the basis vectors are defined as
\begin{equation}\label{eq15}
{\bf{a}}_{\alpha}(\boldsymbol{\xi}) = {{\bf{F}}}_{,\alpha}(\boldsymbol{\xi}), \quad \alpha=1, 2.
\end{equation}
Here $ {{\bf{F}}}_{,\alpha}$ denotes the partial derivative of the spline geometric mapping with respect to the $\alpha$-th curvilinear coordinate. The midsurface normal vector ${\bf{a}}_{3}$ can be constructed as a normalized cross-product of the two in-plane vectors ${\bf{a}}_{\alpha}$, that is 
\begin{equation}\label{eq16}
	{\bf{a}}_{3} = \frac{	{\bf{a}}_{1} \times 	{\bf{a}}_{2}}{\norm{{\bf{a}}_{1} \times{\bf{a}}_{2}} }.
\end{equation}
Then, we can construct the contravariant basis vectors that satisfy the Kronecker relationship as $\bf{a} _{\alpha}\cdot \bf{a}^{\beta} = \delta_{\alpha}^{\beta}$. Note that it holds ${\bf{a}}_3 = {\bf{a}}^3$.  The reader is further referred to \cite{Benzaken2021} for an elaborate discussion on fundamentals of differential geometry that are relevant to the Kirchhoff-Love shell formulation. 

Let us now introduce a discrete space $V_h({\bm{\mu}})$ to approximate our problem:
\begin{equation}\label{eq17a}
V_h({\bm{\mu}})  = \{ \boldsymbol{v}_h \in S^p({{{{\Omega}}({\bm{\mu}})}})^3  : \left.\boldsymbol{v}_h\right\rvert_{\Gamma_D({\bm{\mu}})} = \bm{0}\},
\end{equation}
where the spline space $S^p({{{{\Omega}}({\bm{\mu}})}})^3$ has at least $C^1$-continuity, as typically required by isogeometric Kirchoff-Love formulations (see, e.g., \cite{Kiendl2009}). Higher continuity requirements will be discussed in Section~\ref{sec4b} for the case of Nitsche's interface coupling for multi-patch domains.

Now let us define the discrete weak formulation of the parameterized Kirchhoff-Love shell problem as: find ${\boldsymbol{u}}_h({\bm{\mu}}) \in V_h({\bm{\mu}})$ such that 
\begin{equation}\label{eq17}
	a({\boldsymbol{u}}_h,{\boldsymbol{v}}_h;\bm{\mu}) = f({\boldsymbol{v}}_h;\bm{\mu}), \qquad \forall {\boldsymbol{v}}_h \in V_h({\bm{\mu}})
\end{equation}
The parameterized bilinear form $a(\cdot,\cdot;\bm{\mu})$ is given as:
\begin{equation}\label{eq17b}
	a({\boldsymbol{u}}_h,{\boldsymbol{v}}_h;\bm{\mu}) = \int_{{\Omega}(\boldsymbol{\mu})} A(\boldsymbol{u}_h) : \alpha(\boldsymbol{v}_h) \, \textrm{d}\Omega+ \int_{{\Omega}(\boldsymbol{\mu})}  B(\boldsymbol{u}_h): \beta(\boldsymbol{v}_h) \, \textrm{d}\Omega,
\end{equation}
and the parameterized linear functional $f(\cdot;\bm{\mu})$ reads:
\begin{equation}\label{eq17c}
\begin{aligned}
	f({\boldsymbol{v}}_h;\bm{\mu}) &= \int_{{\Omega}(\boldsymbol{\mu})} {\tilde{\bf{f}}} \cdot \boldsymbol{v}_h \, \textrm{d}\Omega+ \int_{{\Gamma_{N,s}}(\boldsymbol{\mu})}  {\tilde{\bf{T}}} \cdot \boldsymbol{v}_h  \, \textrm{d}\Gamma \\&
 +  \int_{{\Gamma_{N,b}}(\boldsymbol{\mu})}  {\tilde{{B}}_{nn}} \theta_n (\boldsymbol{v}_h)  \, \textrm{d}\Gamma + \sum_{C \in \chi_N} (\tilde{S}v_3)\big\vert_C.
 \end{aligned}
\end{equation}
Note that $\theta_n(\boldsymbol{v}_h) = - \boldsymbol{n}\cdot\nabla{\boldsymbol{u}}^{\top}_h{\bf{a}}_3$ is the normal rotation, $ \boldsymbol{n} $ is the unit outward normal vector to the boundary $\Gamma( \boldsymbol{\mu})$, and the membrane and bending strain tensors are defined, respectively, as:
\begin{equation}\label{eq18}
\alpha(\boldsymbol{v}_h) = P\,  \text{sym}(\nabla^*(\boldsymbol{v}_h))\, P,\quad \beta(\boldsymbol{v}_h) = - P\,  \text{sym}({\bf{a}}_3\nabla^*\nabla^*(\boldsymbol{v}_h))\, P,
\end{equation}
where $\text{sym}(\cdot)$ denotes the symmetric part of a tensor, $\nabla^*$ is the surface gradient, and  $P= {\bf{I}} - {\bf{a}}_3 \otimes {\bf{a}}_3$ is the in-plane projector, with $\bf{I}$ the identity tensor.

Let us recall that in Kirchhoff-Love shell kinematics we assume that the transverse shear strains vanish, i.e., the normal vectors remain straight and normal during deformation. The rotations are therefore constrained as follows:
\begin{equation}\label{eq18a}
{\theta}(\boldsymbol{u}_h) = - \nabla{\boldsymbol{u}}^{\top}_h {\bf{a}}_3.
\end{equation}

Now, let us define the energetically conjugate stresses for both the membrane and bending strains. For this purpose, we assume a linear elastic constitutive model and define the following fourth-order elasticity tensor:
\begin{equation}\label{eq19}
	\begin{aligned}
 \mathbb{C} = &\mathbb{C}^{\alpha\beta\lambda\mu} {\bf{a}}_{\alpha}\otimes{\bf{a}}_{\beta}\otimes{\bf{a}}_{\lambda}\otimes{\bf{a}}_{\mu},\\
 &\text{with} \ \mathbb{C}^{\alpha\beta\lambda\mu} = \frac{E}{2(1+\nu)}(a^{\alpha\lambda}a^{\beta\mu}+a^{\alpha\mu}a^{\beta\lambda}+\frac{2\nu}{1-\nu}a^{\alpha\beta}a^{\lambda\mu}),
		\end{aligned}
\end{equation}
where repeated indices imply summation from 1 to 3, and $E$ and $\nu$ are the elasticity modulus and Poisson's ratio, respectively. Assuming a constant thickness $t$ and performing through-thickness integration we obtain the membrane and bending stresses as:
\begin{equation}\label{eq20}
A(\boldsymbol{u}_h) = t\,\mathbb{C} : \alpha(\boldsymbol{u}_h), \quad B(\boldsymbol{u}_h) = \frac{t^3}{12}\mathbb{C}:\beta(\boldsymbol{u}_h).
\end{equation}
For a rigorous derivation of the weak formulation the reader is further referred to \cite{Benzaken2021}.

\section{Multi-patch geometries}\label{sec4}
Let us now consider that the parameterized domain $\Omega(\bm{\mu})$ is split into non-overlapping subdomains, i.e., patches, such that
\begin{equation}\label{eq22}
	\overline{\Omega}(\bm{\mu}) = \bigcup_{k=1}^{N_p} {\overline{\Omega}}_{k}(\bm{\mu}_{k}),
\end{equation}
where $\Omega_{k}(\bm{\mu}_{k}) \cap \Omega_{l}(\bm{\mu}_{l}) = \varnothing$ for $l \ne k$. Here the patches ${\Omega}_{k}(\bm{\mu}_{k})$ are in principle trimmed and their definition follows the setting introduced in Section~\ref{sec:paramtrimming}. For ease of exposition, we assume that the parameters associated to the $k$-th patch coincide with the parameters describing the global computational domain such that $\bm{\mu}_{k} = \bm{\mu}$ for $k=1,\dots,N_p$, although in principle different choices are possible (see \cite{Chasapi2022} for more details).

Furthermore, we define the common interface $\gamma_{j}({\bm{\mu}})$, between two adjacent patches such that 
\begin{equation}\label{eq27}
 \gamma_{j}({\bm{\mu}})=\partial{\Omega_{k}}({\bm{\mu}})\cap\partial{\Omega_{l}}({\bm{\mu}}) \ne \varnothing, \quad k \ne l,\ \text{for } j=1,\dots,N_{\Gamma}
\end{equation}
where $N_{\Gamma}$ is the number of interfaces. We note that the interface can be both trimmed or non-trimmed.
The B-spline space~\eqref{eq14a} can be now rewritten  as 
\begin{equation}\label{eq24}
	S^p_{k,0} = \text{span}\{\mathcal{B}_{{\bf{i}},p}^{k} \circ {\bf{F}}_k^{-1} \},
\end{equation}
where $\mathcal{B}_{{\bf{i}},p}^{k}$ and  ${\bf{F}}_k$ are the B-spline functions and the geometric map associated to the $k$-th patch, respectively. The dimension of the associated space is denoted as $\mathcal{N}_{k,0}=\text{dim}(S^p_{k,0})$.
And, analogously to~\eqref{eq14}, the B-spline space of degree $p$ corresponding to the trimmed patch $\Omega_{k}(\bm{\mu})$ now reads
\begin{equation}\label{eq25}
	S^p({{{{\Omega}_{k}}({\bm{\mu}})}}) = \text{span}\{\mathcal{B}_{{\bf{i}},p}^{k} \circ {\bf{F}}_k^{-1}, \text{supp}({\mathcal{B}_{{\bf{i}},p}^{k}}) \cap  \hat{\Omega}_{k}({\bm{\mu}})  \ne \varnothing \}, 
\end{equation}
where $ \hat{\Omega}_{k}({\bm{\mu}})$ is the $k$-th parametric trimmed domain in~\eqref{eq13b}.
Following the choice~\eqref{eq17a}, the multi-patch discrete space can be then written as
\begin{equation}\label{eq26}
V_{h}({\bm{\mu}}) = \bigoplus_{k=1}^{N_p} V_{h,k}({\bm{\mu}}),\quad\text{with}\quad V_{h,k}({\bm{\mu}}) = \{ \boldsymbol{v}_h \in S^p({{{{\Omega_k}}({\bm{\mu}})}})^3  : \left.\boldsymbol{v}_h\right\rvert_{\Gamma_D({\bm{\mu}})} = \bm{0}\}.
\end{equation}

As in the single-patch case, and unless stated otherwise, we assume the spaces $V_{h,k}({\bm{\mu}})$ to be $C^1$-continuous within each patch.
Nevertheless, the multi-patch space $V_{h}({\bm{\mu}})$ is discontinuous across the interfaces $\gamma_{j}({\bm{\mu}})$.
The continuity required by problem~\eqref{eq17} will be imposed in a weak sense.
Thus, let us now define the coupling conditions for each interface $\gamma_{j}({\bm{\mu}})$.
First, we denote the displacement fields restricted to $\partial{\Omega_{k}}({\bm{\mu}})$ and $\partial{\Omega_{l}}({\bm{\mu}})$ as $\boldsymbol{u}_{k}({\bm{\mu}})$ and $\boldsymbol{u}_{l}({\bm{\mu}})$, respectively.  Then the coupling conditions can be expressed as
\begin{equation}\label{eq28}
	\begin{aligned}
	\boldsymbol{u}_{k}({\bm{\mu}}) - \boldsymbol{u}_{l}({\bm{\mu}}) &= 0 \quad \text{on} \  \gamma_{j}({\bm{\mu}}), \\
	\theta_n(\boldsymbol{u}_{k}({\bm{\mu}})) - 	\theta_n(\boldsymbol{u}_{k}({\bm{\mu}})) & = 0 \quad \text{on} \  \gamma_{j}({\bm{\mu}}),
		\end{aligned}
\end{equation}	
that, by means of the standard jump operators we can rewrite as
\begin{equation}\label{eq29}
\begin{aligned}
	\llbracket\boldsymbol{u}({\bm{\mu}}) \rrbracket= 0 \quad \text{on} \  \gamma_{j}({\bm{\mu}}), \\
\llbracket	\theta_n(\boldsymbol{u}({\bm{\mu}})) \rrbracket= 0 \quad \text{on} \  \gamma_{j}({\bm{\mu}}).
\end{aligned}	
\end{equation}

There exist several  coupling strategies for the imposition of the above continuity constraints in a weak sense. In the numerical experiments of Section~\ref{sec7}, we use both a super-penalty~\cite{Coradello2021} and Nitsche's~\cite{Benzaken2021} methods to achieve displacement and rotational continuity, both described below.
In principle, the reduced order modeling techniques to be applied (Section~\ref{sec5}) are also suitable for other coupling strategies, such as the mortar method where Lagrange multipliers can be eliminated in the context of the reduced basis method~\cite{Horger2017}.

The weak formulation~\eqref{eq17} of the parameterized problem can be extended to the multi-patch case by enriching the bilinear form with additional terms that weakly enforce the coupling conditions. Thus, the shell multi-patch problem reads:  find ${\boldsymbol{u}}_h({\bm{\mu}}) \in V_{h}({\bm{\mu}})$ such that 
\begin{equation}\label{eq44}
{\sum_{k=1}^{N_p}}	a_k({\boldsymbol{u}}_h,{\boldsymbol{v}}_h;\bm{\mu}) +  {\sum_{j=1}^{N_{\Gamma}}} a_j^{\Gamma}({\boldsymbol{u}}_h,{\boldsymbol{v}}_h;\bm{\mu}) = {\sum_{k=1}^{N_p}}f_k({\boldsymbol{v}}_h;\bm{\mu}), \qquad \forall {\boldsymbol{v}}_h \in V_h({\bm{\mu}}),
\end{equation}
where $a_k$ and $f_k$ are the bilinear and linear parameterized forms~\eqref{eq17b} and~\eqref{eq17c}, respectively, restricted to the patch $k$, while the coupling terms $a_j^{\Gamma}$ will be further  discussed in Sections~\ref{sec4a} and~\ref{sec4b}. The discretization yields the following parameterized linear system of dimension $\mathcal{N}_h({\bm{\mu}})=\text{dim}\left(V_h({\bm{\mu}})\right)$
\begin{equation}\label{eq21}
	{\bf{K}}(\bm{\mu}){\bf{u}}(\bm{\mu}) = {\bf{f}}(\bm{\mu}),
\end{equation}
where ${\bf{K}}$ $\in \mathbb{R}^{\mathcal{N}_h({\bm{\mu}})\times \mathcal{N}_h({\bm{\mu}})}$ is the global stiffness matrix, ${\bf{f}} \in \mathbb{R}^{\mathcal{N}_h({\bm{\mu}})}$ is the force vector, and $\mathcal{N}_h({\bm{\mu}})$ is the global number of degrees of freedom. The problem~\eqref{eq44}-\eqref{eq21} is the high-fidelity or full order model (FOM) upon which we aim to construct a reduced order model (ROM) for fast parametric simulations. 

In what follows, we describe the super-penalty \cite{Coradello2021} and Nitsche's \cite{Benzaken2021} coupling methods (Sections~\ref{sec4a} and~\ref{sec4b}, respectively).
While the first approach is simpler and more efficient, its application is limited to the case in which each interface can be associated to a parametric face for at least one of the two contiguous patches at the interface.
On the other hand, Nitsche's method overcomes such limitation, but is more convoluted as it involves further terms and third-order derivatives, and, consequently, requires $C^2$-discretization spaces.

\subsection{Projected super-penalty approach}\label{sec4a}
 In what follows, we will briefly present the projected super-penalty approach that we will use to impose the continuity constraints~\eqref{eq29} for some of the numerical experiments in Section~\ref{sec7}. The reader is further referred to \cite{Coradello2021a,Coradello2021} for a more detailed overview of the coupling approach. 

At each interface $\gamma_j({\bm{\mu}})$, for $j=1,\dots,N_{\Gamma}$, let us first choose arbitrarily one of the neighboring patches as active.
This active patch must be chosen such that $\gamma_j({\bm{\mu}})$ is fully contained in one of its four patch parametric faces.
We extract the patch knot vector associated to that face and remove the first and last knots, denoting the resulting vector as $\Xi_j$.
Then, using $\Xi_j$ we construct a one-dimensional lower degree spline space $S^{p-2}\left(\gamma_j(\boldsymbol{\mu})\right)$.
Finally, we denote as $\Pi_j$ the $L^2$-projection, related to the interface $\gamma_j(\boldsymbol{\mu})$, onto the reduced vector space $S^{p-2}\left(\gamma_j(\boldsymbol{\mu})\right)^d$ for displacements, and onto $S^{p-2}\left(\gamma_j(\boldsymbol{\mu})\right)$ for normal rotations. The discretized bilinear form is enriched with penalty terms that weakly enforce the coupling conditions~\eqref{eq29} and exploit the properties of the $L^2$-projection.
Thus, these coupling terms read
\begin{equation}\label{eq30}
\begin{aligned}
a_j^{\Gamma}({\boldsymbol{u}}_h,{\boldsymbol{v}}_h;\bm{\mu}) &=
 c_{\textrm{disp}}^{(j)} \int_{{\gamma_j}(\boldsymbol{\mu})} \Pi_j  \llbracket\boldsymbol{u}_h({\bm{\mu}}) \rrbracket \cdot  \Pi_j  \llbracket\boldsymbol{v}_h({\bm{\mu}}) \rrbracket\, \mathrm{d}\Gamma  \\
&+  c_{\textrm{rot}}^{(j)}  \int_{{\gamma_j}(\boldsymbol{\mu})} \Pi_j  \llbracket\theta_n(\boldsymbol{u}_h({\bm{\mu}})) \rrbracket \,  \Pi_j  \llbracket\theta_n(\boldsymbol{v}_h({\bm{\mu}})) \rrbracket\, \mathrm{d}\Gamma.
\end{aligned}
\end{equation}
The parameters of the penalty method are chosen as:
\begin{equation}\label{eq31}
	\begin{aligned}
		c_{\textrm{disp}}^{(j)} &= \left(\lvert\gamma_j(\boldsymbol{\mu})\rvert\right)^{c_{\textrm{exp}}-1}\frac{Et}{(h_j)^{{c_{\textrm{exp}}}}(1-\nu^2)}, \\ 
		c_{\textrm{rot}}^{(j)} &= \left(\lvert\gamma_j(\boldsymbol{\mu})\rvert\right)^{c_{\textrm{exp}}-1}\frac{Et^3}{(12h_j)^{c_{\textrm{exp}}}(1-\nu^2)},
	\end{aligned}
\end{equation}
where $h_j$ denotes the interface mesh size and the measure $\lvert\gamma_j(\boldsymbol{\mu})\rvert$ is the length of the coupling interface $\gamma_j(\boldsymbol{\mu})$. The exponent factor $c_{\textrm{exp}}$ is chosen as $c_{\textrm{exp}}=p-1$ in the numerical experiments of Section~\ref{sec7}, which yields optimal convergence of the method in the $H^2$ norm. A detailed discussion on the choice of $c_{\textrm{exp}}$ in relation to optimal convergence and conditioning of the underlying system of equations is given in \cite{Coradello2021}. 

We remark that this coupling approach is locking-free at interfaces and the choice of penalty parameters can be done automatically based on the problem setup. Nevertheless, the computational cost grows for higher order splines ($p > 3$) and the condition number related to the chosen penalty parameters may affect the accuracy of the solution. In addition, the coupling at interfaces where both sides are trimmed is still an open issue.
To this end, we will consider the Nitsche's method for the coupling of patches in more general cases of complex geometries.
\begin{remark}\label{rmk:coradello2021}
To overcome such problem, in \cite{Coradello2021} the internal knots of the trimming interfaces are neglected for the computation of the intersection mesh at the interface, which can yield sub-optimal results in particular for non-smooth interfaces.
\end{remark}

\subsection{Nitsche's method}\label{sec4b}
Let us now introduce the Nitsche's method for recovering $C^1$-continuity at the multi-patch interfaces.
A detailed overview of the method is given in~\cite{Benzaken2021}.
The coupling conditions are imposed in a weak sense by augmenting the discretized bilinear form with penalty, consistency, and symmetry terms.
Thus, the coupling terms in~\eqref{eq44} read
\begin{equation}\label{eq42}
	\begin{aligned}
 a_j^{\Gamma}&({\boldsymbol{u}}_h,{\boldsymbol{v}}_h;\bm{\mu}) =
  \underbrace{\int_{{\gamma_j}(\boldsymbol{\mu})}\langle{\bf{T}}(\boldsymbol{u}_h)\rangle_{\delta} \cdot \llbracket{\boldsymbol{v}_h}\rrbracket +
  \langle{B_{nn}(\boldsymbol{u}_h)}\rangle_{\delta} \llbracket\theta_n(\boldsymbol{v}_h)\rrbracket \mathrm{d}\Gamma}_{\text{Consistency terms}} \\
&+\underbrace{\int_{{\gamma_j}(\boldsymbol{\mu})}
  \langle{\bf{T}}(\boldsymbol{v}_h)\rangle_{\delta} \cdot \llbracket{\boldsymbol{u}_h}\rrbracket+ \langle{B_{nn}(\boldsymbol{v}_h)}\rangle_{\delta} \llbracket\theta_n(\boldsymbol{u}_h)\rrbracket \mathrm{d}\Gamma}_{\text{Symmetry terms}} \\
 &+\underbrace{
 c_{\textrm{disp}}\int_{{\gamma_j}(\boldsymbol{\mu})}\llbracket\boldsymbol{v}_h\rrbracket\cdot\llbracket{\boldsymbol{u}_h}\rrbracket\mathrm{d}\Gamma
 +c_{\textrm{rot}}\int_{{\gamma_j}(\boldsymbol{\mu})}
     \llbracket\theta_n(\boldsymbol{v}_h)\rrbracket
     \llbracket\theta_n(\boldsymbol{u}_h)\rrbracket
     \mathrm{d}\Gamma}_{\text{Penalty terms}}.
	\end{aligned}
\end{equation}
Here, $\langle\cdot\rangle_{\delta}$ denotes the average operator, $B_{nn}=\boldsymbol{n} \cdot   B(\boldsymbol{u}_h) \boldsymbol{n} $ is the bending moment and $\bf{T}$ is the ersatz force, defined in \cite{Benzaken2021} as
\begin{equation}\label{eq42b}
	{\bf{T}} = A(\boldsymbol{u}_h)  \boldsymbol{n}  - \boldsymbol{b} \left( B(\boldsymbol{u}_h) \boldsymbol{n} + B_{nt} (\boldsymbol{u}_h)   \boldsymbol{t}\right) + \left[ \left(\nabla \cdot B(\boldsymbol{u}_h)\right) \cdot \boldsymbol{n} + \frac{\partial B_{nt}(\boldsymbol{u}_h)}{\partial \boldsymbol{t}} \right] {\bf{a}}_3,
\end{equation}
where $\boldsymbol{b} = - \nabla {\bf{a}}_3$ is the curvature tensor, $\boldsymbol{t}$ is the counter-clockwise positively-oriented unit tangent vector to $\Gamma(\boldsymbol{\mu})$, and $B_{nt}(\boldsymbol{u}_h)  = \boldsymbol{n} \cdot B(\boldsymbol{u}_h) \boldsymbol{t}$ is the twisting moment.
The average operator is defined as:
\begin{equation}\label{eq:averageop}
\langle a \rangle_{\delta} = {\delta}\,\left. a\right\vert_{\Omega_k(\bm{\mu})} + (1-{\delta})\,\left. a\right\vert_{\Omega_l(\bm{\mu})},\quad\text{with }{\delta}\in[1/2,1],
\end{equation}
where $a$ is an arbitrary function defined over $\Omega(\bm{\mu})$, and $\left. a\right\vert_{\Omega_k(\bm{\mu})}$ and $\left. a\right\vert_{\Omega_l(\bm{\mu})}$ its restrictions to the patches $\Omega_k(\bm{\mu})$ and $\Omega_l(\bm{\mu})$ at the interface $\gamma_j(\bm\mu)$.
In the numerical experiments of Section~\ref{sec7}, the penalty constants are chosen as
\begin{equation}\label{eq43}
		c_{\textrm{disp}} = 10^3 \frac{Et}{h}, \qquad
		c_{\textrm{rot}}= 10^3 \frac{Et^3}{h}.
\end{equation}

We remark that this coupling approach is variationally consistent and stable (see Remark~\ref{rmk:stability} below).
In addition, its discretized weak formulation results in a well-conditioned system of linear system and is more robust with respect to the chosen parameters compared to penalty approaches.  
Nevertheless, it is easy to realize that the ersatz force $\bf{T}$ requires third-order derivatives (terms $\nabla \cdot B(\boldsymbol{u}_h)$ and $\partial B_{nt}(\boldsymbol{u}_h)/\partial \boldsymbol{t}$ in Equation~\eqref{eq42b}), what implies the necessity of $C^2$-continuous discretization spaces.
Furthermore, the number of terms and the order of the involved derivatives make this method's implementation more convoluted than the super-penalty approach (recall Section~\ref{sec4a}).

\begin{remark}
The ersatz force $\bf{T}$ in Equation~\eqref{eq42b} involves third-order derivatives (terms $\nabla \cdot B(\boldsymbol{u}_h)$ and $\partial B_{nt}(\boldsymbol{u}_h)/\partial \boldsymbol{t}$), what implies the necessity of $C^2$-continuous discretization spaces.
In particular, and as detailed in~\cite{Benzaken2021}, the patch approximation spaces $V_{h,k}({\bm{\mu}})$ in~\eqref{eq26} must be chosen such that 
\begin{equation}
V_{h,k}({\bm{\mu}})  \subset \{ (u_1,u_2) \in H^1({{{{\Omega_k}}({\bm{\mu}})}})^2 \  \text{and} \ u_3 \in H^2({{{{\Omega_k}}({\bm{\mu}})}})  \},
\end{equation}
where $u_1, u_2, u_3$ are contravariant components of the displacement, i.e., ${\boldsymbol{u}_h} = u_1{\bf{a}}^1 + u_2{\bf{a}}^2 + u_3{\bf{a}}^3$.
Consequently, in the numerical examples included in Section~\ref{sec7}, whenever Nitsche's method is applied, $C^2$-continuous spline spaces are considered.
\end{remark}
\begin{remark}\label{rmk:stability}
As discussed, for instance, in~\cite{antolin2021overlapping}, in the case of trimmed domains that present small elongated cut elements at the interface, instability effects may arise in the evaluation of the normal fluxes in the formulation~\eqref{eq42}.
This problem can be easily overcome in the case in which one of the patches is not trimmed at the interface, by just selecting $\delta=1$ in~\eqref{eq:averageop}, i.e., computing the flux on the non-trimmed side only.
However, this is not possible when both patches are trimmed and both present small elongated cut elements at the interface.

A possible alternative is the use of stabilization techniques that have been recently proposed for the Nitsche method in the context of isogeometric discretizations~\cite{Elfverson2019,antolin2021overlapping,wei2023stabilized}.
However, and to the best of our knowledge, no mathematically sound stabilization method has been proposed yet for isogeometric Kirchhoff-Love shells.
It is our belief that the minimal stabilization technique proposed in our previous work~\cite{antolin2021overlapping} may be handy in stabilizing the problem.
Nevertheless, this study is out of the scope of this paper and has not been addressed yet.
In the numerical experiments included in Section~\ref{sec7}, no numerical instabilities were found.
\end{remark}

\section{Local Reduced Basis method}\label{sec5}
We are interested in problems where a large amount of solution evaluations are required for problem~\eqref{eq21} for different values of the parameters vector $\boldsymbol{\mu}$. The use of reduced order modeling techniques, such as the reduced basis method, can speedup the computation of parameterized problems.  The main idea is based on an offline/online split: In the offline phase, the first step is to compute snapshots of the FOM and extract a linear combination of reduced basis functions using techniques such as the greedy algorithm or the Proper Orthogonal Decomposition method. Then, a ROM is constructed by orthogonal projection into the subspace spanned by the reduced basis. In the online phase,  a reduced problem is solved to obtain the solution for any given parameter. The reader is further referred to \cite{QMN_RBspringer, Hesthaven2016} for more details on the reduced basis method. Nevertheless, the application of standard reduced order modeling techniques on parameterized trimmed domains entails several challenges:
\begin{itemize}
	\item The spline space $V_h({\bm{\mu}})$ and its dimension $\mathcal{N}_h({\bm{\mu}})$ depend on the geometric parameters ${\bm{\mu}}$. In particular, the set of active basis functions may change for different snapshots depending on the value of the parameters ${\bm{\mu}}$. This impedes the construction of snapshots matrices to extract a reduced basis, since snapshots may be vectors of different length, and requires suitable snapshots extension \cite{Karatzas2020}. 
	\item The efficiency of the offline/online decomposition is based on the assumption that operators depend affinely on the parameters. In case of geometrically parameterized problems, this assumption does not always hold and has to be recovered by constructing affine approximations with efficient hyper-reduction techniques \cite{Barrault2004,Negri2015}. 
	\item The solution manifold obtained by extended snapshots is highly nonlinear with respect to the parameters $\bm{\mu}$.  The same holds also for the affine approximations of the parameter-dependent operators. The approximation of a nonlinear manifold with a single, linear reduced basis space may yield a very high dimension of the basis. This requires tailored strategies to ensure the efficiency of the method.
\end{itemize}
In the following we will briefly review a local reduced basis method to construct efficient ROMs on parameterized trimmed domains based on clustering strategies and the Discrete Empirical Interpolation Method (DEIM). The reader is further referred to our previous work \cite{Chasapi2023} for a more detailed exposition.  The main steps of the local method that will be discussed in the following sections are summarized as follows:
\begin{enumerate}
\item perform a trivial extension of snapshots and define the extended FOM on a common, background mesh,
\item cluster the parameters and associated snapshots in order to construct DEIM and reduced basis approximations,
\item in the offline phase, train local DEIM approximations and reduced bases for each cluster combination and perform the projection,
\item choose the cluster with the smallest distance to a given parameter during the online phase and solve a reduced problem of sufficiently small dimension. 
\end{enumerate}

\subsection{Snapshots extension}\label{sec:extension}
Let us first discuss the construction of snapshots. The domain at hand $\Omega(\bm{\mu})$ depends on $\bm{\mu}$, and thus the support of the B-spline basis functions is also parameter-dependent. Since both the spline space $V_h({\bm{\mu}})$ and its dimension $\mathcal{N}_h({\bm{\mu}})$ depend on the geometric parameters ${\bm{\mu}}$, the dimension of each solution snapshot may differ from one parameter to the other. Therefore, the first step is to perform a suitable extension of the solution vector ${\bf{u}}(\bm{\mu})$ such that all solution vectors have the same length. As the original non-trimmed domain $\Omega_0$ is independent of the parameters, the natural choice is to use it as a common background domain where snapshots are extended. In this work, we consider a zero extension although other choices are also possible \cite{Karatzas2020}.
This extension is performed over the non-trimmed, multi-patch space $V_{h,0}=\oplus_{k=1}^{N_p} (S^p_{k,0})^3$ (recall definition~\eqref{eq24}), whose dimension is $\mathcal{N}_{h,0}=\text{dim}\left(V_{h,0}\right)$.

With these definitions at hand, the extended full order problem in~\eqref{eq21} becomes:
\begin{equation}\label{eq46}
	\widehat{{\bf{K}}}(\bm{\mu})\widehat{{\bf{u}}}(\bm{\mu}) = \widehat{{\bf{f}}}(\bm{\mu}),
\end{equation}
where $\widehat{\bf{K}}$ $\in \mathbb{R}^{\mathcal{N}_{h,0}\times \mathcal{N}_{h,0}}$ is the extended stiffness matrix, $\widehat{\bf{u}}(\bm{\mu})  \in \mathbb{R}^{\mathcal{N}_{h,0}}$ is the extended solution vector and $\widehat{\bf{f}}(\bm{\mu}) \in \mathbb{R}^{\mathcal{N}_{h,0}}$ is the extended right-hand side vector. Therefore, the size of the above extended problem is $\bm{\mu}$-independent and is suitable for the computation of snapshots.

\subsection{Clustering strategy}\label{sec:cluster}
In what follows we will briefly review the clustering strategy to construct local ROMs. The main idea is to build multiple approximations based on smaller subspaces instead of one global space. Then in the online phase, the \textit{closest} cluster is selected for a given parameters vector $\bm{\mu}$ and a local reduced problem is solved. 

For geometrically parameterized problems on trimmed domains we opt for a parameter-based clustering \cite{Chasapi2023}. We seek for a partition of the parameter space $\mathcal{P}$ in $N_c$ subspaces, such that 
\begin{align}\label{eq47a}
	\mathcal{P} = \bigcup_{k=1}^{N_c} \mathcal{P}_k.
\end{align}
In this work, we will use the \emph{k-means} clustering algorithm \cite{Likas2003} to partition
the parameter space
although other partitioning strategies are also possible \cite{Eftang2010a,Eftang2011,Haasdonk2011}. The main idea is to assign a given parameter vector ${\bm{\mu}}$ to the cluster that minimizes the maximum distance  $\mathcal{D}({\bm{\mu}},k)$ between boundaries of the trimmed regions $\Omega_i(\bm{\mu}), \ i=1,\dots,K$, in~\eqref{eq13}, where $\exists C > 0$ 
\begin{align}\label{eq47}
	\mathcal{D}({\bm{\mu}},k) = \max_i \text{dist}(\partial{\hat{\Omega}}_i(\bar{\bm{\mu}}_k)\ \cap \ \partial{\hat{\Omega}}_i({\bm{\mu}})) \le C \norm{{\bm{\mu}} - {\bar{\bm{\mu}}}_k}_2^2, \quad k=1,\dots,N_c,
\end{align}
and $\bar{\bm{\mu}}_k$ is the centroid of the $k$-th cluster\footnote{\label{note1}Hereinafter, and for the sake of conciseness, whenever it is clear from the context, the range $1,\dots,N_c$ of the cluster index $k$ will be omitted.}.
Then, the parameter space $\mathcal{P}$ is partitioned into $N_c$ subspaces as
\begin{equation}\label{eq48}
	{\mathcal{P}_k}= \{{\bm{\mu}}\in\mathcal{P} \ :\ \underset{i=1,\dots,N_c}{\arg \min} \norm{{\bm{\mu}} - {\bar{\bm{\mu}}}_{i}}_2^2=k\}.
\end{equation}
Thus, the k-means clustering minimizes the distance between each parameter vector and the cluster's centroid with respect to the Euclidean norm $\norm{\cdot}_2$.
This strategy partitions trimmed discretizations that have similar active and inactive regions.

For performing such partition, we work with a discrete counterpart of the continuous space $\mathcal{P}$.
Thus, we create a sufficiently fine and properly selected training sample set $\mathcal{P}^s=\{\bm{\mu}_1,\dots,\bm{\mu}_{N_s}\} \subset \mathcal{P}$ of dimension $N_s=\text{dim}(\mathcal{P}^s)$, for which the $N_c$ centroids are sought and updated iteratively until the algorithm converges.
We refer the interested reader to \cite[Algorithm~5]{Amsallem2012} for a detailed overview of the k-means algorithm. Once the parameter space partition~\eqref{eq47a} is created, the training set $\mathcal{P}^s$ can be clustered accordingly as
\begin{equation} \label{trainpartition}
\mathcal{P}^s = \bigcup_{k=1}^{N_c}\mathcal{P}^s_{k}\quad\text{with}\quad \mathcal{P}^s_{k}=\mathcal{P}^s \cap \mathcal{P}_{k}.
\end{equation}

Note that a suitable number of clusters $N_c$ should be chosen in advance to perform the k-means algorithm.
The suitability of this choice can be evaluated a posteriori by considering the k-means variance as
\begin{equation}\label{eq51a}
	{\bf{\mathcal{V}}} = \sum_{k=1}^{N_c} \sum_{{\bm{\mu}} \in \mathcal{P}^s_{k}} \norm{{\bm{\mu}} - \bar{{\bm{\mu}}}_k}_2^2.
\end{equation}
In particular, the k-means variance is expected to decrease with increasing number of clusters and the smallest integer can be chosen as $N_c$ at the transition between a steep slope and a plateau \cite{Hess2019}. Note that in the numerical experiments of Section~\ref{sec7}, this transition occurs at $N_c \le 10$.

Once the clustering is performed, local ROMs are constructed in the offline phase. The construction of the localized ROM will be discussed in Sections~\ref{sec:localRB} and~\ref{sec:deim}, while a detailed overview of the offline phase of the algorithm is given in \cite[Algorithms~1-2]{Chasapi2023}.
Afterwards, in the online phase, for a given parameters vector $\bm\mu$, we determine the corresponding cluster $\mathcal{P}_k$ according to~\eqref{eq48} and select its associated local ROM for solving a reduced problem. This online phase of the algorithm is also presented in more detail in \cite[Algorithm~3]{Chasapi2023}.

\subsection{Proper Orthogonal Decomposition}\label{sec:pod}
In the following we will briefly recall the Proper Orthogonal Decomposition (POD) that we will use for the construction of the ROM later on. The POD is based on the singular value decomposition algorithm (SVD) and aims to extract a set of orthonormal basis functions \cite{Quarteroni2017}. The SVD of a matrix $	{\mathbb{S}} \in \mathbb{R}^{m \times n}$ reads
\begin{equation}\label{eq59}
	{\mathbb{S}}= \mathbb{U}\boldsymbol{\Sigma}\mathbb{Z}^{\top},
\end{equation}
where the orthogonal matrices $\mathbb{U}=[\bm{\zeta}_1,\dots,\bm{\zeta}_{m}] \in \mathbb{R}^{m \times m}$, $\mathbb{Z}=[\bm{\psi}_1,\dots,\bm{\psi}_{n}] \in \mathbb{R}^{n \times n}$ have columns containing the left and right singular vectors of ${\mathbb{S}}$, respectively, $\boldsymbol{\Sigma}\in \mathbb{R}^{m \times n}$ is a rectangular diagonal matrix that contains the singular values $\sigma_1 \ge \sigma_2 \ge \dots \sigma_r$, and $r \le \min(m, n)$ is the rank of $\mathbb{S}$. The POD basis of dimension $N$ is then defined as the set of the first $N$ left singular vectors of $\mathbb{S}$, i.e., the $N$ largest singular values. We can choose the dimension of the basis $N$ such that the error in the POD basis is smaller than a prescribed tolerance $\varepsilon_{\text{POD}}$ \cite{QMN_RBspringer}, namely $N$ is the smallest integer such that
\begin{equation}\label{eq61}
	1 - \frac{\sum_{i=1}^{N}\sigma_i^2}{\sum_{j=1}^{r}\sigma_j^2} \le \varepsilon^2_{\text{POD}}.
\end{equation}
Note that the POD basis is orthonormal by construction and the basis functions can be understood as modes  that retain most of the energy of the original system. The dimension of the latter is then reduced such that the energy captured by the neglected modes in smaller than or equal to $\varepsilon_{\text{POD}}$ in~\eqref{eq61}.

\subsection{Local Reduced Basis problem}\label{sec:localRB}
In order to construct an efficient ROM, we seek for local reduced bases ${\bf{V}}_k \in \mathbb{R}^{\mathcal{N}_{h,0}\times N_k}$ for every cluster $k$, where $N_k$ is the local reduced space dimension that is ideally of sufficiently small dimension, i.e., $N_k \ll \mathcal{N}_{h,0}$. The reduced basis is constructed in the offline phase separately for each cluster based on the strategy discussed in Section 
\ref{sec:cluster}. In this work we will consider the POD to construct each reduced basis ${\bf{V}}_k$, while its construction was briefly discussed in Section~\ref{sec:pod}. Note that other techniques can in principle be also chosen for the construction, as, e.g., the greedy algorithm \cite{Chasapi2022}.

To construct a POD basis, let us again consider a fine training sample set $\mathcal{P}^s = \{\bm{\mu}_1,\dots,\bm{\mu}_{N_s}\} \subset \mathcal{P}$ with dimension $N_s=\text{dim}(\mathcal{P}^s)$ introduced in Section~\ref{sec:cluster}.
Then we form the snapshots matrix ${\mathbb{S}} \in \mathbb{R}^{\mathcal{N}_{h,0} \times N_s}$ as
\begin{equation}\label{eq58}
	{\bf{S}}_u= [\widehat{\bf{u}}_1,\dots,\widehat{\bf{u}}_{N_s}],
\end{equation}
where the vectors $\widehat{\bf{u}}_j \in \mathbb{R}^{\mathcal{N}_{h,0}}$ denote the extended solutions $\widehat{\bf{u}}{(\bm{\mu}_j)}$  for $j=1,\dots,N_s$.
These snapshots are also partitioned into $N_c$ submatrices as $\lbrace {\bf{S}}_u^1,\dots,{\bf{S}}_u^{N_c}\rbrace$, according to $\mathcal{P}^s = \cup_{k=1}^{N_c}\mathcal{P}^s_{k}$ in~\eqref{trainpartition}.
The local reduced basis ${\bf{V}}_k$ is then extracted from each cluster ${\bf{S}}_u^k$, separately applying the POD as discussed in Section~\ref{sec:pod}. Thus, the basis reads:
\begin{equation}\label{eq60}
	{\bf{V}}_k = [\bm{\zeta}_1,\dots,\bm{\zeta}_{N_k}] \in \mathbb{R}^{\mathcal{N}_{h,0} \times N_k}.
\end{equation}

Let us now derive the local reduced basis problem. For any ${\bm\mu} \in \mathcal{P}_k$, the solution $\widehat{{\bf{u}}}(\bm{\mu})$ can be approximated using  the local reduced basis ${\bf{V}}_k$, as 
\begin{equation}\label{eq52}
	\widehat{{\bf{u}}}(\bm{\mu}) \approx {\bf{V}}_k {\bf{u}}_N(\bm{\mu}),
\end{equation}
where ${\bf{u}}_N(\bm{\mu})\in \mathbb{R}^{N_k}$ is the solution vector of the reduced problem. A projection-based ROM can be obtained from~\eqref{eq46} by enforcing the residual to be orthogonal to the subspace spanned by ${\bf{V}}_k$, such that
\begin{equation}\label{eq53}
	{\bf{V}}_k^{\top} (\widehat{{\bf{K}}}(\bm{\mu}){\bf{V}}_k  {\bf{u}}_N(\bm{\mu})-{\widehat{{\bf{f}}}(\bm{\mu})}) = {\bf{0}}.
\end{equation}
Thus, the local reduced basis problem reads:
\begin{equation}\label{eq54}
	{\bf{K}}_N(\bm{\mu}) {\bf{u}}_N(\bm{\mu}) = {\bf{f}}_N(\bm{\mu}),
\end{equation} 
while the reduced matrix and vector are defined as
\begin{equation}\label{eq55}
	{\bf{K}}_N = {\bf{V}}_k^{\top}\widehat{{\bf{K}}}(\bm{\mu}){\bf{V}}_k , \qquad {\bf{f}}_N = {\bf{V}}_k^{\top}\widehat{{\bf{f}}}(\bm{\mu}).
\end{equation}
The size of the reduced problem~\eqref{eq54} is $N_k \ll \mathcal{N}_{h,0}$, which makes it suitable for fast online computation given many different parameters $\bm{\mu} \in \mathcal{P}_k$. However, the solution of the reduced problem requires the assembly of the parameter-dependent operators $\widehat{{\bf{K}}}(\bm{\mu})$ and $\widehat{{\bf{f}}}(\bm{\mu})$. Therefore, an important assumption for the efficiency of the reduced basis method in general is that the operators depend affinely on the parameters ${\bm{\mu}}$. This assumption is not always fulfilled in the presence of geometrical parameters. Therefore, we will build affine approximations to recover the affine dependence. Since we aim to approximate a manifold of extended operators that is nonlinear with respect to $\bm{\mu}$, the dimension of the approximation space may be high. Therefore, the clustering strategy of Section~\ref{sec:cluster} will be also considered to construct local affine approximations. Note that from now on we assume for ease of exposition that the clustering is performed only once for constructing both the reduced bases and affine approximations, although in principle this could be chosen differently \cite{Chasapi2023}.  We now introduce the following local affine approximation for any ${\bm\mu} \in \mathcal{P}_k$:
\begin{equation}\label{eq56}
	\widehat{{\bf{K}}}(\bm{\mu}) \approx \sum_{q=1}^{Q_{a}^k} \theta_{a,q}^{k} (\bm{\mu})\widehat{{\bf{K}}}_q^k,
 \qquad \widehat{{\bf{f}}}(\bm{\mu})\approx\sum_{q=1}^{Q_f^k} \theta_{f,q}^k(\bm{\mu})\widehat{{\bf{f}}}_q^k,
\end{equation}
where $\theta_{a,q}^{k} : \mathcal{P}_k \to \mathbb{R} $, for $q=1,\dots,Q_{a}^k$, and $\theta_{f,q}^k : \mathcal{P}_k \to \mathbb{R} $, for $q=1,\dots,Q_f^k$, are $\bm{\mu}$-dependent functions, whereas $\widehat{{\bf{K}}}_q^k \in \mathbb{R}^{\mathcal{N}_{h,0}\times \mathcal{N}_{h,0}}$ and  $\widehat{{\bf{f}}}_q^k \in \mathbb{R}^{\mathcal{N}_{h,0}}$ are ${\bm{\mu}}$-independent forms\footnote{Henceforward, and for the sake of clarity, whenever it is clear from the context, the ranges $1,\dots,Q_{a}^k$ and $1,\dots,Q_{f}^k$ of the index $q$ referred to the terms of the local affine approximation~\eqref{eq56} will be omitted.}.
Since the latter forms do not depend on the parameters ${\bm{\mu}}$, they can be pre-computed and stored in the offline phase. Then, the online assembly  requires only the evaluation of $\theta_{a,q}^{k}, \theta_{f,q}^k$, which is inexpensive assuming that $Q_{a}^k, Q_f^k \ll \mathcal{N}_{h,0}$. To obtain the affine approximation in the form of~\eqref{eq56}, we will employ the Discrete Empirical Interpolation Method in combination with Radial Basis Functions Interpolation. This hyper-reduction strategy will be further discussed in Section~\ref{sec:deim}.  Once the affine approximation is recovered, inserting~\eqref{eq56} into~\eqref{eq55} yields, for any given parameter ${\bm\mu} \in \mathcal{P}_k$:
\begin{equation}\label{eq57}
	{\bf{K}}_N(\bm{\mu}) \approx\sum_{q=1}^{Q_{a}^k} \theta_{a,q}^{k} (\bm{\mu}){\bf{V}}_k^{\top}\widehat{{\bf{K}}}_q^k{\bf{V}}_k,
 \qquad {\bf{f}}_N(\bm{\mu}) \approx \sum_{q=1}^{Q_f^k} \theta_{f,q}^k(\bm{\mu}){\bf{V}}_k^{\top}\widehat{{\bf{f}}}^k_q,
\end{equation}
where ${\bf{K}}_N(\bm{\mu}) \in \mathbb{R}^{N_k \times N_k}$ and ${{\bf{f}}}_N(\bm{\mu}) \in  \mathbb{R}^{N_k}$ are the reduced matrix and right-hand side vector, respectively. We remark that in~\eqref{eq57} only the coefficients $\theta_{a,q}^{k}, \theta_{f,q}^k$ depend on the parameters $\bm{\mu}$ and are evaluated online, while all other quantities are assembled and stored in the offline phase.  Finally, during the online phase, for any given parameter $\bm{\mu}$ the respective cluster is selected as in~\eqref{eq48} and the local reduced problem in~\eqref{eq54} is solved considering the approximation assembly in~\eqref{eq57}. Finally, the high-fidelity approximation of the solution can be recovered through~\eqref{eq52}. We remark that the efficiency of the overall method depends on the size of the local reduced problem $N_k$, on the number of local affine terms  $Q_{a}^k, Q_f^k $, as well as the efficient online evaluation of the coefficients $\theta_{a,q}^{k}, \theta_{f,q}^k$. The latter aspects will be further elaborated in Section~\ref{sec:deim}.

\subsection{Local Discrete Empirical Interpolation Method}\label{sec:deim}
In this section we will briefly present the hyper-reduction strategy based on the Discrete Empirical Interpolation Method (DEIM) for matrices and vectors. The reader is further referred to \cite{Negri2015} for a detailed presentation of the method. 

As discussed before, the first crucial step for the efficiency of the ROM involves constructing local affine approximations in the form of Equation~\eqref{eq56}. These are constructed separately for each cluster during the offline phase. For ease of exposition, we consider the same training sample set $\mathcal{P}^s = \{\bm{\mu}_1,\dots,\bm{\mu}_{N_s}\} \subset \mathcal{P}$ of dimension $N_s$ as the one in~\eqref{eq58}, although other choices are also possible \cite{Chasapi2023}. Then we form the snapshots matrices ${{\bf{S}}_a} \in \mathbb{R}^{{{\mathcal{N}^2_{h,0}}} \times N_{s}}$ and ${{\bf{S}}_f} \in \mathbb{R}^{\mathcal{N}_{h,0} \times N_{s}}$ 
\begin{equation}\label{eq63}
	{\bf{S}}_a= [\widehat{\bf{k}}_1,\dots,\widehat{\bf{k}}_{N_{s}}], \qquad {\bf{S}}_f= [\widehat{\bf{f}}_1,\dots,\widehat{\bf{f}}_{N_{s}}],
\end{equation}
where the vectors $\widehat{\bf{k}}_i = \text{vec}(\widehat{{\bf{K}}} (\bm{\mu}_i))\in\mathbb{R}^{\mathcal{N}^2_{h,0}}$ and $\widehat{\bf{f}}_i=\widehat{\bf{f}}(\bm{\mu}_i)\in\mathbb{R}^{\mathcal{N}_{h,0}}$, with $i=1,\dots,N_s$,  denote the vectorization of the extended  stiffness matrix and the extended right-hand side vector, respectively.
Following the training sample set partitioning $\mathcal{P}^s = \cup_{k=1}^{N_c}\mathcal{P}^s_{k}$ in~\eqref{trainpartition}, the snapshots matrices ${\bf{S}}_a$ and ${\bf{S}}_f$ are also partitioned into $N_c$ submatrices $\lbrace {\bf{S}}_a^1,\dots,{\bf{S}}_a^{N_c}\rbrace$ and $\lbrace {\bf{S}}_f^1,\dots,{\bf{S}}_f^{N_c}\rbrace$, accordingly.
Then, we apply the POD to each submatrix to obtain the matrices $\widehat{{\bf{K}}}_q^k$ and vectors $\widehat{{\bf{f}}}_q^k$ in~\eqref{eq56}. Here, the number of affine terms $Q_a^k$ and $Q_f^k$ can be determined by prescribing a tolerance $\varepsilon_{\text{POD}}$ as in~\eqref{eq61}. It should be remarked that the latter should in general be lower than the tolerance used to construct the local reduced basis, so that the accuracy of the DEIM approximation does not impede the overall accuracy of the ROM \cite{QMN_RBspringer}. 

Now, let us discuss how to efficiently compute the parameter-dependent coefficients $\theta_{a,q}^{k} (\bm{\mu})$ and $\theta_{f,q}^{k} (\bm{\mu})$ in~\eqref{eq56} for each cluster. For this purpose, we will use the known as \emph{magic points} \cite{Maday2009} according to the empirical interpolation procedure \cite{Barrault2004}. For each local affine approximation $k$, a collection of $Q_a^k, Q_f^k$ entries is selected based on a greedy algorithm that minimizes the interpolation error over the snapshots \cite{QMN_RBspringer}. In what follows, the selected magic points are denoted as $\mathcal{J}_a^k, \mathcal{J}_f^k$ for the stiffness matrix and right-hand side vector, respectively. These entries fulfill exactly the following interpolation constraints for the stiffness matrix and right-hand side vector for each $\bm{\mu}\in\mathcal{P}_k$:
\begin{equation}\label{eq64}
	\begin{aligned}
	\sum_{q=1}^{Q_{a}^k} \theta_{a,q}^{k} (\bm{\mu})[\widehat{{\bf{K}}}^k_q]_{i,j} &= [\widehat{{\bf{K}}}(\bm{\mu})]_{i,j}, \quad \forall (i,j) \in \mathcal{J}_a^k, \\
	\sum_{q=1}^{Q_{f}^k} \theta_{f,q}^{k} (\bm{\mu})[\widehat{{\bf{f}}}^k_q]_{i} &= [\widehat{{\bf{f}}}(\bm{\mu})]_{i}, \qquad \forall i \in \mathcal{J}_f^k.
	\end{aligned}
\end{equation}
We remark that the two right-hand sides in the equations above require the online assembly of a collection of $Q_a^k$/$Q_f^k$ FOM matrix/vector entries for a given ${\bm{\mu}} \in \mathcal{P}_k$, which can be costly and  intrusive. Therefore, we will opt for interpolating with radial basis functions (RBFs) following our previous work \cite{Chasapi2023}. Similar to the construction of local affine approximations, local RBF-interpolants are constructed separately for each cluster $k$. The main idea is to compute offline the values of $\theta_{a,q}^{k} (\bm{\mu})$ and $\theta_{f,q}^{k} (\bm{\mu})$ as in~\eqref{eq64}, for each training sample ${\bm{\mu}} \in \mathcal{P}^s_k$, and train a fast interpolant using these computations. Then, in the online phase the local interpolants can be evaluated rapidly for any given ${\bm{\mu}} \in \mathcal{P}_k$, being the coefficients $\theta_{a,q}^{k} (\bm{\mu})$ in~\eqref{eq56} approximated as
\begin{equation}\label{eq66}
	{\theta}^k_{a,q}(\bm{\mu})\approx \sum_{j=1}^{N_s^k} \omega_{a,q,j}^{k} \phi_{q,j}^k (\norm{\bm{\mu}-\bm{\mu}_j}_2).
\end{equation}
where $N_s^k=\text{dim}(\mathcal{P}^s_{k})$ is the number of training samples associated to the $k$-th cluster, $\phi_{q,j}^k$ is the radial basis function associated to the $j$-th center parameter point $\bm{\mu}_j$ and $\norm{\cdot}_2$ denotes the Euclidean norm. For the numerical experiments in Section~\ref{sec7} we will use cubic RBFs, although other types of functions can be also chosen \cite{Buhmann2003}. During the offline phase, the unknown weights $\omega_{a,q,j}^{k}$ are computed separately for each cluster such that they fulfill the interpolation constraint exactly for ${\bm{\mu}_k} \in \mathcal{P}^s_{k}$
\begin{equation}\label{eq67}
	\sum_{j=1}^{N_{s}^k} \omega_{a,q,j}^{k} \phi_{q,j}^k (\norm{\bm{\mu}_k-\bm{\mu}_j}_2) = \theta_{a,q}^{k}(\bm{\mu}_k),\quad\text{with }\bm{\mu}_j\in\mathcal{P}^s_k.
\end{equation}
The procedure is identical for the coefficients ${\theta}^k_{f,q}$ associated to the right-hand side vector and is therefore omitted here. 

So far we have presented two approximations with respect to the FOM in~\eqref{eq21}, namely the local reduced problem as well as affine approximations of the extended stiffness matrix $\widehat{\bf{K}}(\bm{\mu})$ and right-hand side vector $\widehat{\bf{f}}(\bm{\mu})$. The construction of localized reduced bases and local affine approximations via DEIM allows to confine the dimension of the bases and number of affine terms, respectively. Moreover, the RBF-interpolation of the coefficients in~\eqref{eq66} enables a rapid evaluation in the online phase. It should be noted that the localized method requires additional offline effort to construct and store multiple bases, nevertheless, the main advantage is the reduction of the online computational cost.

\section{ROM-based parametric shape optimization}\label{sec6}
In this section, we will briefly review the aplication of the discussed reduction strategies to optimization problems, following closely the notation in \cite{Quarteroni2017}. Such problems require several evaluations of the solution and the objective function to be minimized, which can be expensive. Therefore, these can benefit from reduced basis approximations. 

First, let us assume that $\boldsymbol{\mu}$ controls the shape of the computational domain at hand  $\Omega(\boldsymbol{\mu})$ for the optimization process. As a first step, a reduced model is constructed offline for the design variables $\boldsymbol{\mu}$ following the approach presented in Section~\ref{sec5}. Then, the optimization is performed online. In what follows we will directly use the discrete reduced approximation ${\mathbf{u}}_N={\mathbf{u}}_N(\boldsymbol{\mu})$ of Equation~\eqref{eq52} for our exposition.

The parametric optimization problem  that we will consider in this paper is the compliance minimization under a given volume constraint, which is a common choice in structural optimization. The minimization of the compliance implies that the structure deforms less, i.e., it becomes stiffer. Let us now formulate the optimization problem at hand as:
\begin{equation}\label{eq33}
	\boldsymbol{\mu}_{opt} =  \underset{\boldsymbol{\mu} \in \mathcal{P}}{\operatorname{arg min}}\, J_N({\bf{u}}_N,\boldsymbol{\mu})  \ \ \text{such that} \ \ \  V(\bm{\mu}) \le V_0 , \ \boldsymbol{\mu}_{min} \le  \boldsymbol{\mu} \le \boldsymbol{\mu}_{max},
\end{equation}
where $V(\bm{\mu})$ is the volume of the domain $\Omega(\bm{\mu})$, $V_0$ a prescribed maximum volume, and $ \boldsymbol{\mu}_{min}, \boldsymbol{\mu}_{max}$ the lower and upper bounds for the design variable $\boldsymbol{\mu}$, respectively.
Here, the cost functional $J_N({\bf{u}}_N,\boldsymbol{\mu})$ is obtained by evaluating the reduced basis approximation of the problem solution.
For the compliance case, the reduced objective function reads 
\begin{equation}\label{eq32}
	J_N({\bf{u}}_N,\boldsymbol{\mu}) = 
	\frac{1}{2} {\bf{u}}_N(\boldsymbol{\mu})\cdot{\bf{f}}_N(\boldsymbol{\mu}).
\end{equation}
\begin{remark}
In the case of the extended FOM problem, the compliance functional~\eqref{eq32} can be written as
\begin{equation}\label{eq32bis}
	J(\widehat{\bf{u}},\boldsymbol{\mu}) = 
	\frac{1}{2} \widehat{\bf{u}}(\boldsymbol{\mu})\cdot\widehat{\bf{f}}(\boldsymbol{\mu}),
\end{equation}
that, after introducing the approximations~\eqref{eq52} and~\eqref{eq55}, becomes
\begin{equation}\label{eq32bis2}
	J(\widehat{\bf{u}},\boldsymbol{\mu}) \approx
	\frac{1}{2} ({\bf{V}}_k{\bf{u}}_N(\boldsymbol{\mu}))\cdot({\bf{V}}_k  {\bf{f}}_N(\boldsymbol{\mu}))= \frac{1}{2} {\bf{u}}_N(\boldsymbol{\mu})\cdot{\bf{f}}_N(\boldsymbol{\mu})=J_N({\bf{u}}_N,\boldsymbol{\mu}),
\end{equation}
where the orthogonality of the basis ${\bf{V}}_k$ was considered.
\end{remark}
There are several ways to solve the optimization problem~\eqref{eq33}-\eqref{eq32}. In the following we will consider a gradient-based approach, where the parameters $\boldsymbol{\mu}$ are updated in an iterative fashion depending on the gradient of the cost functional $J_N$. The gradient can be evaluated either analytically or based on a suitable approximation, e.g., with a finite difference scheme. The latter yields a black-box optimization approach that simply requires the reduced solution and the evaluation of the objective function. Nevertheless, the reduced model is perfectly suitable for computing parametric sensitivities due to its differentiability and affine parametric dependence \cite[Proposition~5.3]{QMN_RBspringer}. Let us now reformulate the parameterized reduced problem as
\begin{equation}\label{eq34}
	\mathbf{G}_N({\bf{u}}_N,\boldsymbol{\mu}) =  {\bf{K}}_N(\boldsymbol{\mu}){\bf{u}}_N(\boldsymbol{\mu})-  {\bf{f}}_N(\boldsymbol{\mu}) =\bf{0}.
\end{equation}
With these definitions at hand and following \cite[Proposition~11.3]{QMN_RBspringer},  the gradient of the objective function $\tilde{J}_N(\boldsymbol{\mu})=J_N({\bf{u}}_N,\boldsymbol{\mu})$ reads:
\begin{equation}\label{eq35}
	\nabla_{\boldsymbol{\mu}} \tilde{J}_N(\boldsymbol{\mu}) = 	\nabla_{\boldsymbol{\mu}} J_N({\bf{u}}_N,\boldsymbol{\mu}) + 	\nabla_{{\bf{u}}_N} J_N({\bf{u}}_N,\boldsymbol{\mu}) \frac{\partial {\bf{u}}_N} {\partial 	\boldsymbol{\mu}}.
\end{equation}
The evaluation of these gradient requires the solution of $M$ sensitivity equations 
\begin{equation}\label{eq39}
	\frac{\partial {\bf{u}}_N} {\partial 	\boldsymbol{\mu}}  = - 	\left(D_{{\bf{u}}_N} \mathbf{G}_N({\bf{u}}_N,\boldsymbol{\mu})\right)^{-1} 	\left(D_{\boldsymbol{\mu}} \mathbf{G}_N({\bf{u}}_N,\boldsymbol{\mu})\right) 
\end{equation}
where $D$ denotes the Fr\'echet derivative. Instead of directly solving the above equations at each step of the optimization process, we define an additional adjoint problem and denote its solution as $\overline{\bf{u}}_N =\overline{\bf{u}}_N(\boldsymbol{\mu}) $ such that
\begin{equation}\label{eq40}
	\left(D_{{\bf{u}}_N} \mathbf{G}_N({\bf{u}}_N,\boldsymbol{\mu})\right)^{\top} \overline{\bf{u}}_N =   \nabla_{\boldsymbol{u}_N} J_N({\bf{u}}_N,\boldsymbol{\mu}).
\end{equation}
It should be remarked that this approach requires in general the construction of a reduced model for the adjoint problem, which implies additional offline cost. This can be constructed following the approach in Section~\ref{sec5}.
Given the reduced adjoint approximation $\overline{\bf{u}}_N$, the gradient in Equation~\eqref{eq35} can be reformulated as:
\begin{equation}\label{eq41}
	\nabla_{\boldsymbol{\mu}} \tilde{J}_N(\boldsymbol{\mu}) = 	\nabla_{\boldsymbol{\mu}} J_N({\bf{u}}_N,\boldsymbol{\mu}) - \overline{\bf{u}}_N  \left(	D_{\boldsymbol{\mu}}\mathbf{G}_N({\bf{u}}_N,\boldsymbol{\mu})\right).
\end{equation}
Note that the evaluation of the above derivatives with respect to $\boldsymbol{\mu}$ is inexpensive assuming that the operators depend affinely on the parameters.

Considering the compliance case, the solution of the adjoint problem can be directly obtained from the reduced solution ${\mathbf{u}}_N$. Thus, by inserting Equations~\eqref{eq32} and~\eqref{eq34} into Equation~\eqref{eq40}, the adjoint problem for the compliance case reads:
\begin{equation}\label{eq38}
	{\bf{K}}_N(\boldsymbol{\mu})\overline{\bf{u}}_N(\boldsymbol{\mu}) = \frac{1}{2} {\bf{f}}_N(\boldsymbol{\mu}) \implies \overline{\bf{u}}_N(\boldsymbol{\mu}) = \frac{1}{2} {\bf{u}}_N(\boldsymbol{\mu}).
\end{equation}
Finally, inserting Equation~\eqref{eq38} into~\eqref{eq41}, and expanding the derivative $D_{\boldsymbol{\mu}}\mathbf{G}_N({\bf{u}}_N,\boldsymbol{\mu})$, the gradient reads:
\begin{align}\label{eq36}
 	\nabla_{\boldsymbol{\mu}} \tilde{J}_N(\boldsymbol{\mu})	= {\bf{u}}_N(\boldsymbol{\mu})\cdot \left(\frac{\partial {\bf{f}}_N(\boldsymbol{\mu})}{\partial \boldsymbol\mu}
  - \frac{1}{2}\frac{\partial {\bf{K}}_N(\boldsymbol{\mu})}{\partial \boldsymbol\mu} {\bf{u}}_N(\boldsymbol{\mu})\right)
\end{align}
Under the assumption of affine parametric dependence in~\eqref{eq56} and considering~\eqref{eq47a} and~\eqref{eq52}, we rewrite the gradient for $\boldsymbol{\mu} \in \mathcal{P}_k$, as
\begin{align}\label{eq37}
		\nabla_{\boldsymbol{\mu}} \tilde{J}_N(\boldsymbol{\mu}) &= {\bf{u}}_N(\boldsymbol{\mu})\cdot   \left[\sum_{q=1}^{Q_f^k}
	{\frac{\partial {\theta_{f,q}^k(\boldsymbol{\mu})}}{\partial \boldsymbol\mu}}
	{\bf{V}}_k^{\top}\widehat{{\bf{f}}}_q^k - \frac{1}{2}\left(\sum_{q=1}^{Q_a^k}
	{\frac{\partial {\theta_{a,q}^k(\boldsymbol{\mu})}}{\partial\boldsymbol \mu}}
	{\bf{V}}_k^{\top}\widehat{{\bf{K}}}_q^k{\bf{V}}_k\right){\bf{u}}_N(\boldsymbol{\mu}) \right].
\end{align}

The derivatives $\partial {\theta_{a,q}^k(\boldsymbol{\mu})}/ \partial \boldsymbol\mu$ are simple and inexpensive to evaluate by just differentiating the expression~\eqref{eq66}.
The procedure for $\partial {\theta_{f,q}^k(\boldsymbol{\mu})}/ \partial \boldsymbol\mu$ is identical.
\section{Numerical results}\label{sec7}
In this section we present some numerical experiments for the Kirchhoff-Love shell problem to assess the capabilities of the presented ROM framework for parameterized trimmed and multi-patch geometries.
In what follows we apply the two presented strategies to enforce interface coupling conditions in a weak sense, namely the projected super-penalty (Section~\ref{sec4a}) and Nitsche's method (Section~\ref{sec4b}).
The numerical experiments are carried out using the open-source Octave/Matlab isogeometric package GeoPDEs \cite{Vazquez2016} in combination with the re-parameterization tool for integration of trimmed geometries presented in our previous works \cite{Antolin2019,Wei2021}, while the ROM construction exploits the open-source library redbKIT \cite{redbKIT}. We remark that the stiffness matrix is preconditioned with a diagonal scaling to avoid large conditions numbers and loss of accuracy due to small trimmed elements as discussed in \cite{Deprenter2017}.

\subsection{Scordelis-Lo roof with holes}
\begin{figure}[!h]
	\centering
	\includegraphics[width=1.0\textwidth]{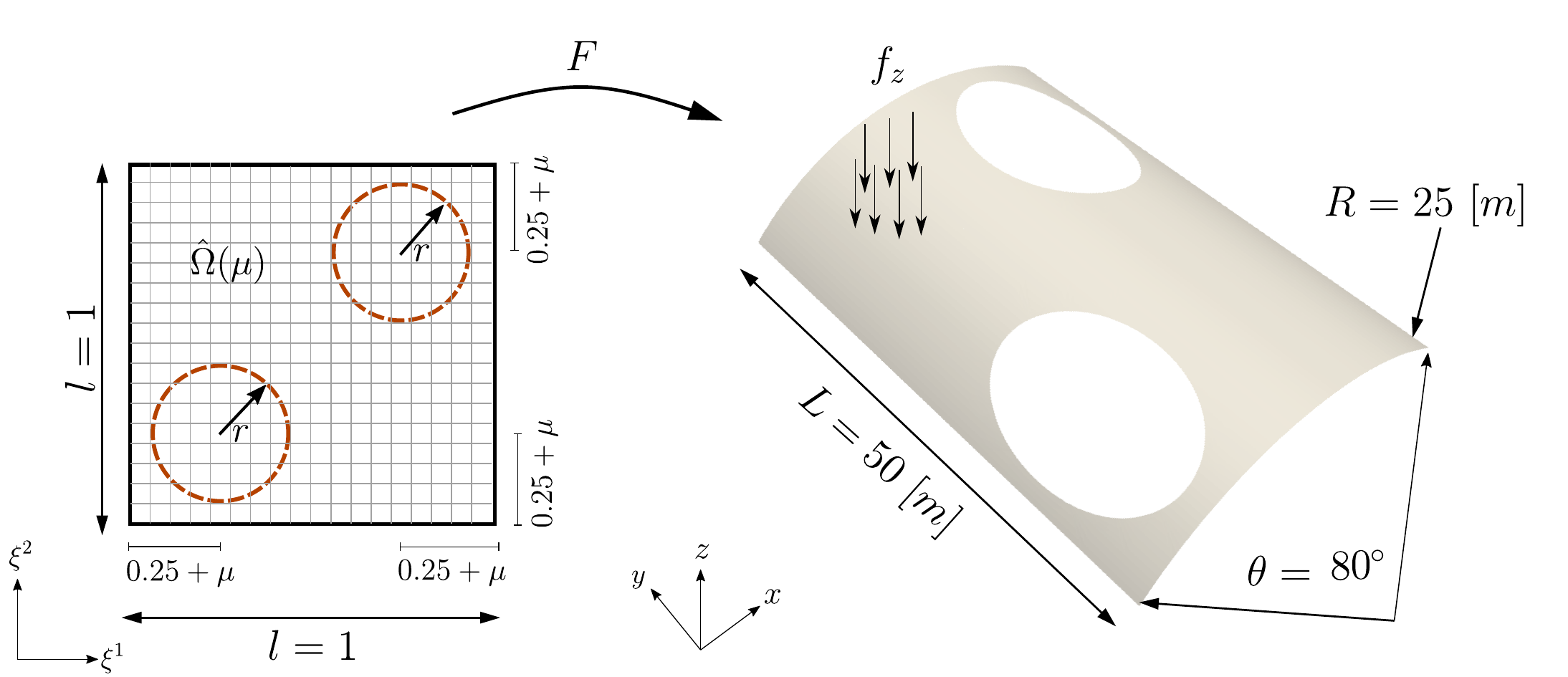} \\
	\caption{Example 7.1: Geometry and parameterization of the Scordelis-Lo roof with holes.}\label{fig:trimmed_scordelis}
\end{figure}

The first numerical example is a single-patch, trimmed variant of the Scordelis-Lo roof. The geometry and material properties are adopted from the well known benchmark, see, e.g., \cite{Belytschko1985} for more details. Thus, the Young's modulus is $E=4.32 \cdot 10^8$ Pa, the Poisson's ratio $\nu=0$, and the thickness $t=0.25$~m. The shell structure is subjected to self-weight with a vertical loading of $f_z = -90$~N/m$^2$ as depicted in Fig.~\ref{fig:trimmed_scordelis}. Rigid diaphragm boundary conditions are imposed on the curved ends of the shell, that is, we fix the displacements in the $xz$-plane. In this example, we cut out two circular holes in the parametric domain as shown in Fig.~\ref{fig:trimmed_scordelis}. The radius of the holes in the parametric domain is fixed as $r = 0.2$. Their location is parameterized, where $\mu \in \mathcal{P}=[0, 0.1]$ is a geometric parameter representing the location of the center of each circle that moves along the diagonal of the unit square. The geometry of the shell structure in the physical space is then obtained by an additional mapping as mentioned in Section~\ref{sec2}. The shell is discretized with cubic $C^2$-continuous B-splines and the dimension of the non-trimmed space is $\mathcal{N}_{h,0} = 1083$.

\begin{figure}[!h]
	\begin{subfigure}[b]{0.49\textwidth}
		\centering
		\includegraphics[width=\textwidth]{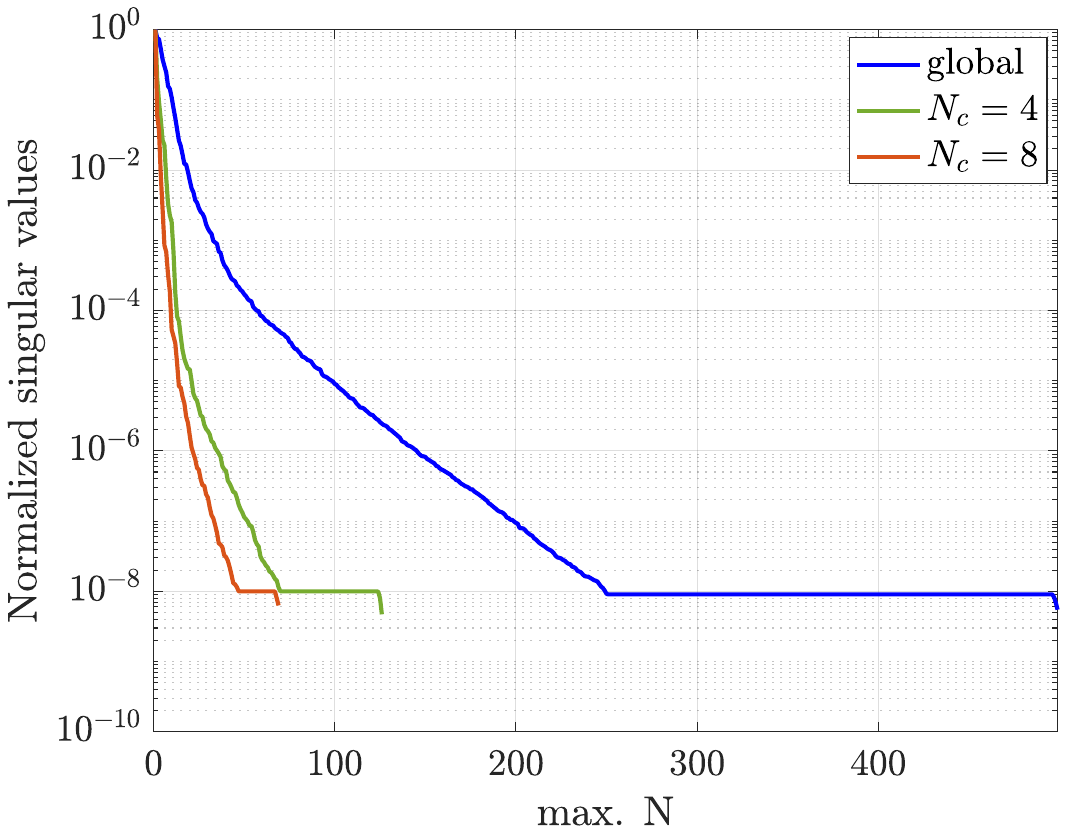}
		\caption{Singular values decay}
		\label{fig:svd_1}
	\end{subfigure}
	\begin{subfigure}[b]{0.49\textwidth}
		\centering
		\includegraphics[width=\textwidth]{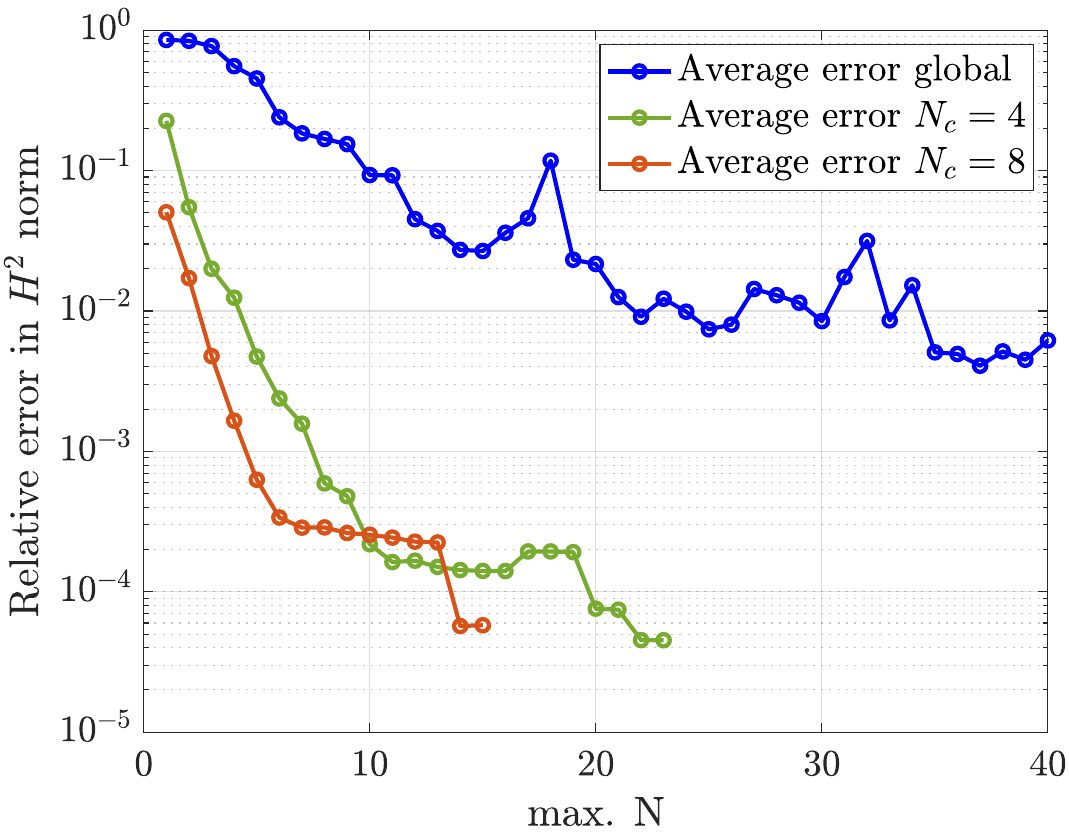}
		\caption{Error analysis}
		\label{fig:error_1}
	\end{subfigure}
	\caption{Example 7.1: Singular values decay and relative error in $H^2$ norm vs.\ maximum number of reduced basis functions $N$ over all the clusters, for different numbers of clusters.}
	\label{fig:svd_error_1}
\end{figure}

\begin{figure}[!h]
	\centering
	\includegraphics[width=0.6\textwidth]{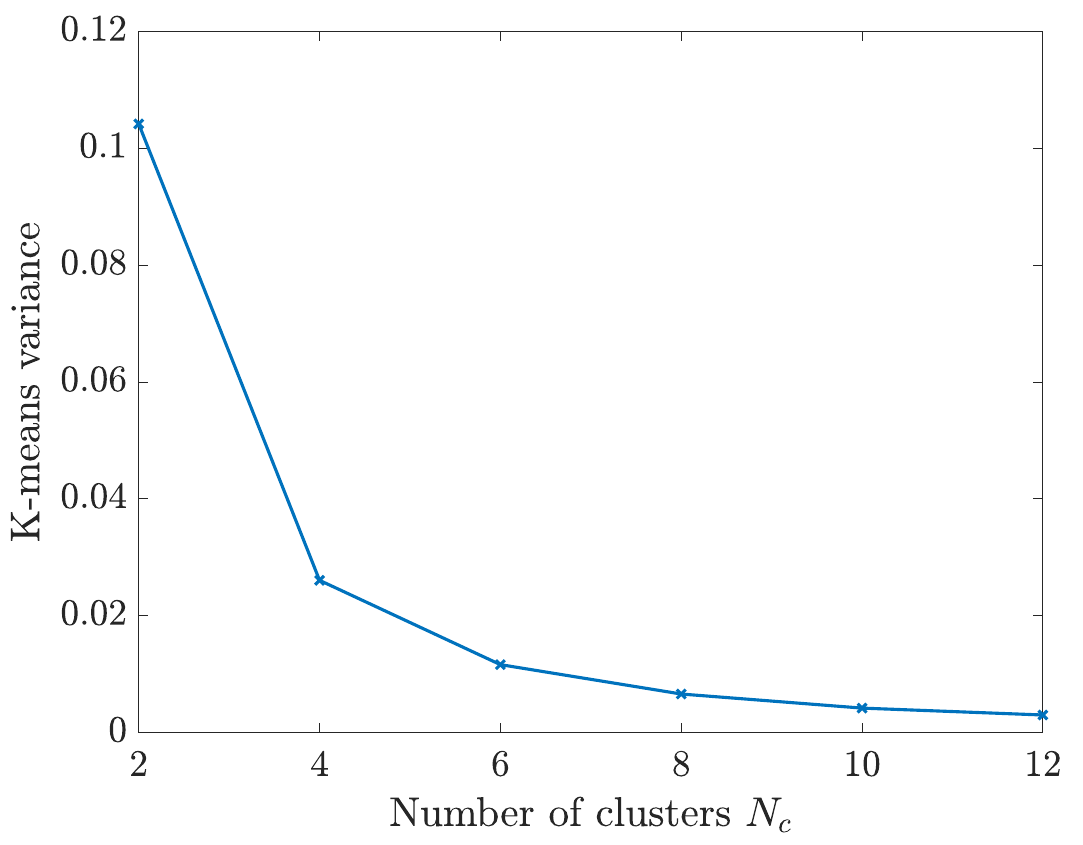} \\
	\caption{Example 7.1: K-means variance over number of clusters $N_c$.}\label{fig:kmeans_scordelis}
\end{figure}

\begin{figure}[h]
	\centering
	\includegraphics[width=0.85\textwidth]{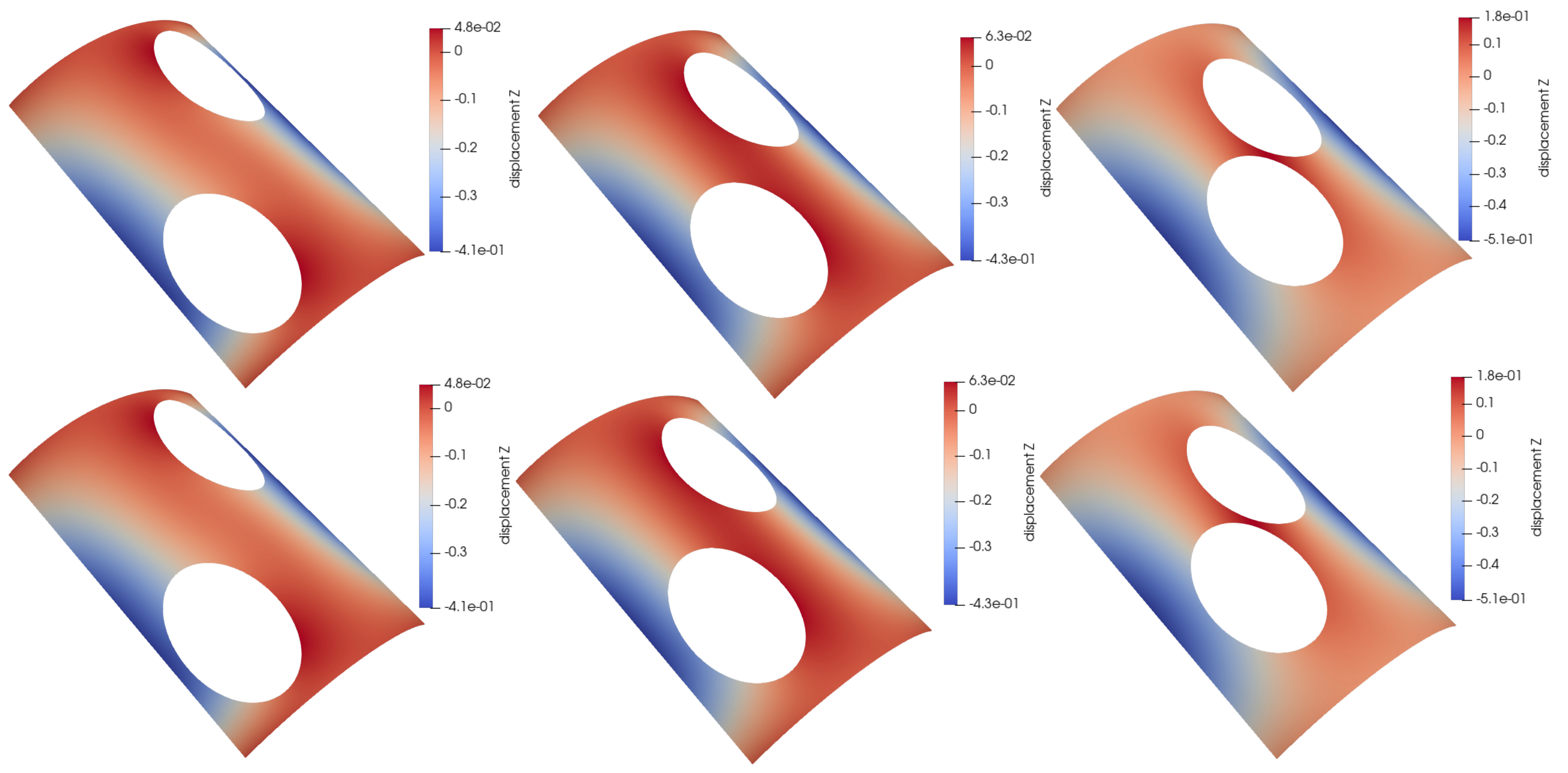} \\
	\caption{Example 7.1: Vertical displacement solutions computed with the FOM (top) and local ROM (bottom) with $N_c = 8$ clusters for three parameter values $\mu= \lbrace0.0, 0.05, 0.1\rbrace$.}\label{fig:scordelis1_plot}
\end{figure}

Let us now construct a ROM with the local reduced basis method presented in Section~\ref{sec5}. For this purpose, we use a training sample of dimension $N_s = 500$ obtained by Latin Hypercube sampling \cite{McKay1979} for the construction of the reduced basis. Note that this sampling is also suitable for more general problems with higher dimension of the parameter space \cite{QMN_RBspringer}. Fig.~\ref{fig:svd_1} depicts the decay of the singular values of the POD for different numbers of clusters. We remark that the POD basis is constructed such that it minimizes the squared projection error with respect to the matrix norm, that represents the algebraic counterpart of the $H^2$ norm. Moreover, the number of clusters for the DEIM approximations is fixed to $16$ for all computations. It can be observed that the decay is more rapid with localized ROMs. Now we employ a test sample of dimension $N_t = 30$ with a uniform random distribution to perform the error analysis whose results are presented in Fig.~\ref{fig:error_1}. We observe that increasing the number of clusters further reduces the dimension of the reduced basis, while an accuracy of $10^{-5}$ is achieved in the $H^2$ norm. The optimal number of clusters is selected as $N_c=8$  based on the k-means variance in  Fig.~\ref{fig:kmeans_scordelis}. Note that the computation of the optimal number of clusters is performed offline and requires 0.21 s for the problem at hand. Fig.~\ref{fig:scordelis1_plot} depicts the vertical displacement solutions, where the local ROM with $N_c = 8$ clusters is compared to the FOM. The online CPU time is for the ROM is $0.12$ s including the closest cluster iterations, whereas the assembly and solution with the FOM requires $2.32$ s. This results in a speedup of $19 \times$. Note that the trimming operation for computing the FOM or vizualizing the results requires  additionally 0.854 s. We remark that the offline time to compute the snapshots for the POD basis is 3 minutes exploiting the affine approximation  \eqref{eq56}. The offline time to compute the snapshot matrices and vectors for the DEIM approximations is 1.75 hours given a training sample of 2000 snapshots. We remark that this  includes the trimming operation for each snapshot and can be further optimized with the use of parallelization.

\begin{figure}[h]
	\begin{subfigure}[b]{0.49\textwidth}
		\centering
		\includegraphics[width=\textwidth]{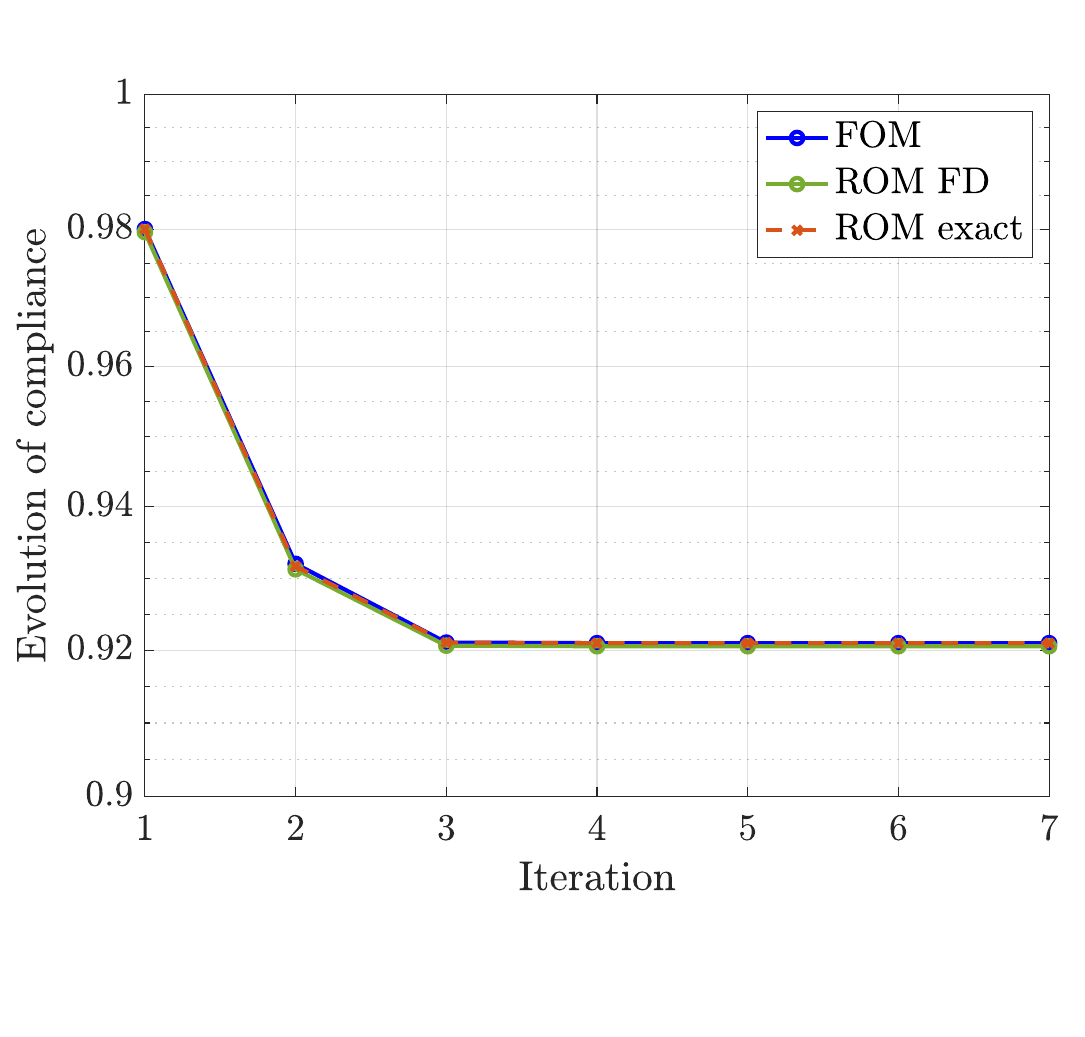}
		\caption{Evolution of the compliance}
		\label{fig:compliance_1}
	\end{subfigure}
	\begin{subfigure}[b]{0.49\textwidth}
		\centering
		\includegraphics[width=\textwidth]{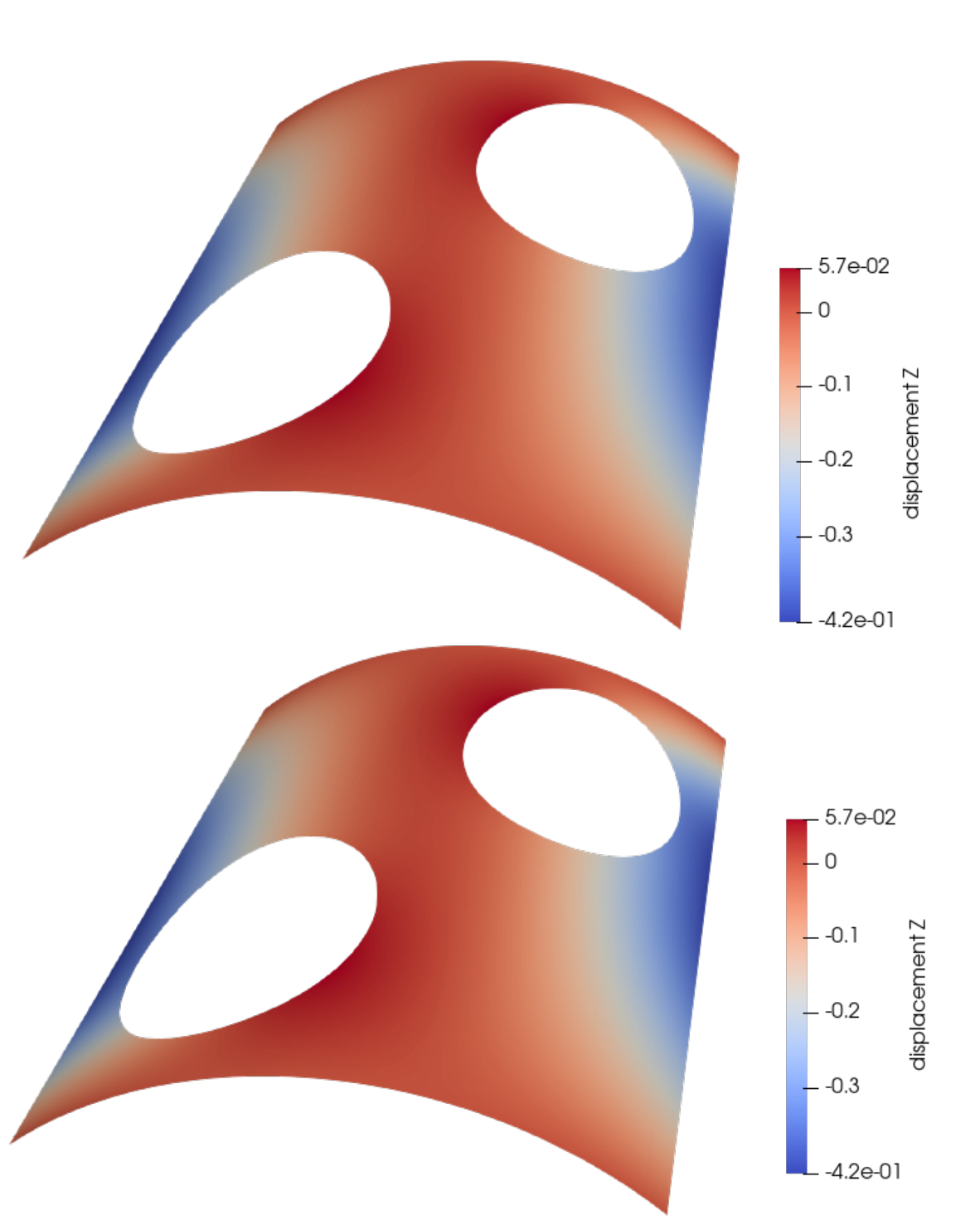}
		\caption{Final optimal shape}
		\label{fig:opti_1}
	\end{subfigure}
	\caption{Example 7.1: Optimization results depicting: (a) the evolution of the relative compliance during the optimization comparing the FOM, the ROM with approximate sensitivities, and the ROM with exact sensitivities; and (b) the vertical displacement for the final optimal shape with $\mu = 0.0319$ using the FOM (top) and the local ROM (bottom).}
	\label{fig:scordelis1_opti}
\end{figure}

Moreover, we perform optimization using both the local ROM and the FOM.
The upper and lower bounds for the design variable are defined by the parameterization at hand, that is we seek the optimal shape within the bounds $0 \le \mu \le 0.1$ during the optimization.
For this problem we do not consider any volume constraints. First, we solve the optimization problem with the FOM using a finite difference scheme to compute sensitivities. The solution requires 7 iterations and 23 function evaluations in total. Then, we employ the local ROM with $N_c=8$ clusters for the optimization. The ROM with approximate sensitivities requires 17 function evaluations versus the ROM with exact sensitivities with 15 function evaluations. Note that a forward finite difference scheme is employed for the computation of the approximate sensitivities during the optimization. The solution is slightly faster for the ROM with exact sensitivities and comprises $76.1$ ms versus $79.5$ ms for the solution with finite difference approximation of the gradient. In all cases the optimal shape is obtained for $\mu = 0.0319$ after 7 iterations. Fig.~\ref{fig:scordelis1_opti} depicts the optimization results. The optimization history is shown in Fig.~\ref{fig:compliance_1} by depicting the evolution of the relative compliance during the optimization. For the optimal solution, the compliance is reduced $8$\% respect to the initial configuration (i.e., for $\mu =0$). Finally, Fig.~\ref{fig:opti_1} shows the vertical displacement solution for the final optimal shape, comparing ROM and FOM solutions.

\subsection{Multi-patch simple geometries}
In this section we will assess the capabilities of the ROM framework for multi-patch geometries. For this purpose we employ two simple geometries, namely the multi-patch Scordelis-Lo roof and two non-conforming planar patches. The coupling of patches is achieved in both cases using the projected super-penalty method discussed in Section~\ref{sec4a}.

\subsubsection{Non-trimmed multi-patch Scordelis-Lo roof}
\begin{figure}[!h]
	\begin{subfigure}[b]{0.3\textwidth}
		\centering
		\includegraphics[width=\textwidth]{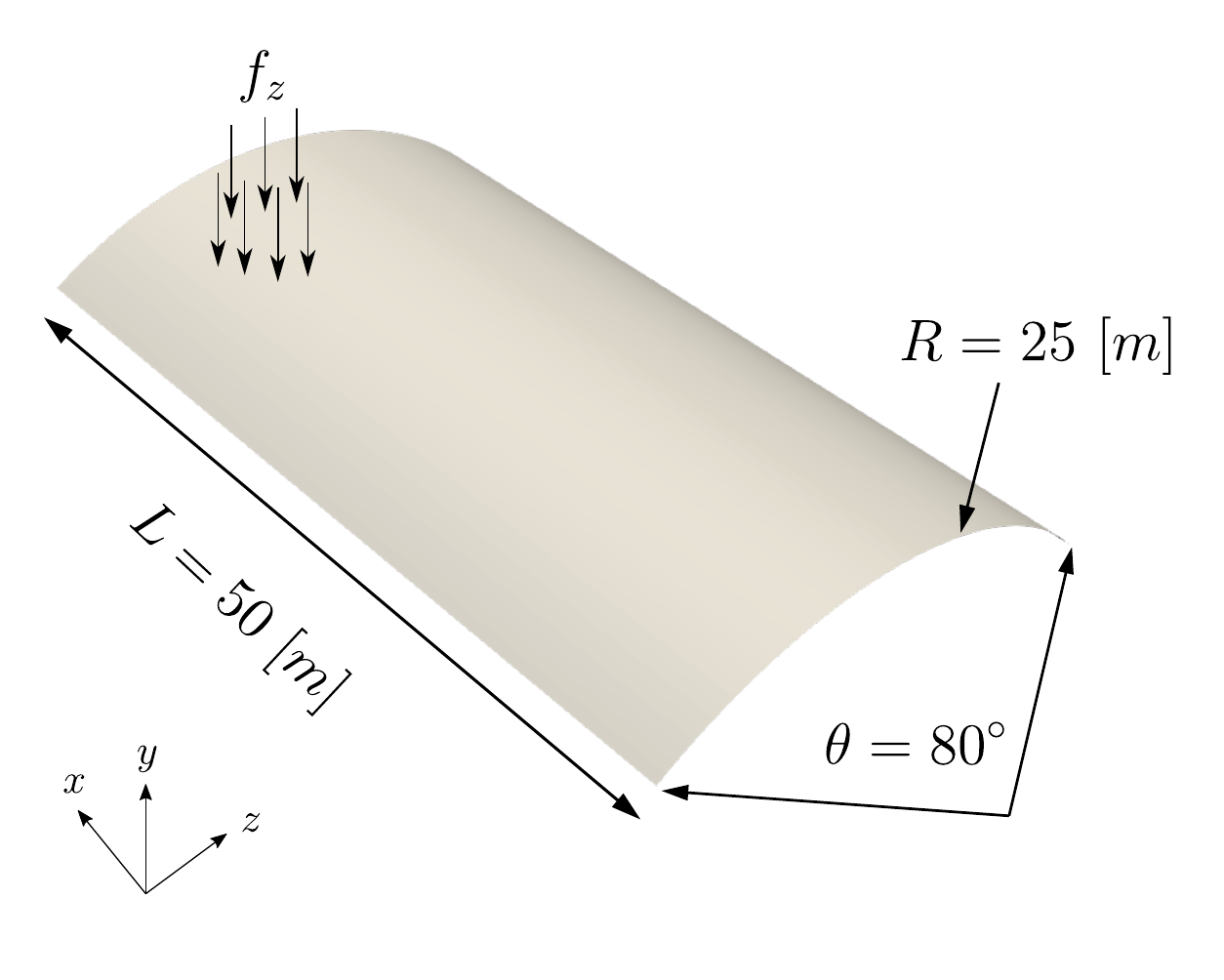}
		\caption{Geometry setup.}
		\label{fig:scordelis_geometry}
	\end{subfigure}
	\begin{subfigure}[b]{0.3\textwidth}
		\centering
		\includegraphics[width=\textwidth]{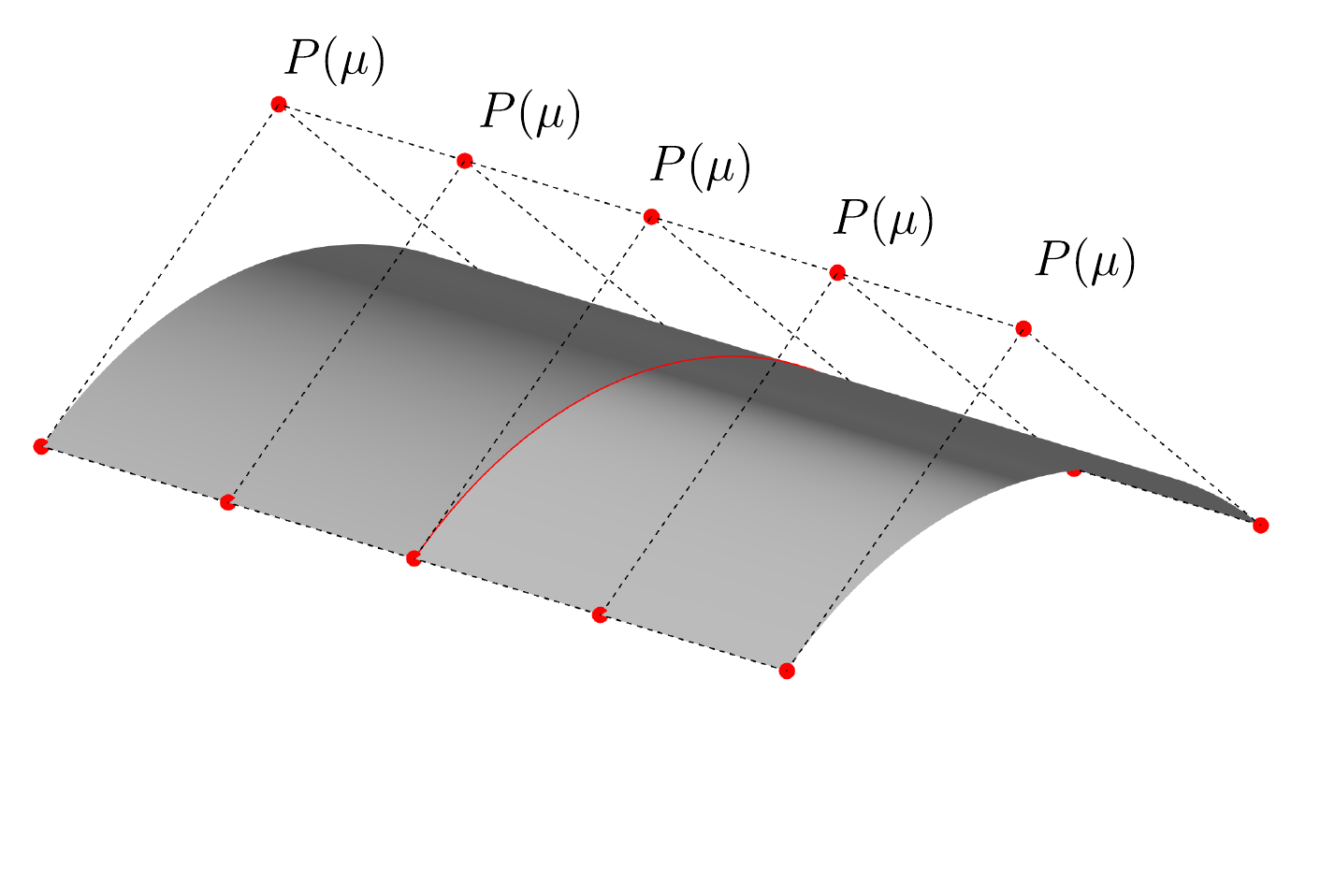}
		\caption{$\mu=0$}
		\label{fig:scordelis_mu1}
	\end{subfigure}
	\begin{subfigure}[b]{0.3\textwidth}
	\centering
	\includegraphics[width=\textwidth]{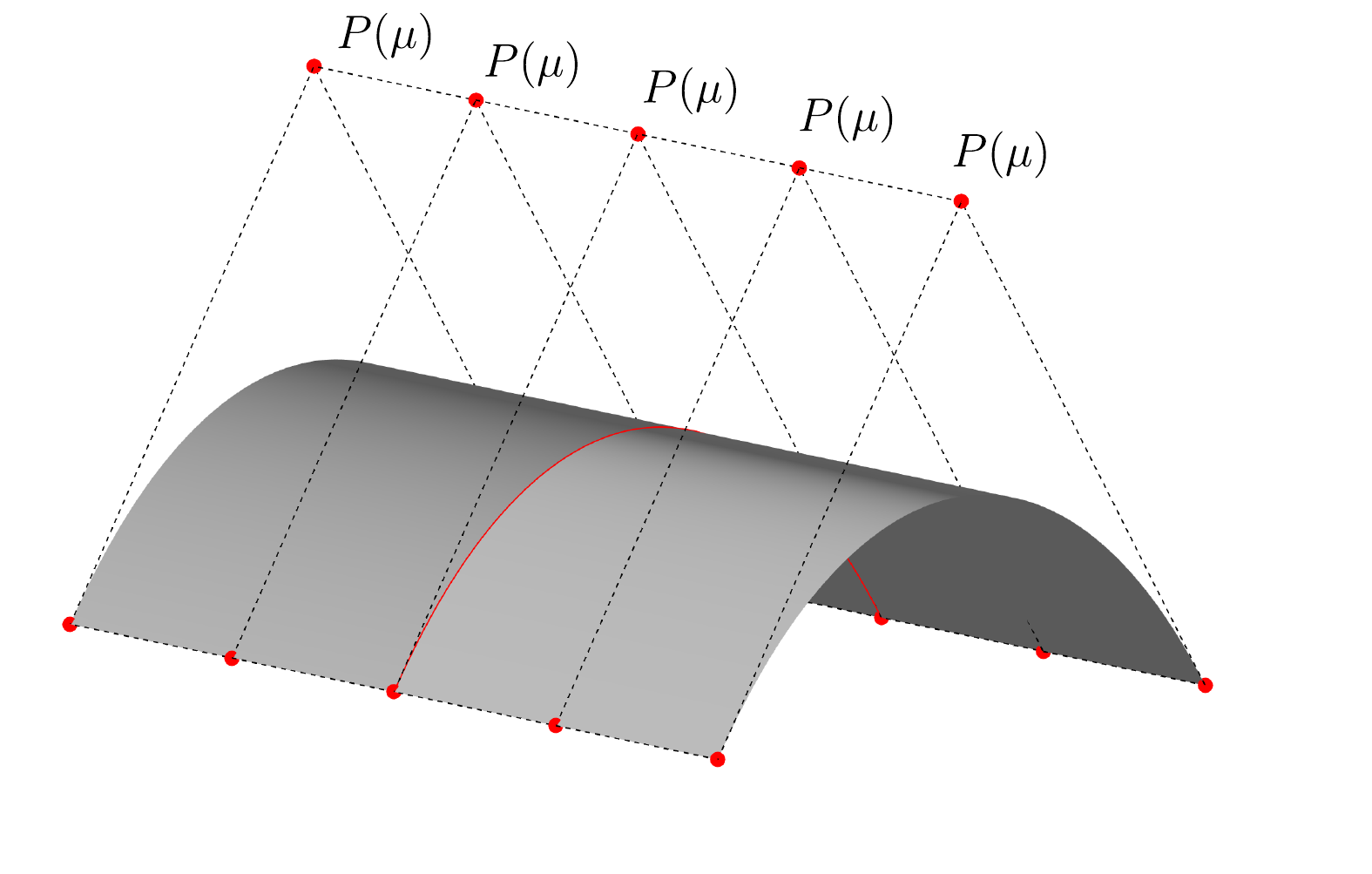}
	\caption{$\mu=10$}
	\label{fig:scordelis_mu2}
\end{subfigure}
	\caption{Example 7.2.1.: Problem setup and geometrical parameterization for different values of $\mu$ for the multi-patch Scordelis-Lo roof.}
	\label{fig:scordelis_setup}
\end{figure}

This example is intended to test the ROM framework on a non-trimmed, multi-patch setting. For this purpose, we employ again the  Scordelis-Lo roof and split the geometry into two subdomains. The geometric setup is depicted in Fig.~\ref{fig:scordelis_setup}. We employ the projected super-penalty method to enforce interface coupling conditions. The material parameters, loading, and boundary conditions are the same as in the previous example. The shell structure is modeled with two conforming subdomains and the parameterization of the multi-patch design is shown in Figs.~\ref{fig:scordelis_mu1}-\ref{fig:scordelis_mu2} for a coarse geometry. The common interface is depicted  in red color. Note that the geometry is parameterized by moving the depicted control points in the vertical direction, where the geometric parameter $\mu \in \mathcal{P} = [0,10]$ defines their position in the y-direction prescribing the curvature of the shell structure. We remark the $\mu=0$ corresponds to the original coordinates of the Scordelis-Lo roof benchmark. In Figs.~\ref{fig:scordelis_mu1}-\ref{fig:scordelis_mu2} we denote for simplicity $P(\mu) = P(x,y+\mu,z)$. The parameters are first defined on a coarse mesh and then the geometry is further refined for the analysis as discussed in \cite{Chasapi2022}. The geometry is discretized with quadratic $C^1$-continuous B-splines for the analysis resulting in $\mathcal{N}_{h}=1944$ degrees of freedom.

\begin{figure}[h]
	\centering
	\includegraphics[width=0.8\textwidth]{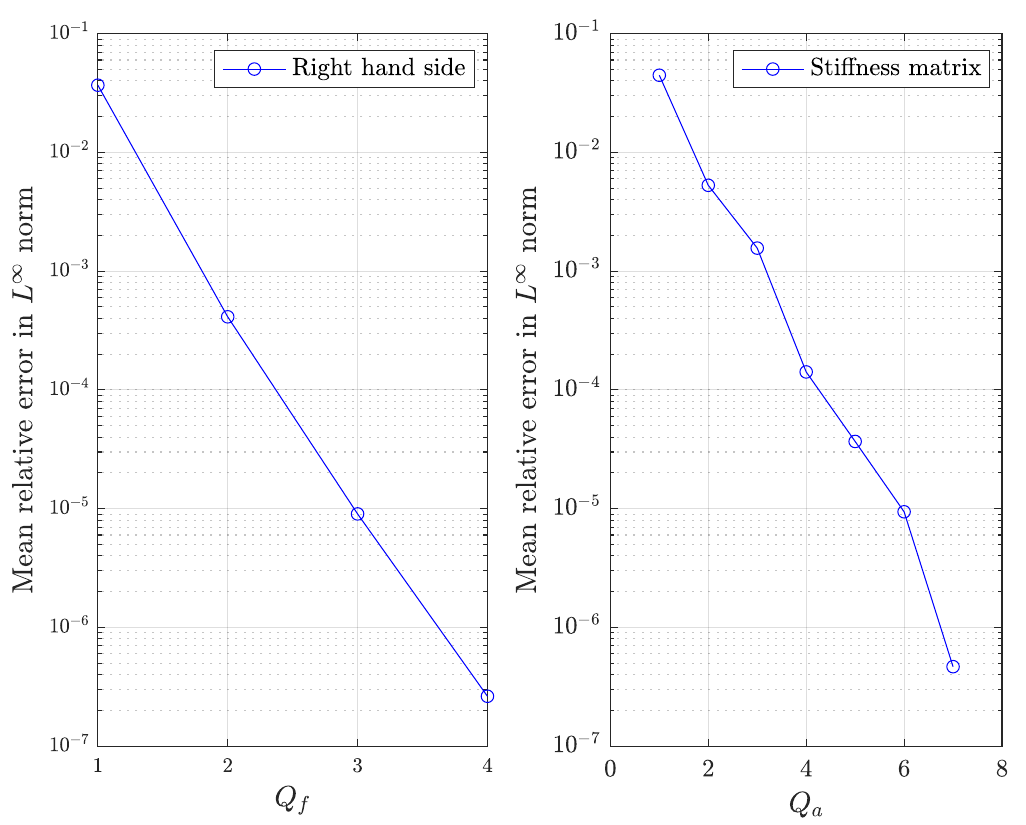} \\
	\caption{Example 7.2.1: Error decay of DEIM approximations in $L^{\infty}$-error norm for right-hand side vector and stiffness matrix using a POD tolerance of $\varepsilon_{\text{POD}}= 10^{-7}$.}\label{fig:scordelis_deim}
\end{figure}

\begin{figure}[h]
	\begin{subfigure}[b]{0.49\textwidth}
		\centering
		\includegraphics[width=\textwidth]{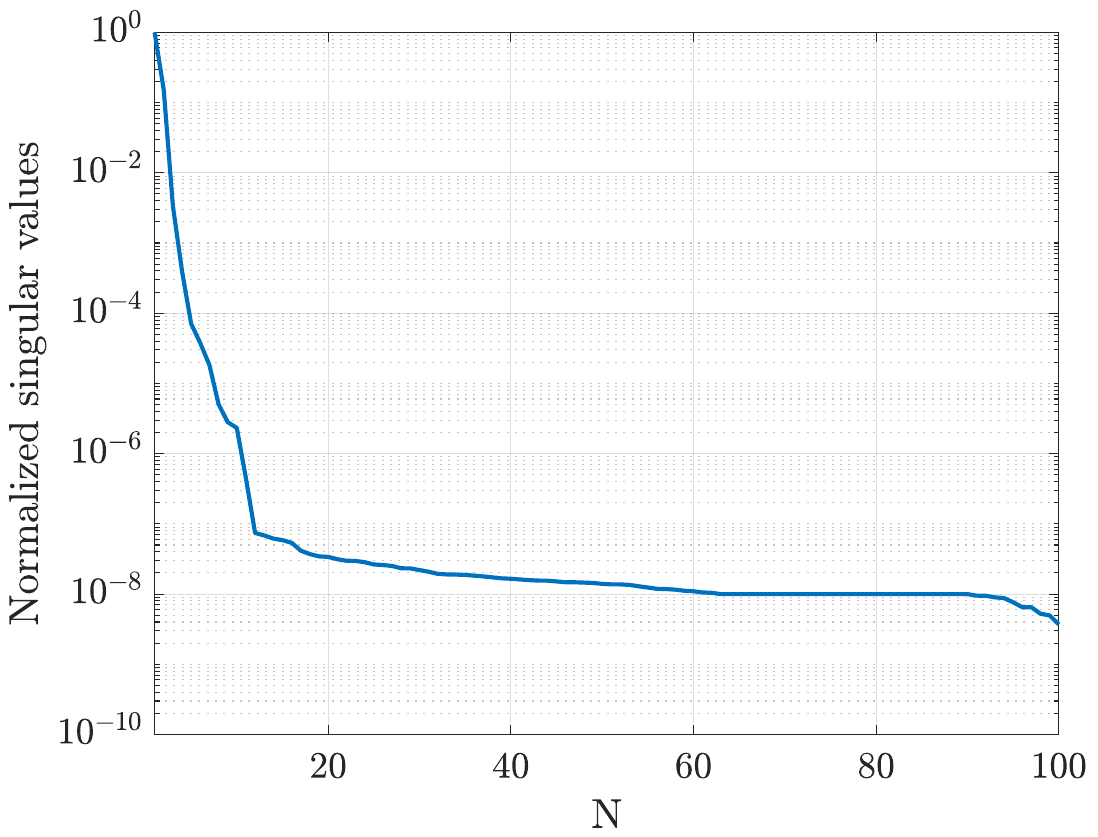}
		\caption{Singular values decay}
		\label{fig:svd_2}
	\end{subfigure}
	\begin{subfigure}[b]{0.49\textwidth}
		\centering
		\includegraphics[width=\textwidth]{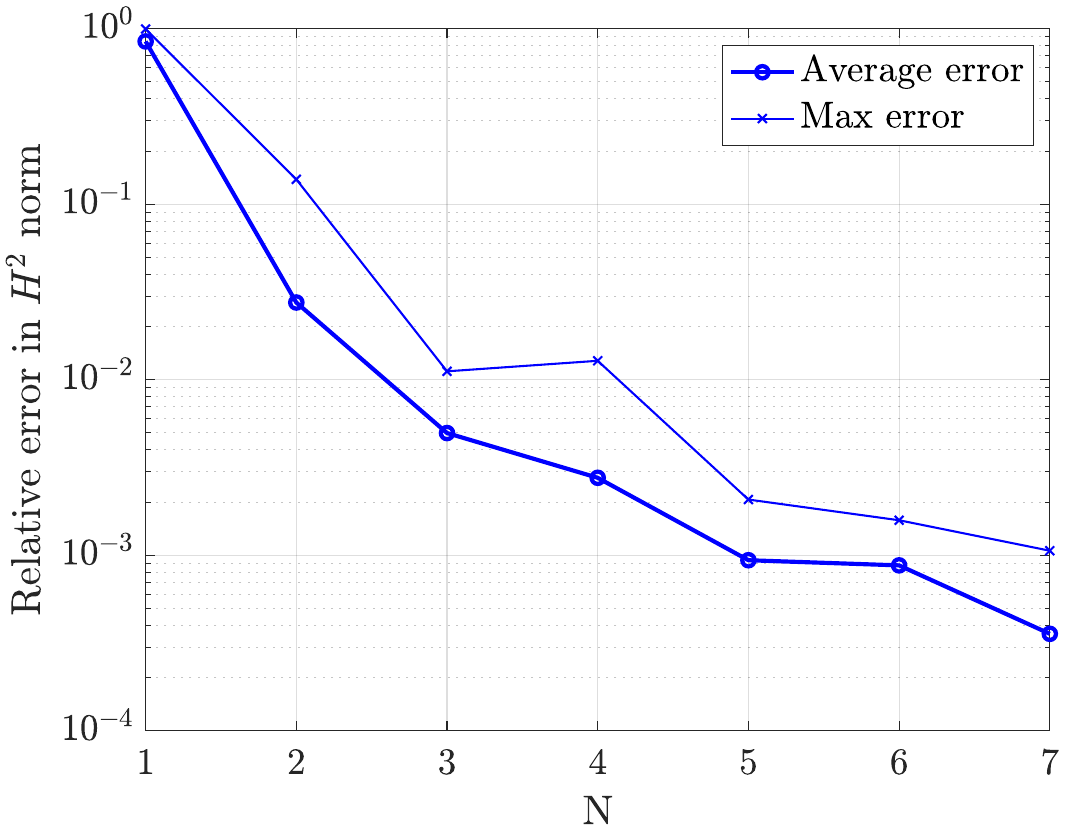}
		\caption{Error analysis}
		\label{fig:error_2}
	\end{subfigure}
	\caption{Example 7.2.1: Singular values decay and relative error in $H^2$ norm vs.\ number of reduced basis functions N.}
	\label{fig:svd_error_2}
\end{figure}

\begin{figure}[h]
	\centering
	\includegraphics[width=1.0\textwidth]{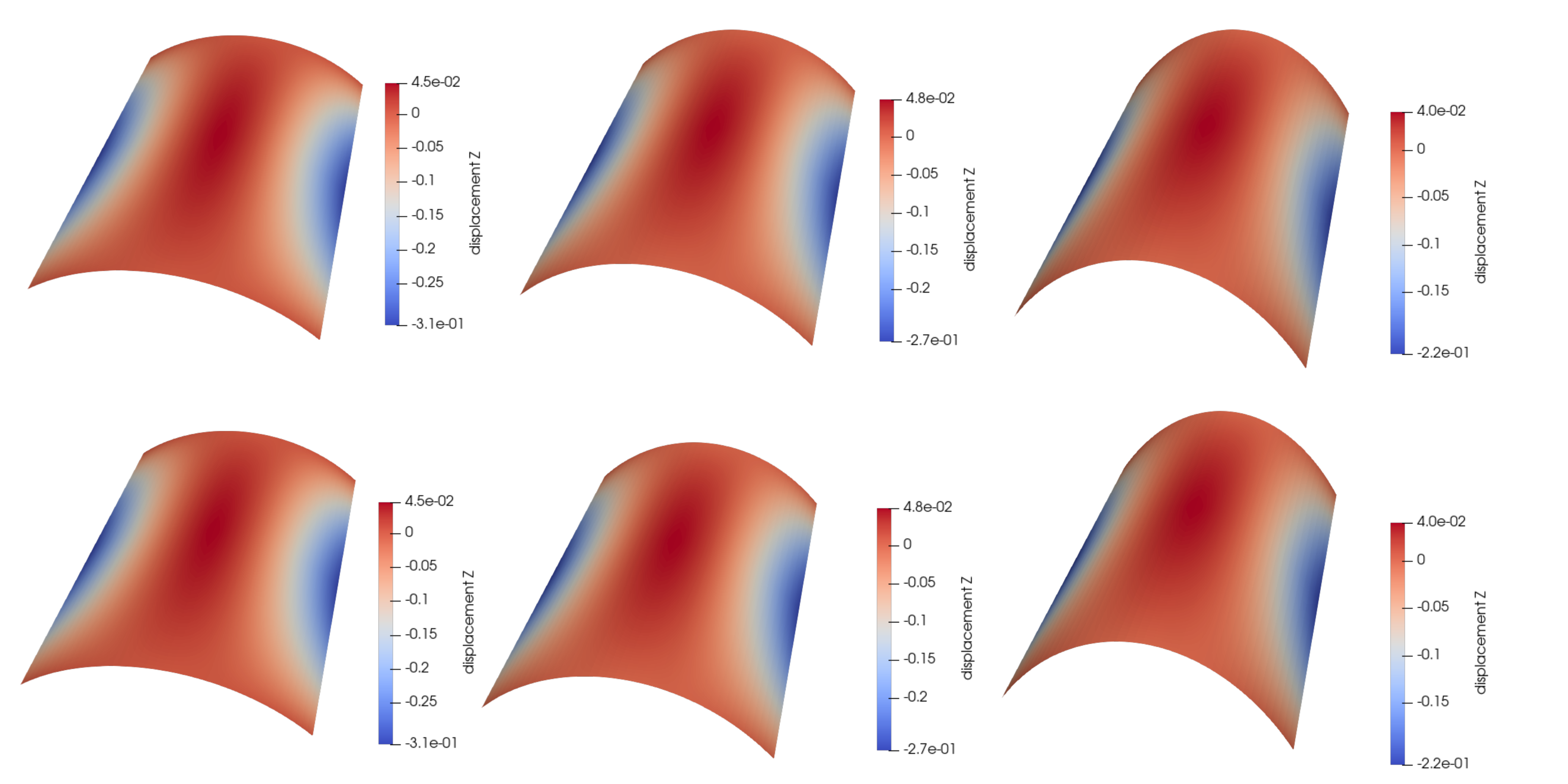} \\
	\caption{Example 7.2.1: Vertical displacement solutions computed with the FOM (top) and ROM (bottom) for three parameter values $\mu= \lbrace0.0, 3.5, 10.0\rbrace$.}\label{fig:scordelis2_plot}
\end{figure}

Now let us construct a ROM for the multi-patch geometry. For this purpose, we employ a training sample of dimension $N_s=500$ obtained by Latin Hypercube sampling. First, the affine approximations are constructed with the DEIM. Fig.~\ref{fig:scordelis_deim} depicts the error of the DEIM approximations in the $L^{\infty}$ norm for both the right-hand side vector and the stiffness matrix. It can be observed that the error already decays rapidly with one global approximation. The error reaches an accuracy of $10^{-7}$ with $Q_a=8$ and $Q_f=5$ affine terms for the stiffness matrix and the right-hand side vector, respectively. As a further step, a reduced basis is constructed with the POD. Fig.~\ref{fig:svd_error_2} shows the decay of the singular values and the relative error of the reduced solution in the $H^2$ norm. The error analysis is performed using a test sample of dimensions $N_t=30$ and a uniform random distribution. Similarly to the DEIM approximations, the decay is already rapid by constructing only one global reduced basis space and the error reaches an accuracy of $10^{-7}$ for a reduced basis space of dimension $N=7$. In Fig.~\ref{fig:scordelis2_plot} we compare the vertical displacement solution of the ROM with the solution of the full order model for three different parameter values. 

\begin{figure}[h]
	\centering
	\includegraphics[width=0.9\textwidth]{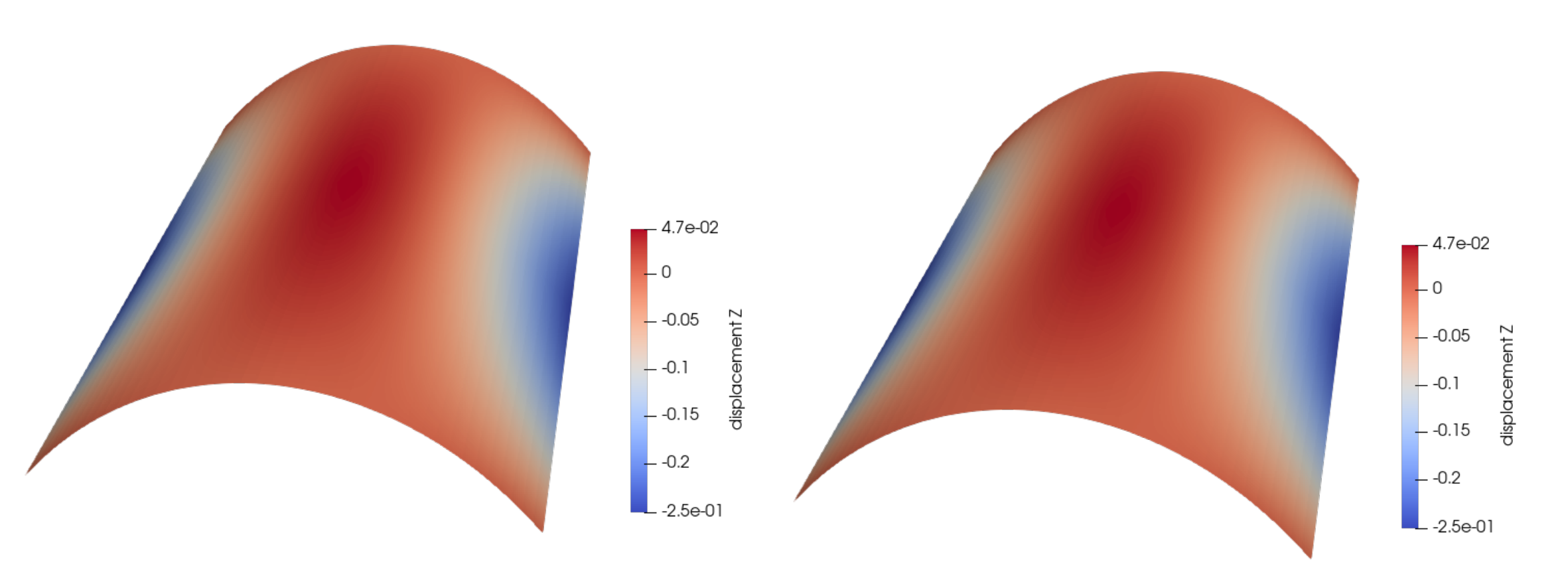} \\
	\caption{Example 7.2.1: Vertical displacements for the final optimal shape with $\mu = 5.3127$ using the FOM (left) and the ROM (right).}\label{fig:scordelis2_opt}
\end{figure}

Moreover, we perform optimization considering a volume constraint. The initial volume is set as $V_0= 1.9 \cdot 10^3$ and we seek the optimal shape for the design variable bounded by $0 \le \mu \le 10$. The optimization is completed after 4 iterations and 14 function evaluations. The ROM-based optimization is computed in $67$~ms, while the optimization with the high fidelity solution is solved in $518$~ms. This implies a speedup of $7.76\times$. In both cases the exact gradients are computed at every optimization step. At the final iteration, the compliance is reduced by $16.45\%$ compared to the initial configuration. Fig.~\ref{fig:scordelis2_opt} depicts the vertical displacement solutions for the optimal shape obtained for $\mu = 5.3127$ with both the ROM and full order solution, while Table~\ref{tab:scordelis_times} summarizes the main results and obtained computation times.  We conclude that this verifies the suitability of the ROM framework for multi-patch geometries coupled along non-trimmed interfaces.

\begin{table}[h]
	\caption{Example 7.2.1: Number of basis functions  and computational cost.}\label{tab:scordelis_times}
      \centering
	\begin{tabular*}{0.65\textwidth}{@{\extracolsep\fill}lc}
		\toprule
		$Q_a$    & 8  \\
		$Q_f$   & 5 \\
		$N$   & 7  \\
		Online CPU time [ms]   & 9.55 \\
		Solution speedup & $14.65\times$ \\
		ROM-based optimization time [ms] & 67 \\
		FOM-based optimization time [ms] & 518 \\
		Optimization speedup & $7.76\times$ \\
		\botrule
	\end{tabular*}
\end{table}

\newpage

\subsubsection{Trimmed non-conforming planar patches }
\begin{figure}[h]
	\begin{subfigure}[b]{0.49\textwidth}
		\centering
		\includegraphics[width=\textwidth]{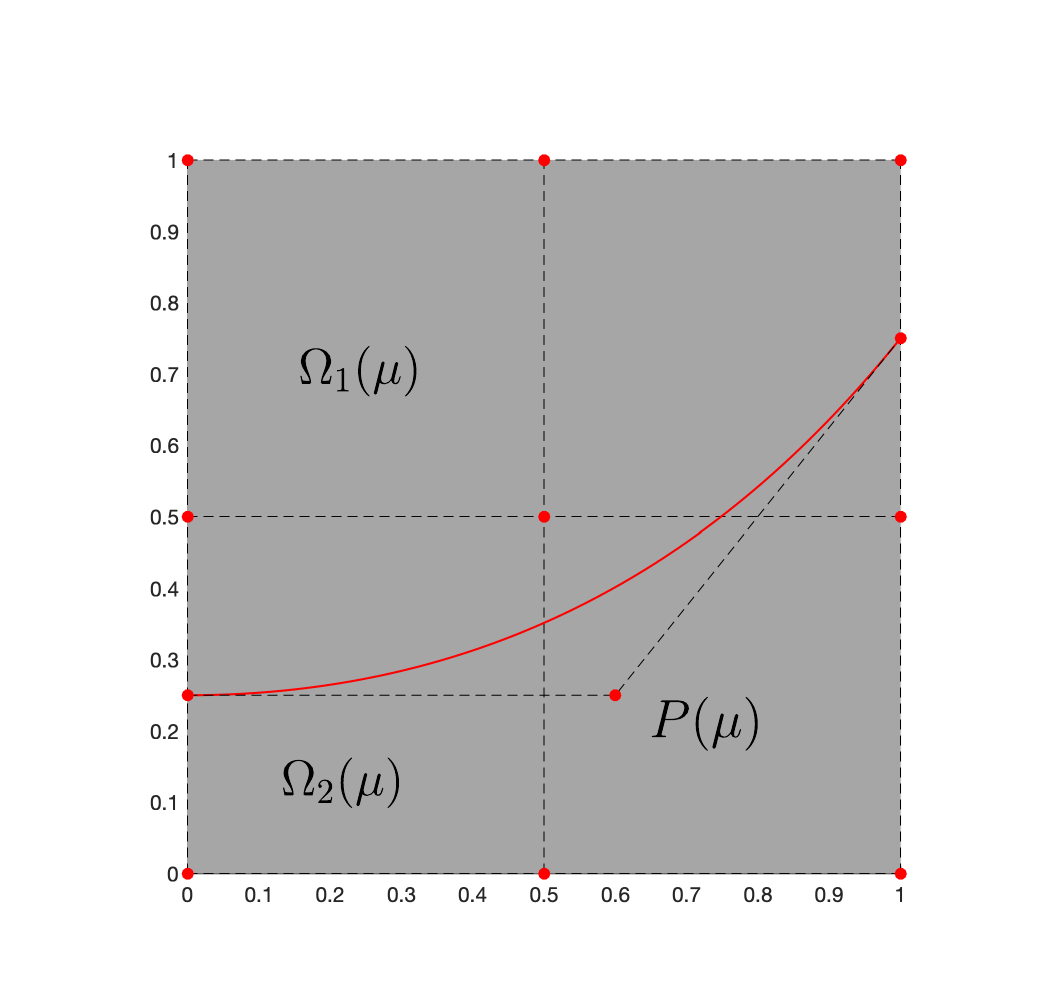}
		\caption{$\mu=0.25$}
		\label{fig:plane_mu1}
	\end{subfigure}
	\begin{subfigure}[b]{0.49\textwidth}
		\centering
		\includegraphics[width=\textwidth]{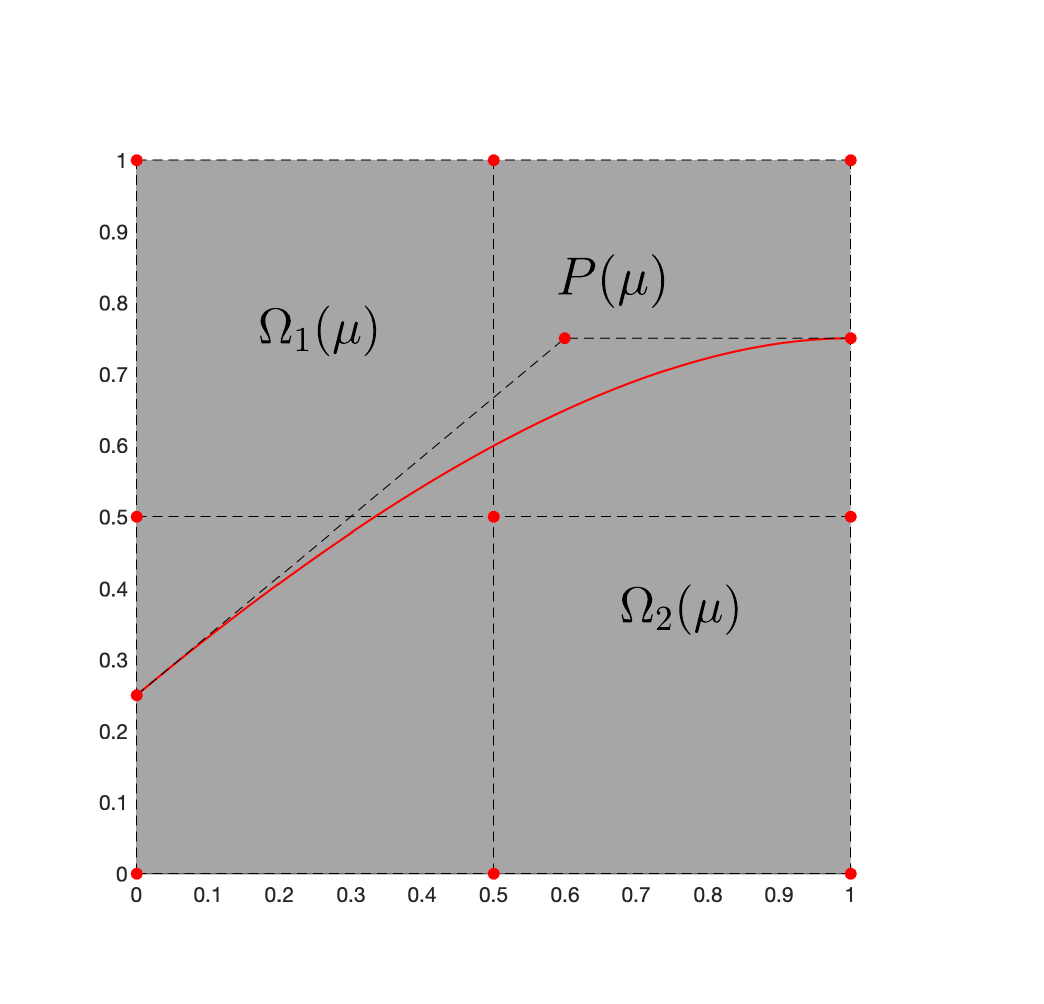}
		\caption{$\mu=0.75$}
		\label{fig:plane_mu2}
	\end{subfigure}
	\caption{Example 7.2.2: Geometry setup and parameterization of the two planar trimmed patches for different values of $\mu$. The trimming interface is a quadratic spline curve and denoted in red color.}
	\label{fig:plane_setup}
\end{figure}

In this example we aim to assess the capabilities of the ROM for trimmed multi-patch geometries with both conforming and non-conforming discretizations. In particular, we investigate the behavior of the local ROM for patches that are coupled along parameterized trimming interfaces and are expected to behave poorly with a standard global reduced basis. For this purpose we consider a planar setting, where the computational domain is a unit square and is subdivided into two patches coupled along a curved trimming interface.  We employ the projected super-penalty method to enforce interface coupling conditions. The geometric setup and parameterization are depicted in Fig.~\ref{fig:plane_setup} for a coarse geometry, while the interface is shown in red color. The geometric parameter $\mu \in \mathcal{P}=[0.25,0.75]$ defines the position of the control point $P(\mu)$ that prescribes the curvature of the trimming interface. The material parameters are the Young's modulus $E=10^6$ Pa, the Poisson ratio $\nu=0.3$, and the thickness $t=0.022$ m. The non-conforming discretization is generated by shifting the internal knots of the original knot vector at the trimming interface by a factor $\frac{\sqrt[]{2}}{100}$. Note that since both patches are trimmed, the internal knots of the coupling curve are neglected for the construction of the interface knot vector $\Xi_j$ following the discussion in Remark~\ref{rmk:coradello2021}. The geometry is discretized with quadratic $C^1$-continuous B-splines and the dimension of the non-trimmed space is $\mathcal{N}_{h,0}=1944$.  The applied boundary conditions and loading are adopted from a manufactured smooth solution given in \cite{Coradello2021} such that $(u_x,u_y,u_z) = (0,0, \sin(\pi x) \sin(\pi y))$.

\begin{figure}[!h]
	\centering
	\includegraphics[width=0.5\textwidth]{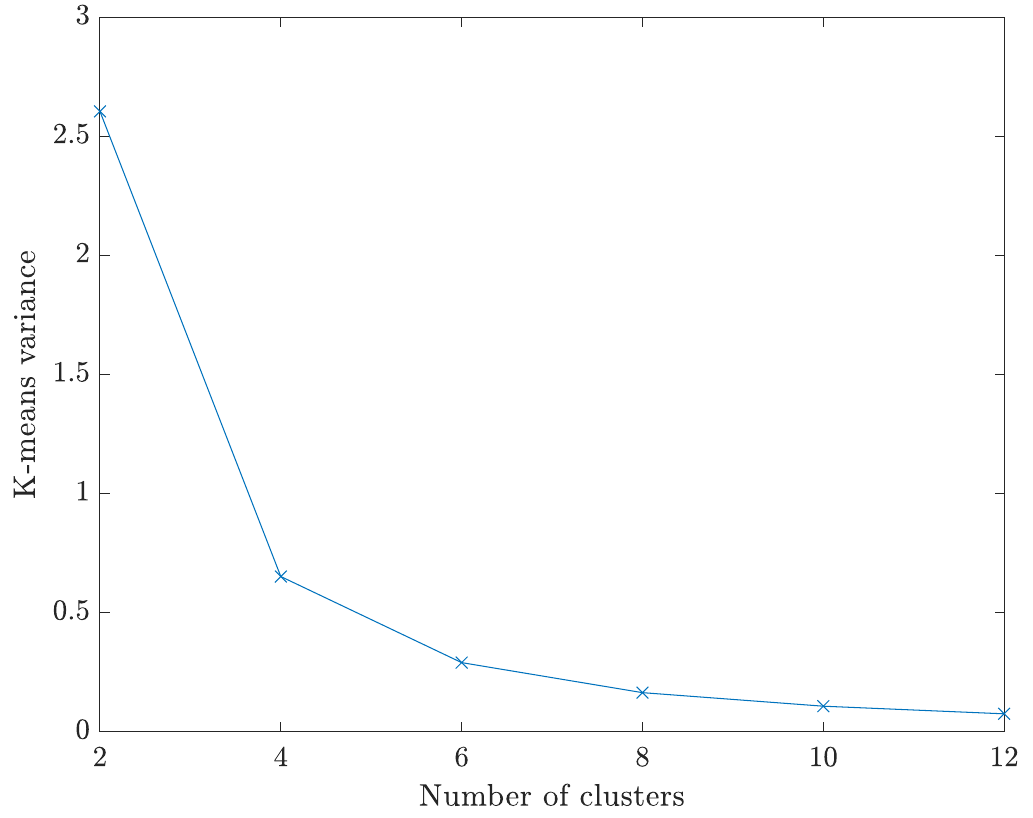} \\
	\caption{Example 7.2.2: K-means variance versus number of clusters $N_c$.}\label{fig:trimmed_plane}
\end{figure}

\begin{figure}[!h]
	\begin{subfigure}[b]{0.49\textwidth}
		\centering
		\includegraphics[width=\textwidth]{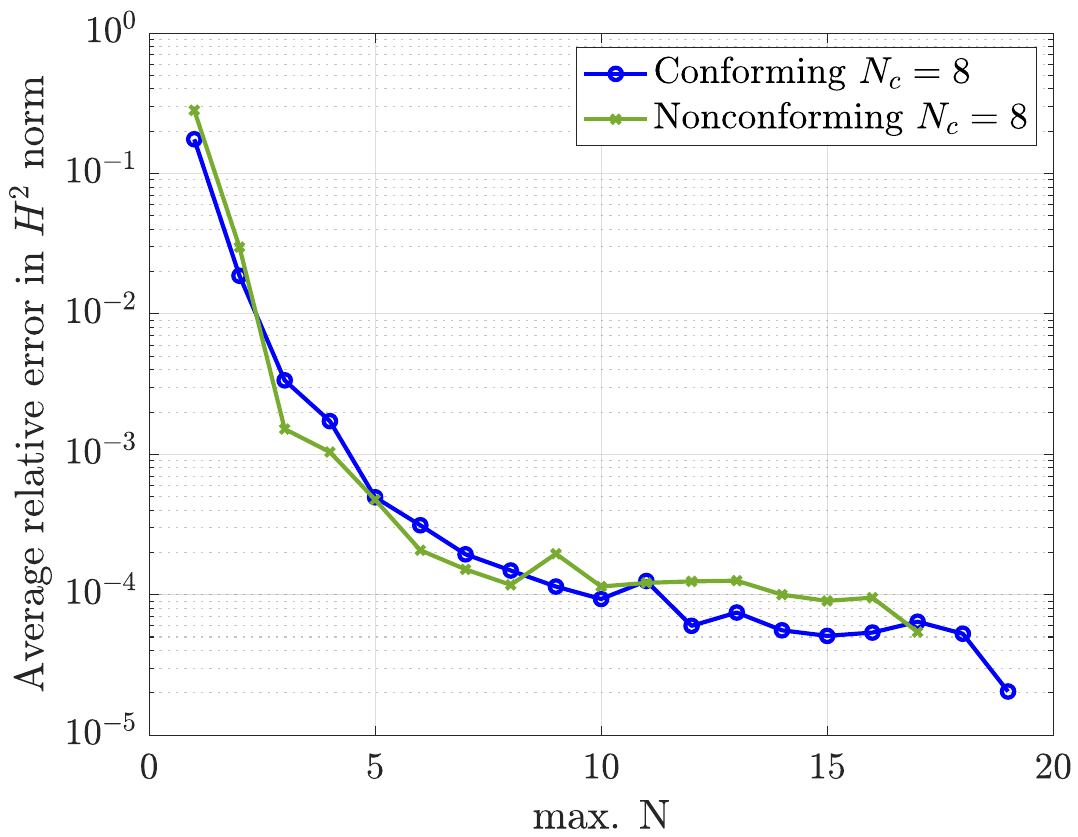}
		\caption{$N_c=8$}
		\label{fig:error_3}
	\end{subfigure}
	\begin{subfigure}[b]{0.49\textwidth}
		\centering
		\includegraphics[width=\textwidth]{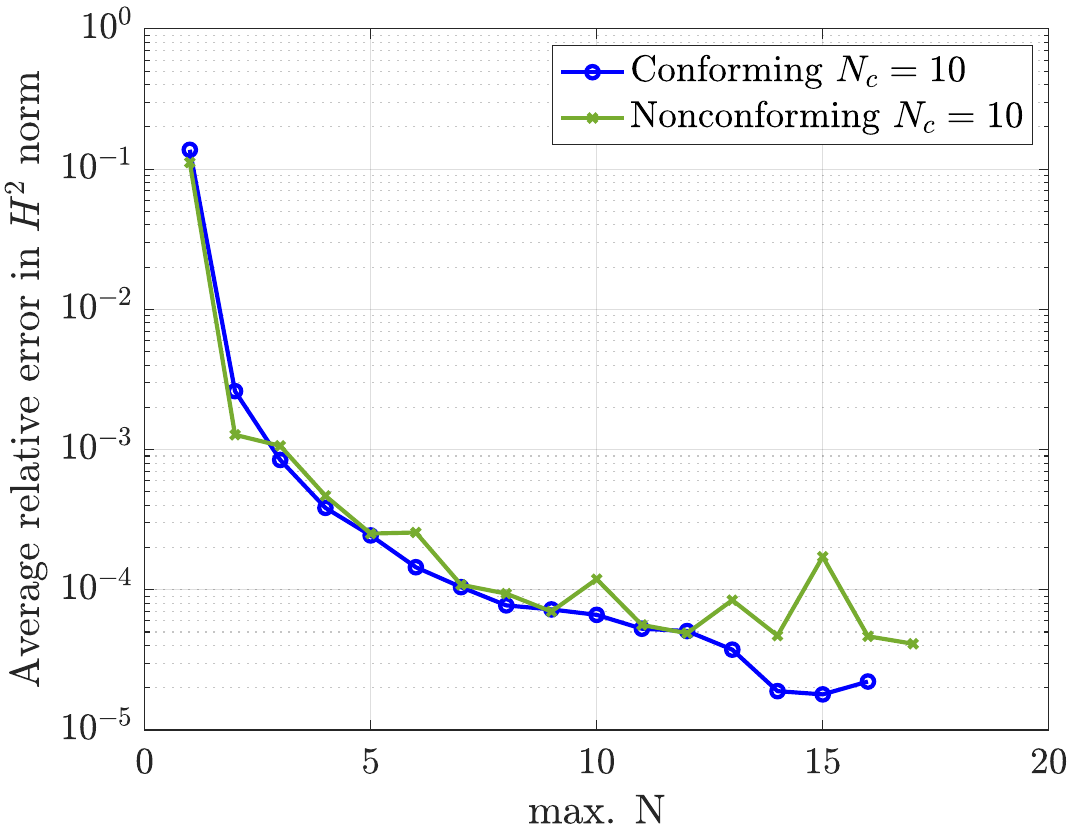}
		\caption{$N_c=10$}
		\label{fig:error_4}
	\end{subfigure}
	\caption{Example 7.2.2: Relative error in $H^2$ norm vs.\ maximum number of reduced basis functions $N$ over all clusters.}
	\label{fig:trimmed_error}
\end{figure}

Now let us construct local ROMs for both the conforming and non-conforming case. For this purpose, we employ a training sample of dimension $N_s = 500$ that we obtain by Latin Hypercube sampling. First, we investigate the k-means variance in Equation~\eqref{eq51a} in order to choose the number of clusters. As it can be seen in Fig.~\ref{fig:trimmed_plane}, the k-means variance does not decrease significantly after $N_c=10$ clusters. To perform the error analysis, we consider a testing sample of dimension $N_t=30$, which is obtained by a uniform random distribution. The relative error in $H^2$ norm is depicted in Fig.~\ref{fig:trimmed_error} for $N_c=8,10$ and the results of both discretizations are compared to each other. As expected, the size of the reduced basis decreases for increasing number of clusters. We observe that the non-conformity slightly affects the maximum number of reduced basis functions $N$, while the error is of the same magnitude in all cases. It can be concluded that the ROM framework is suitable for multi-patch geometries coupled at trimmed interfaces with both conforming and non-conforming discretizations.

\subsection{Joint of intersecting tubes}
The last example aims to demonstrate the capabilities of the presented framework for complex geometries. For this purpose we consider the geometry of three intersecting tubes (see Fig.~\ref{fig:joint_geometry}) that represent a generic configuration for, e.g., joints in steel support truss structures (see also \cite{Guo2017,Prosperio2022} for a similar variant). We remark that optimizing such large-scale structures of industrial relevance is a challenging task, where the shape parameters of each joint element may differ and the ROM can be reused for several online evaluations. The optimization of such problems requires further measures for efficiency (e.g. domain decomposition strategies), which is out of the scope of the present work. For this particular geometry, both sides are trimmed at all interfaces, which hinders the use of the projected super-penalty approach. To this end, we employ the Nitsche's method to enforce interface coupling conditions. The material parameters are the Young's modulus $E= 3 \cdot 10^6$, the Poisson's ratio $\nu=0.3$, and the thickness $t=0.2$. The radius of the main tube is $R_1=1$. Homogeneous Dirichlet boundary conditions are applied on the top and bottom sections of the main tube, while periodic boundary conditions are applied at the cylinders' closure. A vertical load $f_z = 10$ is applied on all tubes. The geometry is discretized with cubic $C^2$-continuous B-splines for the analysis and results in $
\mathcal{N}_{h,0}=2673$ degrees of freedom. In the following we will demonstrate the capabilities of the ROM for two test cases of geometrical parameterization.

\subsubsection{Parameterization of radius}
\begin{figure}[!h]
	\begin{subfigure}[b]{0.3\textwidth}
		\centering
		\includegraphics[width=\textwidth]{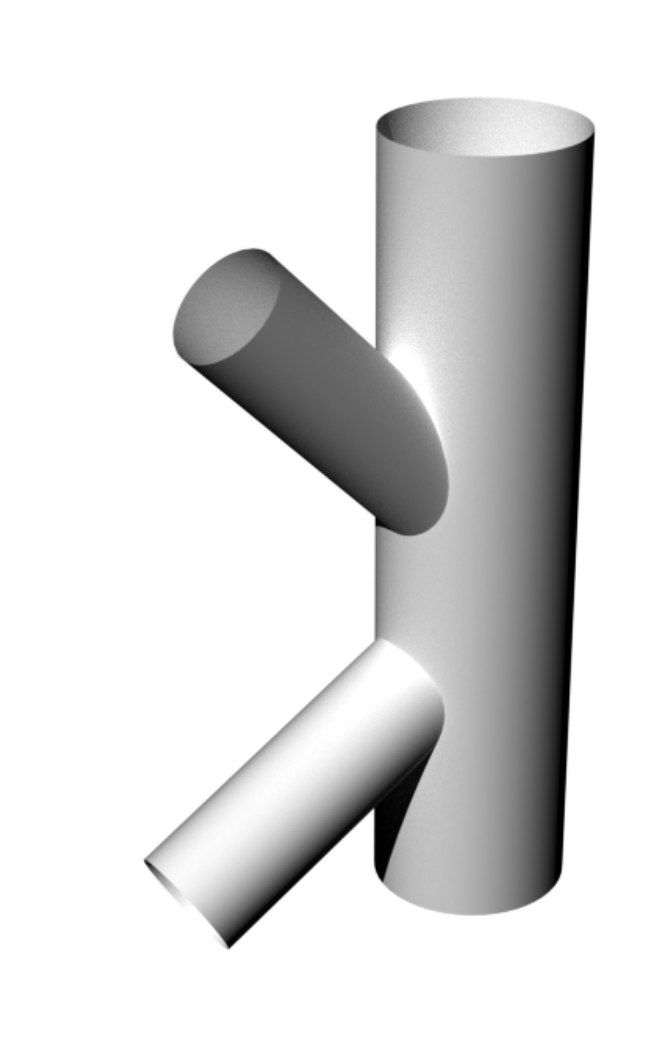}
		\caption{$\mu=0.6$}
		\label{fig:joint_geometry_1}
	\end{subfigure}
	\begin{subfigure}[b]{0.3\textwidth}
		\centering
		\includegraphics[width=\textwidth]{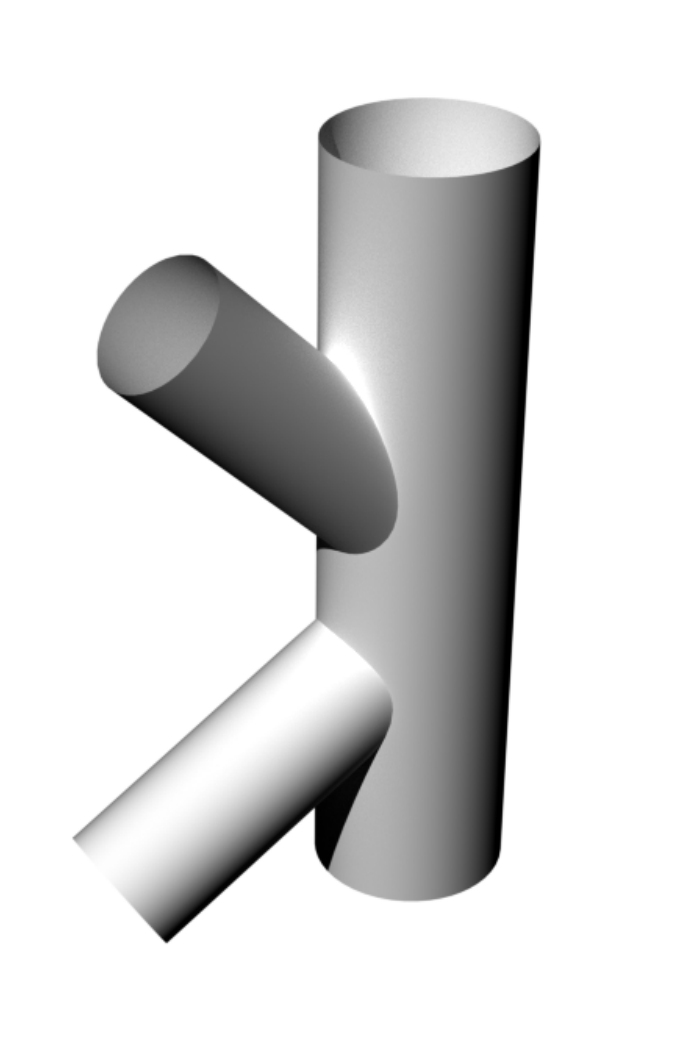}
		\caption{$\mu=0.7$}
		\label{fig:joint_geometry_2}
	\end{subfigure}
	\begin{subfigure}[b]{0.3\textwidth}
		\centering
		\includegraphics[width=\textwidth]{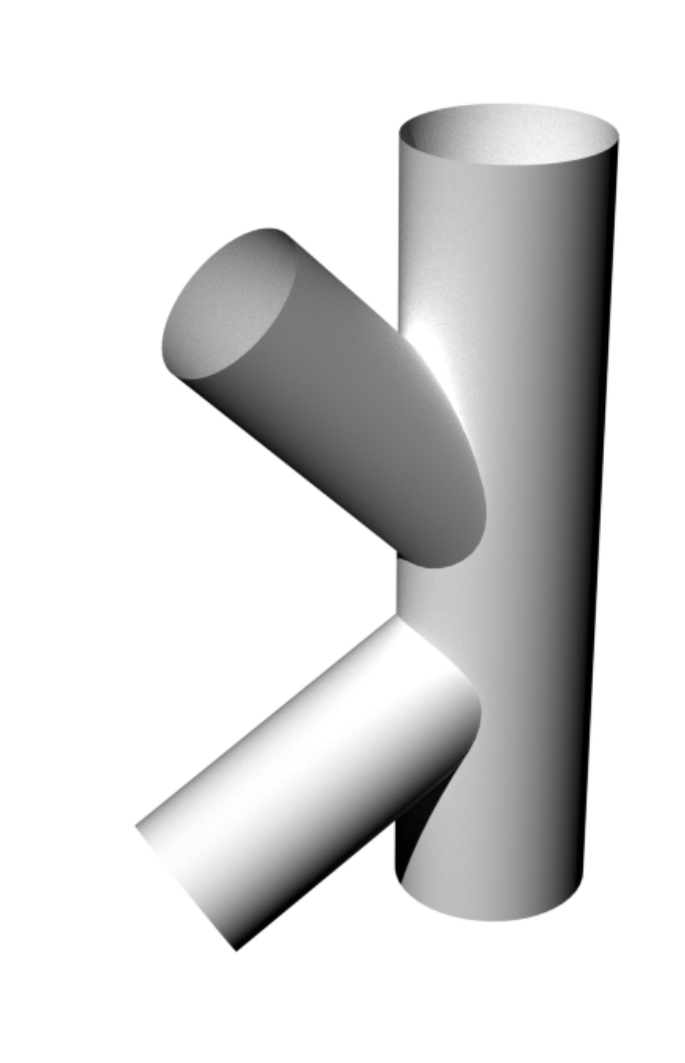}
		\caption{$\mu=0.8$}
		\label{fig:joint_geometry_3}
	\end{subfigure}
	\caption{Example 7.3.1: Geometrical parameterization for different values of $\mu$ (radius of skewed cylinders) for the joint of intersecting tubes.}
	\label{fig:joint_geometry}
\end{figure}

\begin{figure}[!h]
	\centering
	\includegraphics[width=0.5\textwidth]{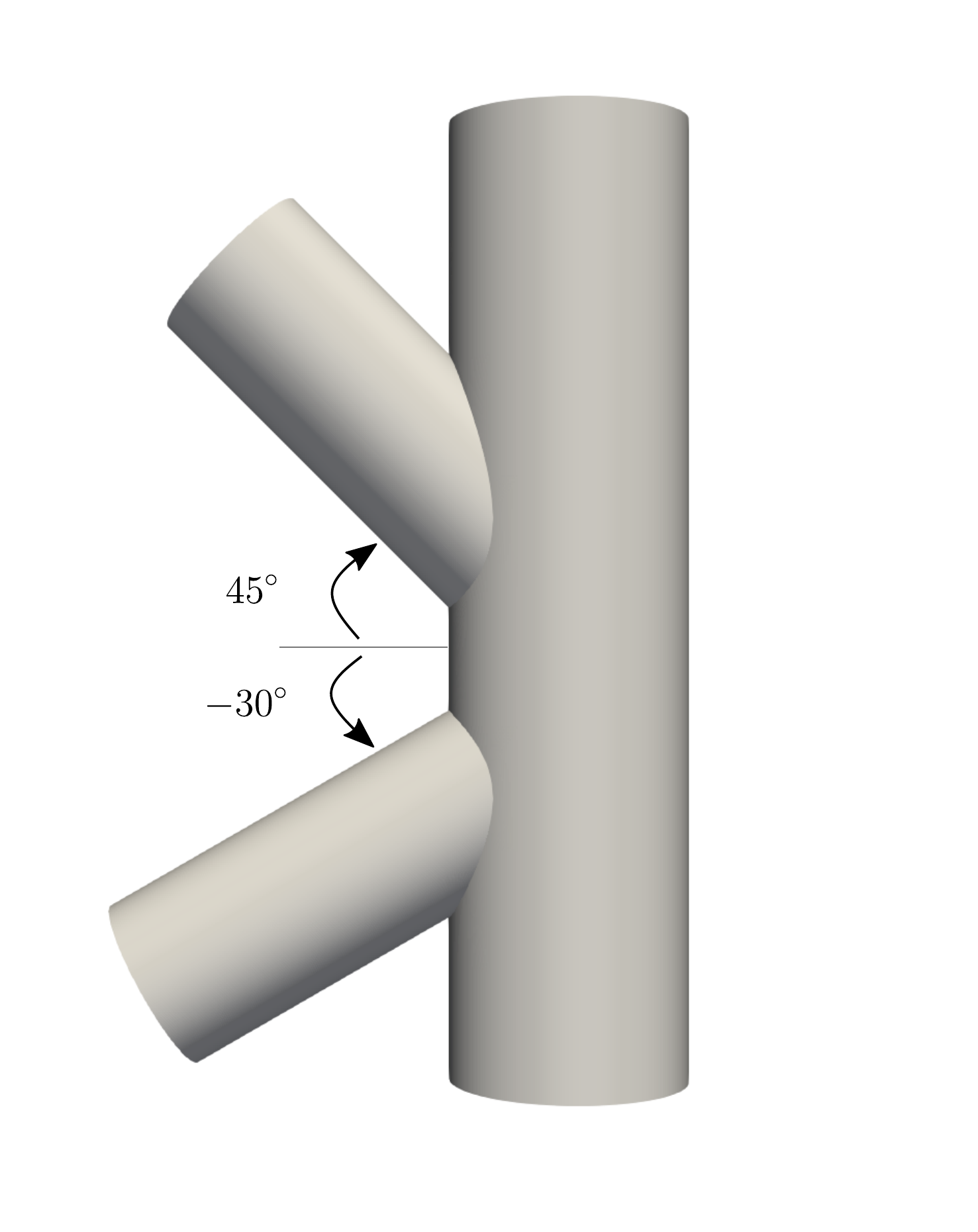} \\
	\caption{Example 7.3.1: Exemplary side view for the joint of intersecting tubes with angle configuration.}
	\label{fig:joint_sideview}
\end{figure}

\begin{figure}[!h]
	\centering
	\includegraphics[width=1.0\textwidth]{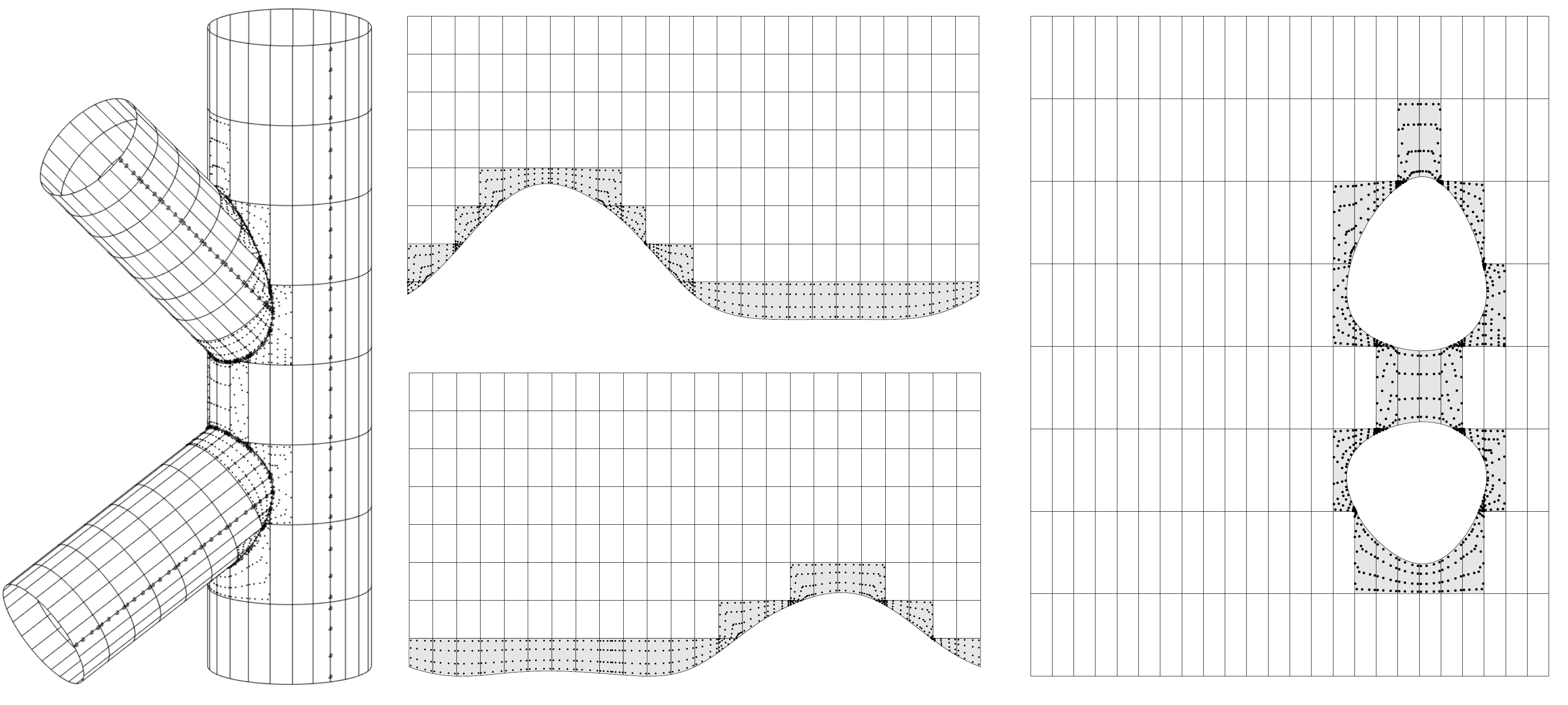} \\
	\caption{Example 7.3.1: Trimming configuration for $\mu=0.7465$. On the left, trimming patches in the physical domain, on the center and right, trimmed parametric domains (trimmed elements are shaded).
 Black dots depict quadrature points in the cut elements, crosses (in physical domain) quadrature points at interfaces.}
 \label{fig:joint_interfaces}
\end{figure}

\begin{figure}[!h]
	\begin{subfigure}[b]{0.49\textwidth}
		\centering
		\includegraphics[width=\textwidth]{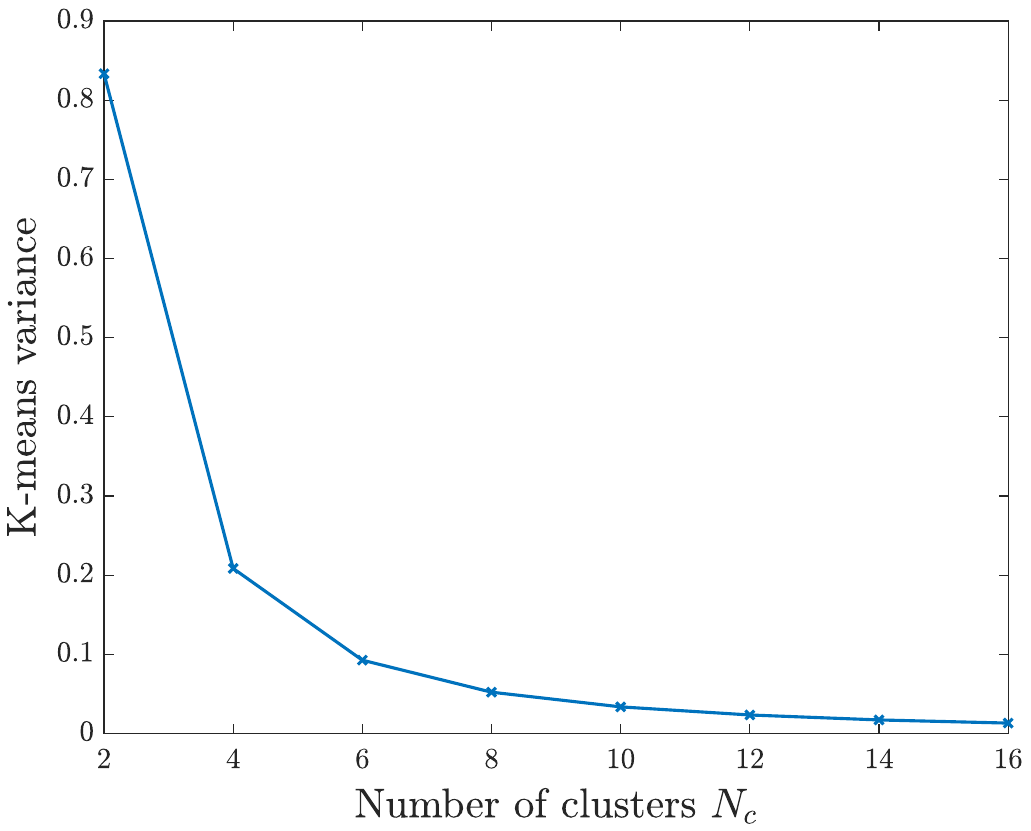}
		\caption{K-means variance}
		\label{fig:joint_variance}
	\end{subfigure}
	\begin{subfigure}[b]{0.49\textwidth}
		\centering
		\includegraphics[width=\textwidth]{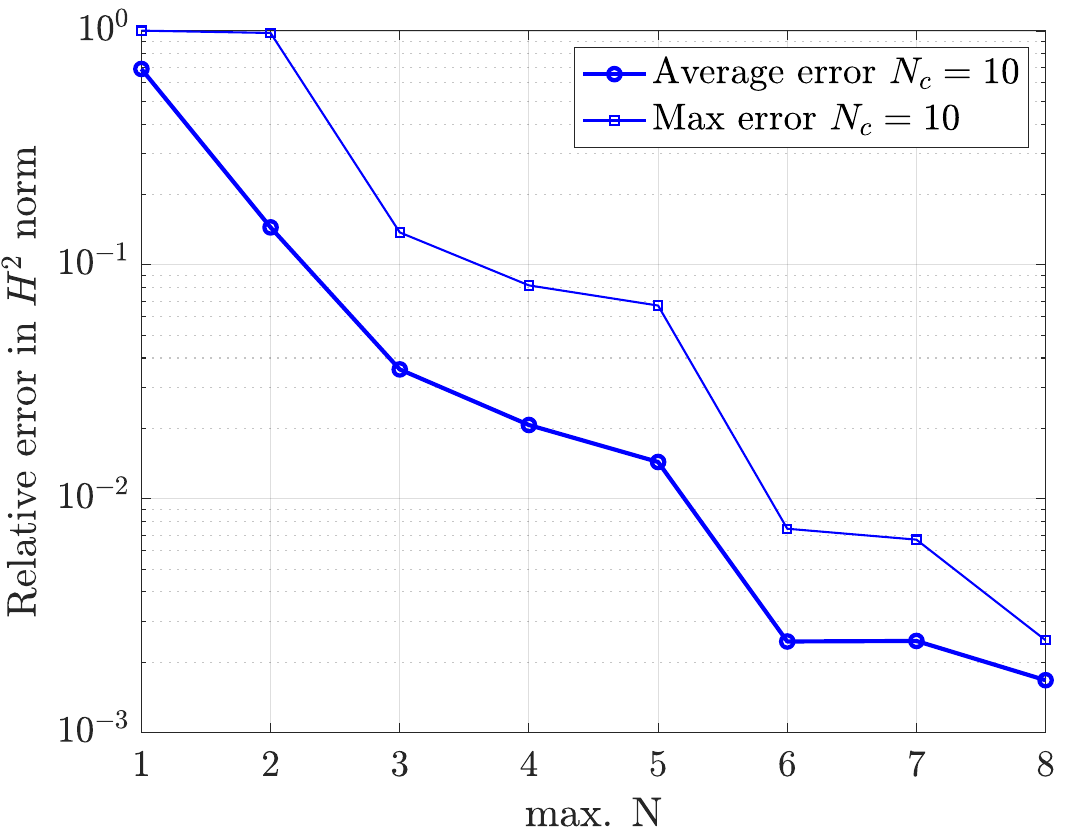}
		\caption{Error analysis}
		\label{fig:joint_error}
	\end{subfigure}
	\caption{Example 7.3.1: K-means variance over number of clusters $N_c$ and relative error in $H^2$ norm vs.\  maximum number of reduced basis functions $N$ over all clusters.}
	\label{fig:joint_results}
\end{figure}

\begin{figure}[!h]
	\centering
	\includegraphics[width=0.8\textwidth]{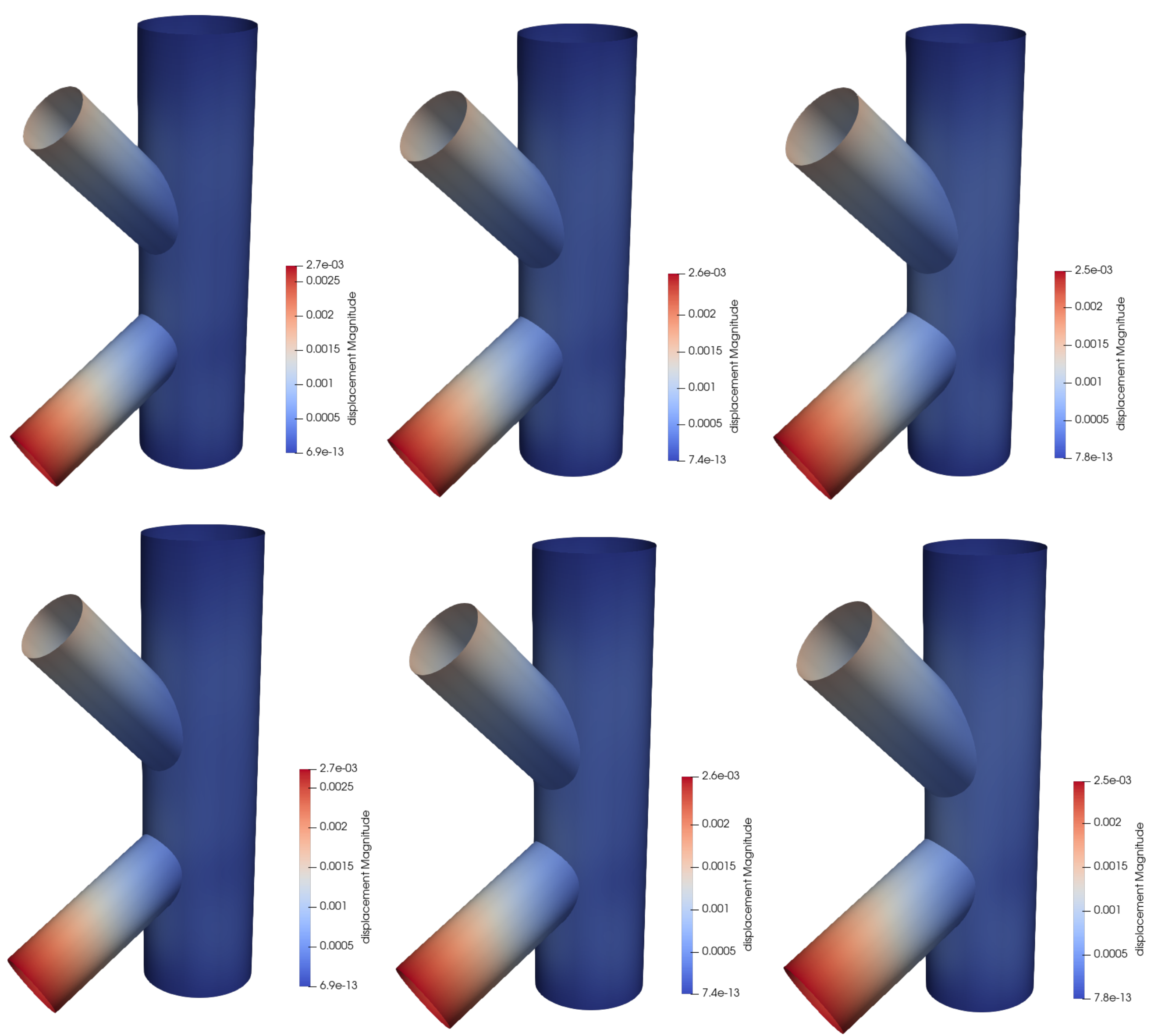} \\
	\caption{Example 7.3.1: Vertical displacement solutions computed with the FOM (top) and ROM (bottom) for three parameter values $\mu= \lbrace0.6033, 0.6789, 0.7465\rbrace$.}\label{fig:joint_plot}
\end{figure}

First, we consider a geometrical parameter $\mu \in \mathcal{P}=[0.6,0.8]$ that represents the radius of the connecting skewed tubes. The geometry setup and parameterization is depicted in Fig.~\ref{fig:joint_geometry}. Note that the angles of the connecting tubes with respect to the horizontal plane are fixed to $45^{\circ}$ and $-30^{\circ}$ for the top and bottom tubes, respectively (see Fig.~\ref{fig:joint_sideview}). In Fig.~\ref{fig:joint_interfaces}, the trimming configuration, including the trimmed parametric domains and quadrature points are shown for a particular case of $\bm\mu$.
To construct the ROM, we employ a training sample of dimension $N_s = 1000$ that we obtain by Latin Hypercube sampling. Note that this refers to the global dimension before clustering, while the snapshots are computed in parallel to speedup the offline phase. The k-means variance is depicted in Fig.~\ref{fig:joint_variance} for increasing number of clusters. We observe that the variance does not change significantly after 10 clusters, thus $N_c=10$ is chosen in what follows. The error analysis is performed using a test sample of dimension $N_t = 30$ obtained by a uniform random distribution. Fig.~\ref{fig:joint_error} shows the error convergence for $N_c=10$ while in Fig.~\ref{fig:joint_plot} the displacement solution obtained with the ROM is compared to the FOM for three  parameter values from the test sample. The main results and computation times are summarized in Table~\ref{tab:joint_times}. 

\begin{table}[!h]
	\caption{Example 7.3.1: Number of basis functions  and computational cost.}\label{tab:joint_times}
        \centering
	\begin{tabular*}{0.65\textwidth}{@{\extracolsep\fill}lc}
		\toprule
		max. $Q_a$    & 25  \\
		max. $Q_f$   & 22 \\
		max. $N$   & 8  \\
		ROM online CPU time [s]   & 0.14 \\
	    ROM offline CPU time [min]   & 79 \\
		FOM time [s]   & 11\\
		\botrule
	\end{tabular*}
\end{table}

\subsubsection{Parameterization of angle}
\begin{figure}[!h]
	\begin{subfigure}[b]{0.49\textwidth}
		\centering
		\includegraphics[width=\textwidth]{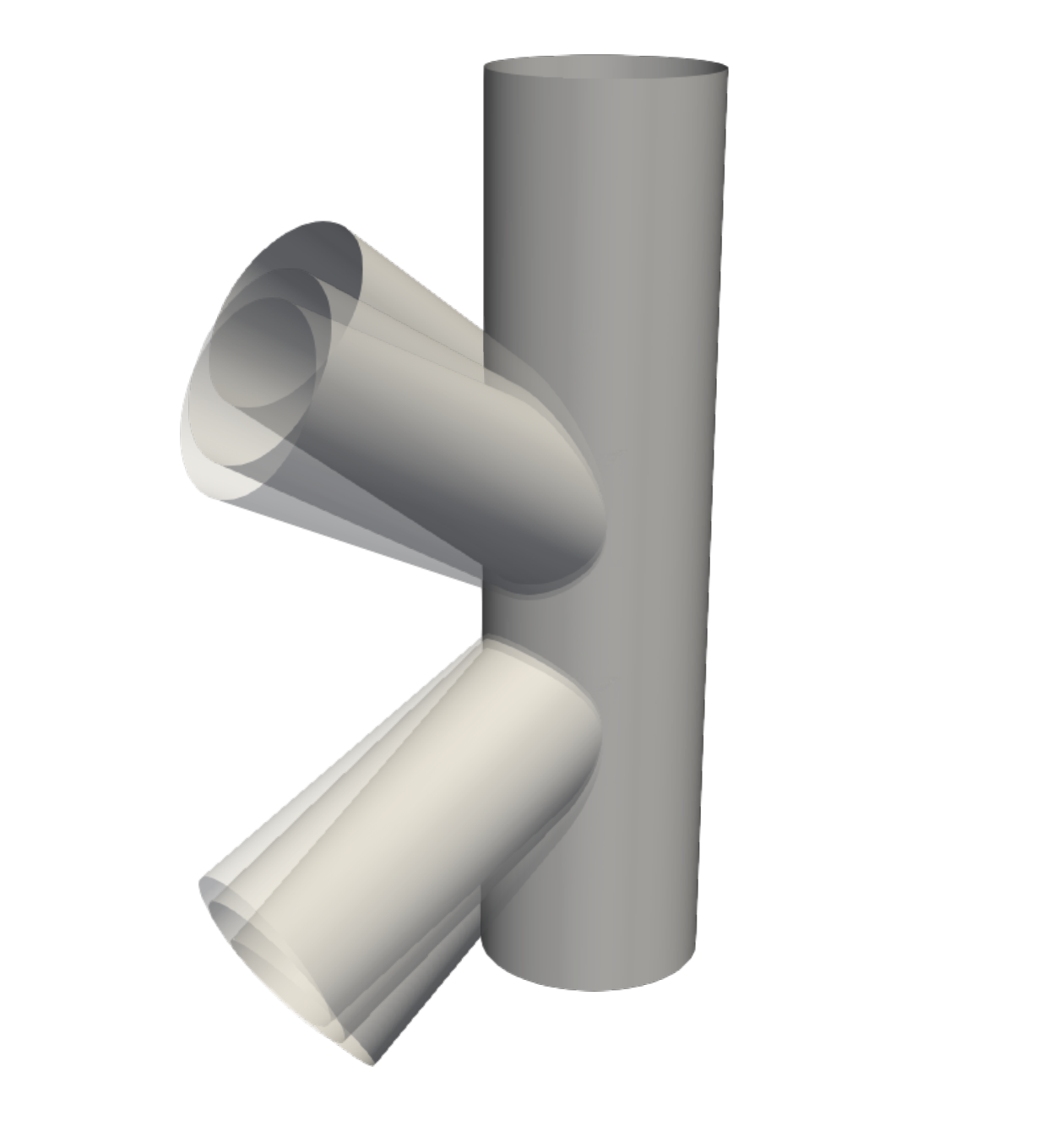}
		\caption{Angle parameterization}
		\label{fig:joint_param}
	\end{subfigure}
	\begin{subfigure}[b]{0.49\textwidth}
		\centering
		\includegraphics[width=\textwidth]{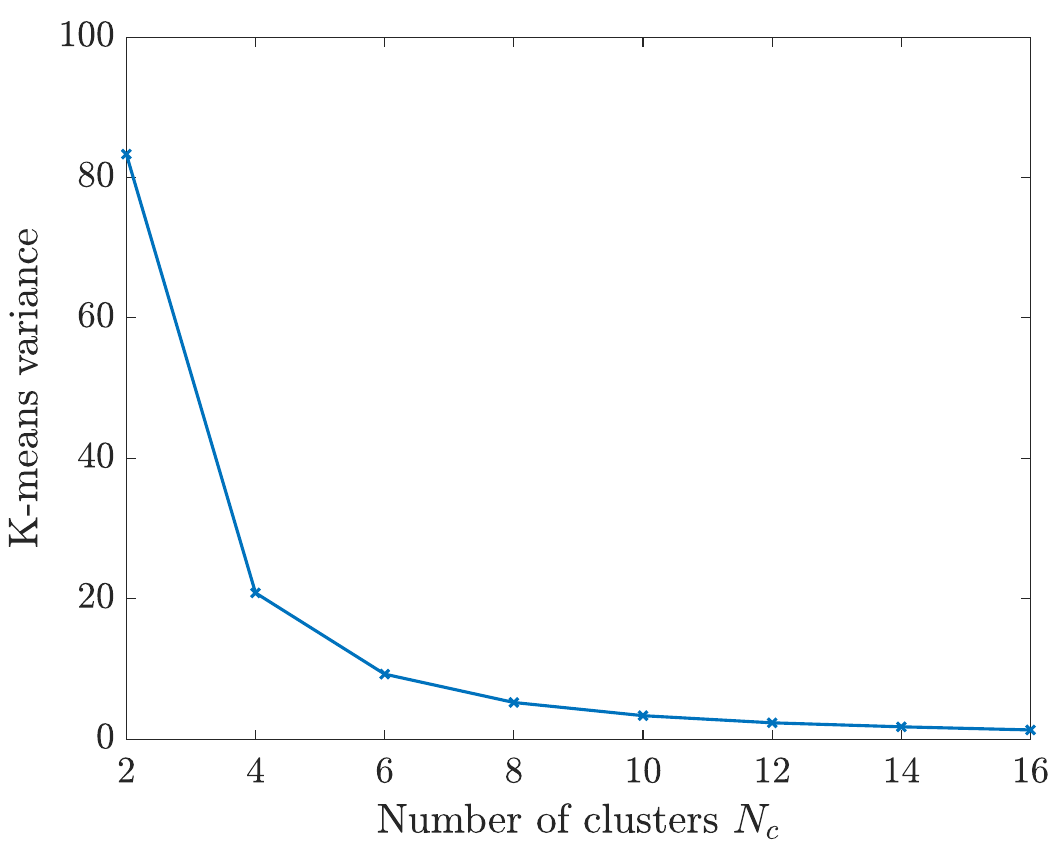}
		\caption{K-means variance}
		\label{fig:joint_variance2}
	\end{subfigure}
	\caption{Example 7.3.2: Geometrical parameterization and k-means variance over number of clusters $N_c$.}
	\label{fig:joint_setup}
\end{figure}

\begin{figure}[!h]
	\begin{subfigure}[b]{0.49\textwidth}
		\centering
		\includegraphics[width=\textwidth]{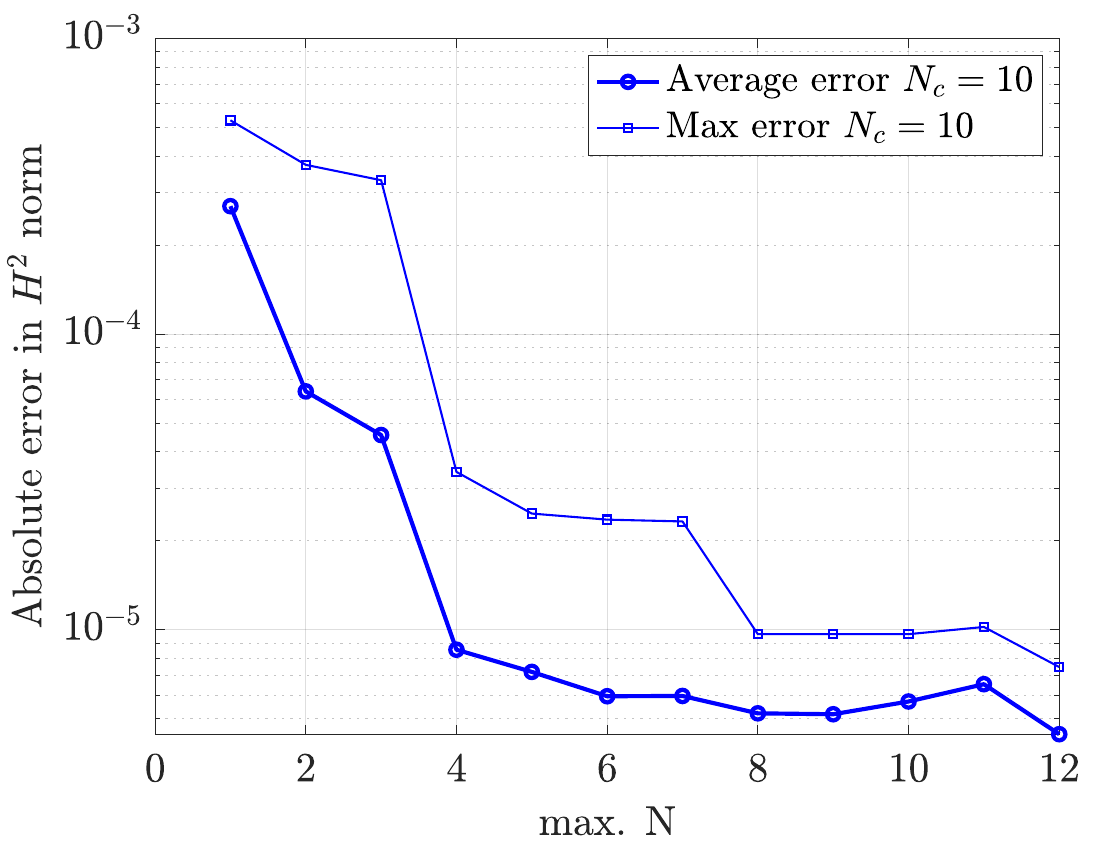}
		\caption{Absolute error}
		\label{fig:joint_absolute}
	\end{subfigure}	
	\begin{subfigure}[b]{0.49\textwidth}
		\centering
		\includegraphics[width=\textwidth]{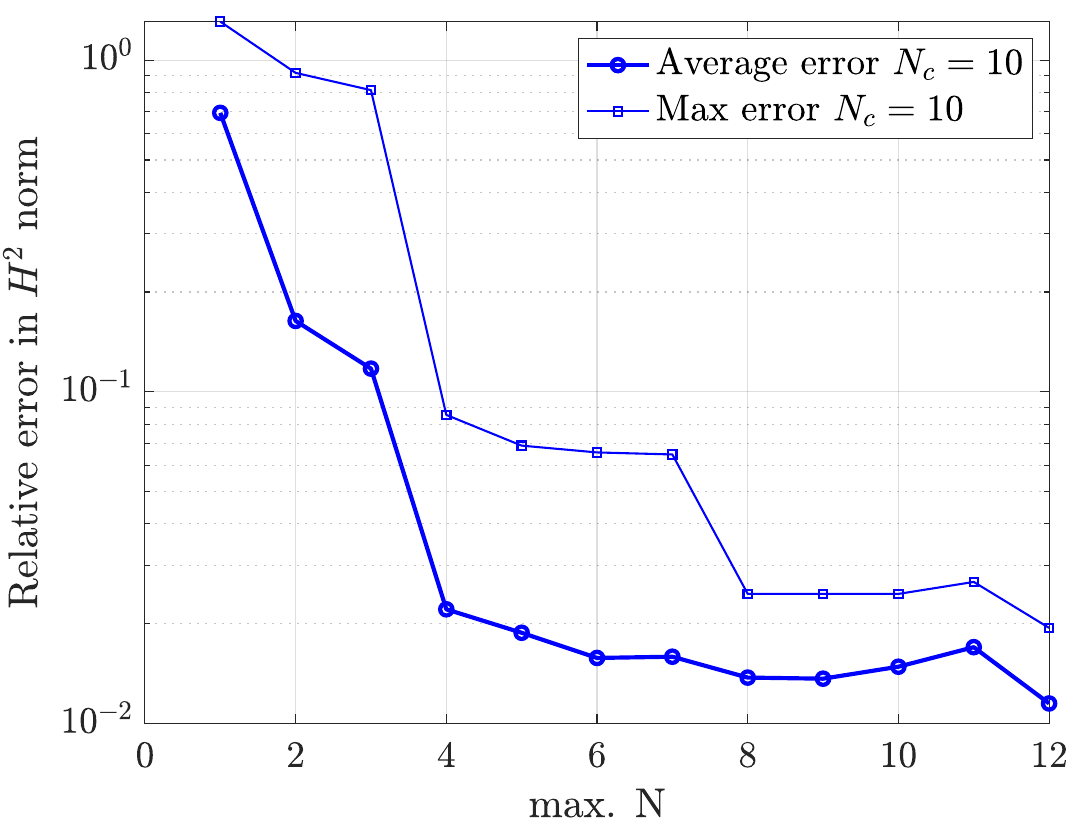}
		\caption{Relative error}
		\label{fig:joint_relative}
	\end{subfigure}
	\caption{Example 7.3.2: Absolute and relative errors in $H^2$ norm vs.\ maximum number of reduced basis functions $N$ over all clusters.}
	\label{fig:joint_error2}
\end{figure}

\begin{figure}[!h]
	\centering
	\includegraphics[width=0.8\textwidth]{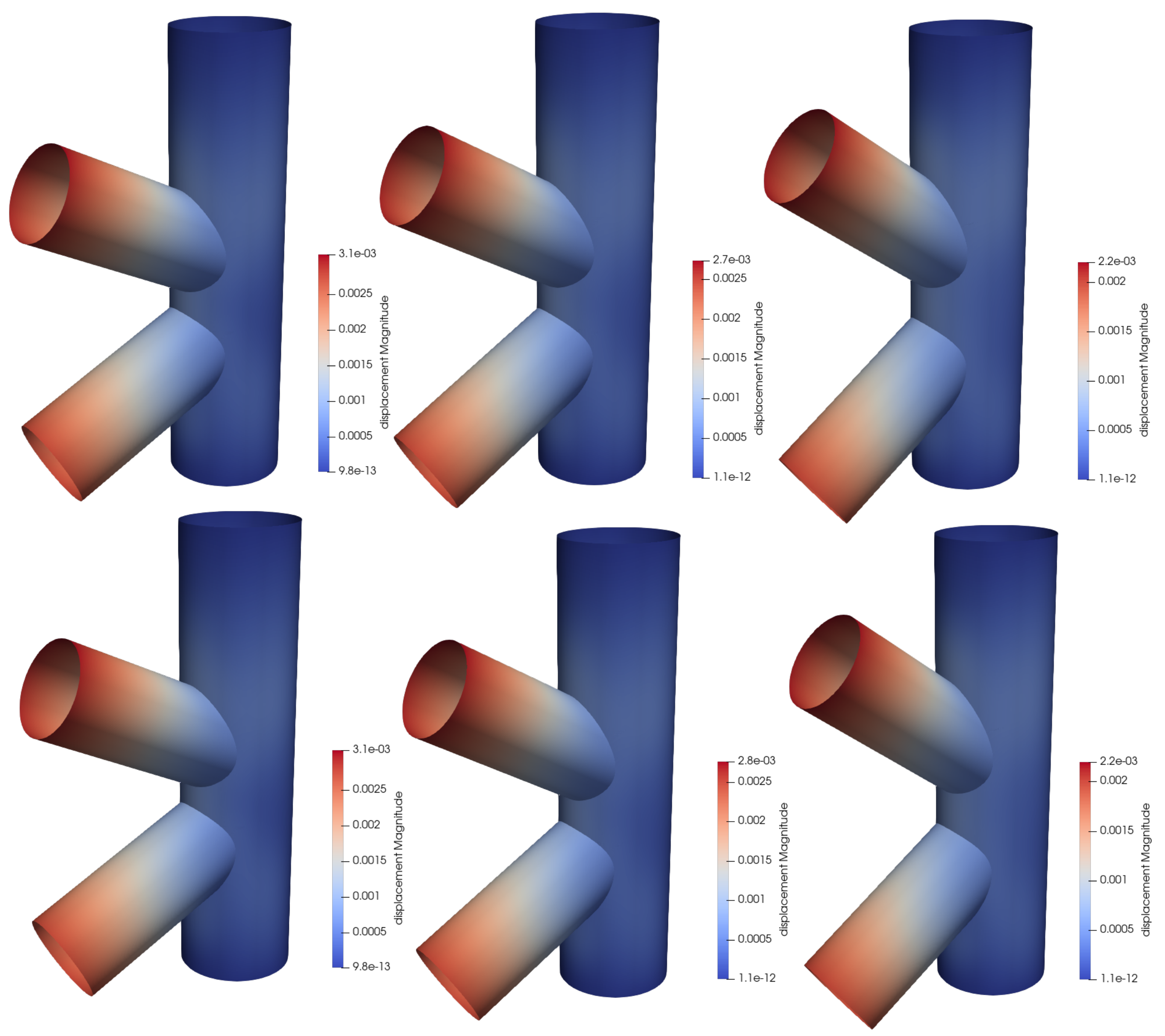} \\
	\caption{Example 7.3.2: Vertical displacement solutions computed with the FOM (top) and ROM (bottom) for three parameter values $\mu= \lbrace25.87^{\circ}, 29.61^{\circ}, 35.78^{\circ} \rbrace$.}\label{fig:joint_plot2}
\end{figure}

\begin{figure}[!h]
	\begin{subfigure}[b]{0.49\textwidth}
		\centering
		\includegraphics[width=\textwidth]{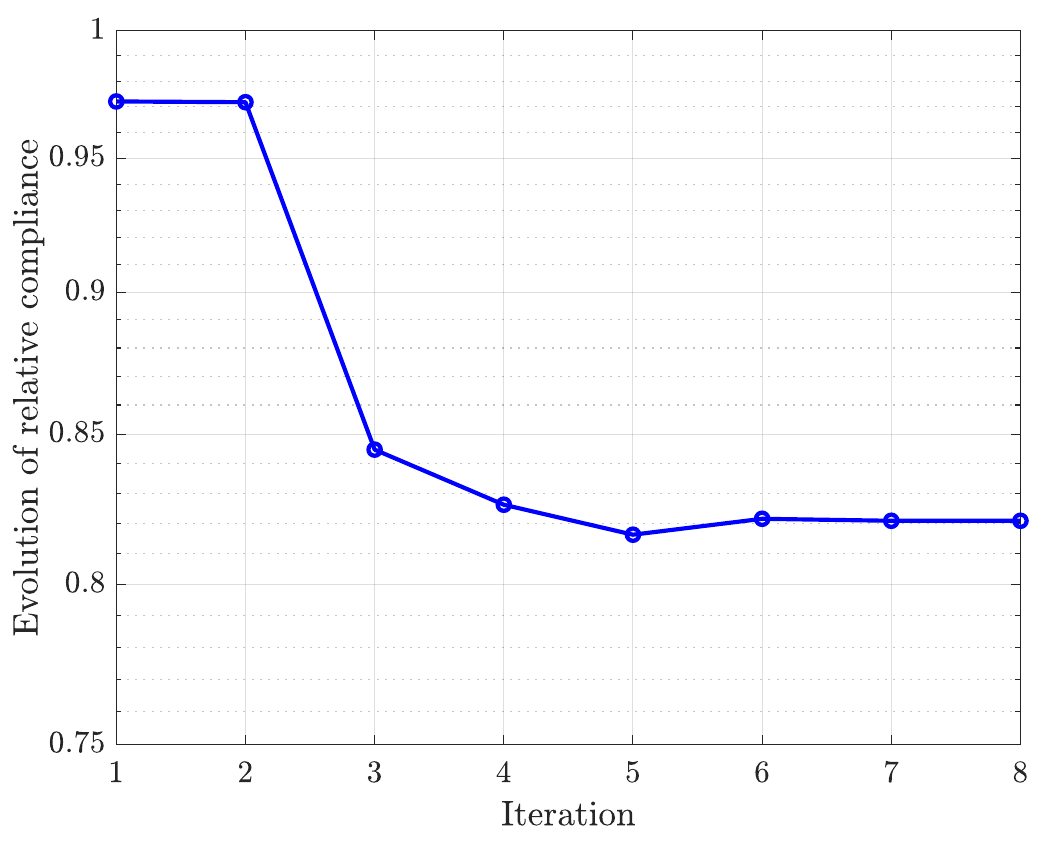}
		\caption{Optimization results}
		\label{fig:joint_angle}
	\end{subfigure}
	\begin{subfigure}[b]{0.4\textwidth}
		\centering
		\includegraphics[width=\textwidth]{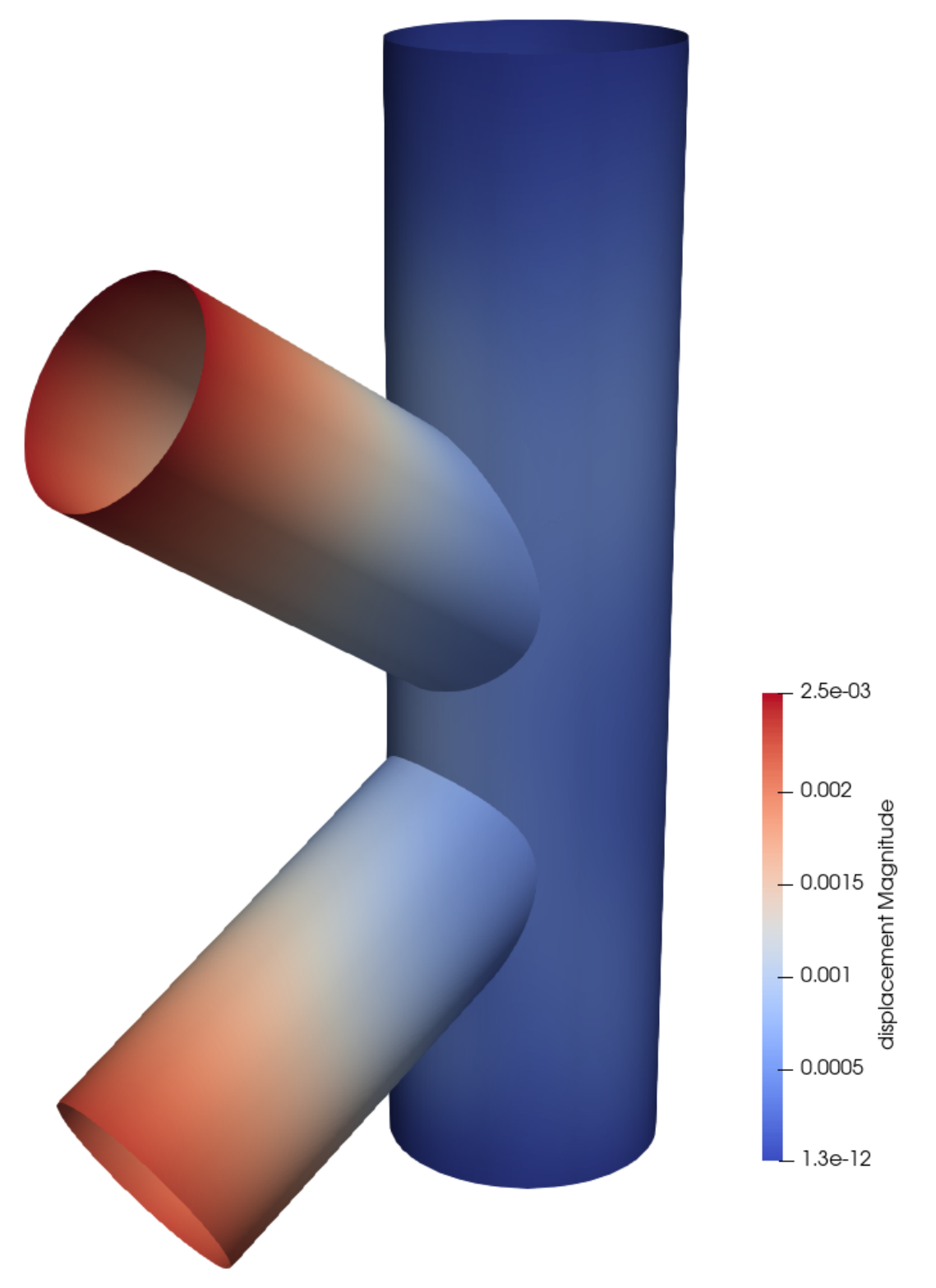}
		\caption{Solution for optimal shape}
		\label{fig:joint_opti}
	\end{subfigure}
	\caption{Example 7.3.2: Optimization results: (a) the evolution of the relative compliance during the optimization for the ROM with exact sensitivities; and (b) the displacement for optimal shape with $\mu=32.22^{\circ}$ using the ROM.}
	\label{fig:joint_optimization}
\end{figure}

Let us now consider a geometrical parameter that represents the angle of the connecting tubes with respect to the horizontal plane. In this case, the radius of the intersecting tubes is fixed as $R_2=0.8$ and the absolute value of the angle varies as $\mu \in \mathcal{P}=[25.7^{\circ},36^{\circ}]$. Note that here the angles of both tubes are symmetric with respect to the horizontal plane. The geometrical variation of the angle is illustrated in Fig.~\ref{fig:joint_param}. Training and test samples are generated similarly to the previous test case for the construction of the ROM and the error analysis, respectively. The number of clusters is chosen as $N_c=10$ based on the k-means variance in Fig.~\ref{fig:joint_variance2}. Fig.~\ref{fig:joint_error2} shows the absolute and relative error in $H^2$ norm with respect to the maximum number of reduced functions over all clusters, while Fig.~\ref{fig:joint_plot2} compares the displacement solutions using the ROM  with the FOM for three parameter values of the test sample. Moreover, we employ the ROM to solve an optimization problem that minimizes the compliance within the parameter space $\mathcal{P}$ without volume constraints. Here, we consider a displacement constraint such that $\text{max}({\bf{u}}(\boldsymbol{\mu)})\le 2.5 \cdot 10^{-3}$. That is, the compliance is minimized such that the maximum displacement does not exceed the prescribed value. We remark that we also employ the ROM for the computation of the displacement constraint to speedup the optimization. Note that the initial shape at the beginning of the optimization corresponds to the minimum angle with $\mu=25.7^{\circ}$. The optimization using exact sensitivities requires 8 iterations and 29 function evaluations.  Fig.~\ref{fig:joint_angle} depicts the evolution of the relative compliance during the optimization.  At the final iteration, the compliance is decreased by $17.9 \% $ compared to the initial configuration.  Note that the optimization with approximate sensitivities using a forward finite difference scheme requires 8 iterations and 21 function evaluations to reach almost identical results. Fig.~\ref{fig:joint_opti} depicts the displacement solution for the optimal shape with $\mu=32.22^{\circ}$ and Table~\ref{tab:joint_times2} summarizes the main results and computation times. It is remarked that the online computation time corresponds to the cluster with the maximum number of basis functions, thus the online cost might differ from one parameter to the other.

\begin{table}[!h]
	\caption{Example 7.3.2: Number of basis functions  and computational cost.}\label{tab:joint_times2}
         \centering
	\begin{tabular*}{0.65\textwidth}{@{\extracolsep\fill}lc}
		\toprule
		max. $Q_a$    & 21  \\
		max. $Q_f$   & 18 \\
		max. $N$   & 12  \\
		Online CPU time [s]   & 0.15 \\
		ROM-based optimization time [s] & 0.174 \\ 
		\botrule
	\end{tabular*}
\end{table}

\section{Conclusion}\label{sec8}
In this work we present a parametric model reduction framework for trimmed, multi-patch isogeometric Kirchhoff-Love shells.
The proposed strategy is suitable for fast simulations on parameterized geometries represented by multiple, non-conforming patches where the trimming interfaces change for different values of the geometric parameters. Following our previous work \cite{Chasapi2022}, efficient ROMs are constructed using a local reduced basis method and hyper-reduction with DEIM. The latter enables an efficient offline/online split with low online cost, which is advantageous for solving parametric optimization problems. 

We have investigated the capabilities of the proposed framework through several numerical experiments.
To this end, we considered trimmed, multi-patch geometries with parameterized interfaces, considering both conforming and non-conforming cases, and that were glued applying super-penalty and Nitsche coupling methods.
We observed a high accuracy of the local ROMs, while the solution evaluation in the online phase was obtained with low computational cost. Moreover, we validated the proposed approach for parametric optimization problems and applied it to a complex geometry.

Overall, the application of the ROM framework to isogeometric Kirchhoff-Love shell analysis of complex geometries and optimization problems is a cost-effective alternative. The application to more complex optimization problems, including higher number of design parameters, and the extension to Reissner-Mindlin shell formulations are future research directions to explore. Regarding the ROM, a further subject of future work is related to error certification using greedy algorithms and tailored a posteriori error estimators.

\section*{Acknowledgments}
The financial support of the Swiss Innovation Agency (Innosuisse) under grant No.\ 46684.1 IP-EE,
of Swiss National Science Foundation through the project No.\ 40B2-0\_187094 (BRIDGE Discovery 2019), and
the European Union Horizon 2020 research and innovation program under grant No.\ 862025 (ADAM$^2$)
is gratefully acknowledged.
We would also like to thank Dr.\ Luca Coradello and  Guiliano Guarino for providing the implementation of multi-patch coupling methods. 

\section*{Declarations}
\subsection*{Competing interests}
The authors have no competing interests to declare that are relevant to the content of this article.

\bibliography{sn-bibliography}

\end{document}